\input amstex
\documentstyle {amsppt}

\pagewidth{32pc} 
\pageheight{45pc} 
\mag=1200
\baselineskip=15 pt

\hfuzz=5pt
\topmatter
\NoRunningHeads 
\title Finite Localities III
\endtitle
\author Andrew Chermak
\endauthor
\affil Kansas State University
\endaffil
\address Manhattan Kansas
\endaddress
\email chermak\@math.ksu.edu
\endemail
\date
October 2016 
\enddate 

\endtopmatter 

\redefine \del{\partial}

\redefine\pce{\preccurlyeq} 

\define\w{\widetilde}

\redefine\norm{\trianglelefteq}

\redefine\bar{\overline}

\redefine\maps{\mapsto}
\redefine\i{^{-1}}

\redefine\l{\lambda}
\redefine\s{\sigma}

\redefine\b{\beta}
\redefine\d{\delta}
\redefine\g{\gamma}

\redefine\r{\rho}

\redefine\G{\Gamma}

\redefine\S{\Sigma}
\redefine\L{\Lambda}

\redefine\<{\langle}
\redefine\>{\rangle}
\redefine\F{\Bbb F}
\redefine\ca{\Cal}

\redefine\D{\Delta}

\redefine\sub{\subseteq} 
\redefine\sup{\supseteq}

\redefine\nset{\emptyset} 

\redefine\1{\bold 1}

\document

\vskip .1in 
\noindent 
{\bf Introduction} 
\vskip .1in 

This paper is the third in a series on finite localities, whose earlier installments are Part I [Ch2] and 
Part II [Ch3]. We continue the convention of referring to results in earlier parts by prefixing a ``I" or 
(now) a ``II" to citations. For example, ``II.2.4" refers to the definition 2.4 of proper locality in [Ch3]. 
Familiarity with the earlier parts is assumed, but we will provide a brief review of some of the core 
material here in section 1. 

This Part III relies more heavily on the language of fusion systems than do Parts I and II. In particular, 
if $\ca F$ is the fusion system of a proper locality $(\ca L,\D,S)$, then the set $\ca F^q$ of 
$\ca F$-quasicentric subgroups of $S$ and, more importantly, the set $\ca F^s$ of $\ca F$-subcentric 
subgroups of $S$ will play a key role here in understanding the structure of $\ca L$. (These 
collections, along with the set $\ca F^c$ of $\ca F$-centric subgroups, and the set $\ca F^{cr}$ of 
$\ca F$-centric radical subgroups of $S$ were defined in II.1.8.)

Throughout this paper 
$(\ca L,\D,S)$ will be a proper locality on $\ca F$. Thus, $(\ca L,\D,S)$ is a 
locality, $\ca F=\ca F_S(\ca L)$ is its fusion system, and it is assumed that $\ca F^{cr}\sub\D$, and 
that the normalizer subgroups $N_{\ca L}(P)$ for $P\in\D$ are of characteristic $p$ (whence $\D\sub\ca F^s$, 
by II.2.8). The main theorems in Part II showed how one may alter the set $\D$ to ``restrict" or 
``expand" $\ca L$, while 
preserving $\ca F$, and while preserving also the structure of the poset of partial normal subgroups of 
$\ca L$. These theorems provide the flexibility whereby $\D$ may be chosen to be whatever 
$\ca F$-closed collection $\D$ of subgroups of $S$ is most convenient for a particular analysis, subject 
only to the requirement $\ca F^{cr}\sub\D\sub\ca F^s$. The principal aim of this paper is to show that there 
is such an $\ca F$-closed collection of subgroups of $S$, 
to be denoted $\d(\ca F)$, which is in many ways preferential - and for at least the following two 
inter-connected reasons. 
\roster 

\item "{(1)}" If $\D=\d(\ca F)$ then each partial normal subgroup $\ca N\norm\ca L$ is itself a proper 
locality. 

\item "{(2)}" If $\D=\d(\ca F)$ then, for each partial normal subgroup $\ca N\norm\ca L$, the set 
$C_{\ca L}(\ca N)$ of all $g\in\ca L$ such that $x^g$ is defined and is equal to $x$, for all $x\in\ca N$, 
is a partial normal subgroup of $\ca L$. 

\endroster 

For any locality $(\ca L,\D,S)$ there is an action of $S$ on $\ca L$ by conjugation, and we may therefore 
speak of $C_S(X)$ for any non-empty subset $X$ of $\ca L$. We shall say that a partial normal subgroup 
$\ca M\norm\ca L$ is {\it large} if $C_S(\ca M)\leq\ca M$. One way 
to define the set $\d(\ca F)$ is to begin by expanding $\ca L$ to a proper locality 
$(\ca L^s,\ca F^s,S)$ on $\ca F$ whose set of objects is as large as possible. That is, we may employ Theorem 
II.A1 in order to obtain the unique (up to isomorphism) proper locality on $\ca F$ whose set of objects is 
the set $\ca F^s$ of $\ca F$-subcentric subgroups of $S$. Then define $F^*(\ca L^s)$ to be the 
intersection of the set of all large partial normal subgroups $\ca M$ of $\ca L^s$ containing 
$O_p(\ca F)$.   
$$ 
F^*(\ca L^s)=\bigcap\{\ca M\norm\ca L^s\mid O_p(\ca F)C_S(\ca M)\leq\ca M\}. 
$$ 
Now define $\d(\ca F)$ to be the set of all subgroups $P$ of $S$ such that 
$P\cap F^*(\ca L^s)\in\ca F^s$. We show in lemma 6.10 that $\ca F^{cr}\sub\d(\ca F)\sub\ca F^s$ 
and that $\d(\ca F)$ is $\ca F$-closed. Thus $\d(\ca F)$ is a permissible choice 
for $\D$; and we say that the proper locality $(\ca L,\D,S)$ on $\ca F$ is {\it regular} if 
$\D=\d(\ca F)$. 

\proclaim {Theorem C} Let $(\ca L,\D,S)$ be a regular locality on $\ca F$, and let $\ca N\norm\ca L$ 
be a partial normal subgroup of $\ca L$. Set $T=S\cap\ca N$ and let $\ca E=\ca F_T(\ca N)$ be the 
fusion system on $T$ generated by the conjugation maps $c_g:V\to T$ with $g\in\ca N$ and with $V\leq T$. 
Define $C_{\ca L}(\ca N)$ to be the set of all $g\in\ca L$ such that, for all $x\in\ca N$, $x^g$ is defined 
and is equal to $x$. Then the following hold. 
\roster 

\item "{(1)}" $\d(\ca E)=\{V\leq T\mid C_S(\ca N)V\in\d(\ca F)\}$, and $\ca N$ is 
itself a regular locality on $\ca E$. 

\item "{(2)}" $C_{\ca L}(\ca N)\norm\ca L$, and $C_{\ca L}(\ca N)=O^p(C_{\ca L}(T))C_S(\ca N)$ 
(cf. section II.7). Moreover, for all $g\in C_{\ca L}(\ca N)$ and all $x\in\ca N$ the product 
$[x,y]:=x\i g\i xg$ is defined and is equal to $\1$. 

\endroster 
\endproclaim

For any proper locality $(\ca L,\D,S)$ on $\ca F$, define $F^*(\ca L)$ to be the partial normal 
subgroup $\ca L\cap F^*(\ca L^s)$ of $\ca L$ corresponding to the partial normal subgroup 
$F^*(\ca L^s)$ of $\ca L^s$ via Theorem II.A2. A partial subgroup $\ca K$ of $\ca L$ is 
{\it subnormal} in $\ca L$ if there is a sequence $(\ca N_0,\cdots,\ca N_m)$ of partial subgroups of 
$\ca L$ with $\ca K=\ca N_0\norm\cdots\norm\ca N_m=\ca L$. The partial subnormal subgroup 
$\ca K$ of $\ca L$ is a {\it component} of $\ca L$ if $\ca K=O^p(\ca K)$ and $\ca K/Z(\ca K)$ is simple 
(where $Z(\ca K)=C_S(\ca K)\cap\ca K))$ ).

For subsets $X$ and $Y$ of a partial group, write $[X,Y]=\1$ to indicate that $[x,y]$ is defined and is 
equal to $\1$ for all $x\in X$ and all $y\in Y$.

\proclaim {Theorem D} Let $(\ca L,\D,S)$ be a regular locality on $\ca F$ and let 
$(\ca K_1,\cdots,\ca K_n)$ be a non-redundant list of all of the components of $\ca L$. Let 
$E(\ca L)$ be the smallest partial subgroup of $\ca L$ containing all of the components $\ca K_i$. 
Then $E(\ca L)\norm\ca L$, and 
$$
F^*(\ca L)=O_p(\ca F)E(\ca L)=O_p(\ca F)\ca K_1\cdots\ca K_n, 
$$ 
where $[O_p(\ca F),E(\ca L)]=\1$, and where $[\ca K_i,\ca K_j]=\1$ for all $i$ and $j$ with $i\neq j$. 
Further, we have $C_S(F^*(\ca L))\leq F^*(\ca L)$, and $\ca L=F^*(\ca L)H$, 
where $H:=N_{\ca L}(S\cap F^*(\ca L))$ is a subgroup of $\ca L$ which acts on $\ca L$ 
by everywhere-defined conjugation.  
\endproclaim 

The last of our main theorems concerns the category of regular localities as a full subcategory of 
the category of partial groups.  

\proclaim {Theorem E} Let $\Bbb L_\d$ be the full sub-category of the category of partial groups, whose 
class of objects is the class of regular localities. Let $(\ca L,\D,S)$ be a 
regular locality on $\ca F$.    
\roster 

\item "{(a)}" If $\ca N\norm\ca L$ is a partial normal subgroup of $\ca L$ then the inclusion map 
$\ca N\to\ca L$ is a $\Bbb L_\d$-homomorphism. 

\item "{(b)}" If $X\leq S$ is fully normalized in $\ca F$, and $X\leq F^*(\ca L)$, then there is a unique 
regular locality $\ca L_X$ on $N_{\ca F}(X)$ such that $\ca L_X$ is a subset of $\ca L$, and such that the 
inclusion map $\ca L_X\to\ca L$ is an $\Bbb L_\d$-homomorphism. 

\item "{(c)}" If $Z\leq C_S(O^p(\ca L))$ is a normal subgroup of $\ca L$ then the quotient locality 
$\ca L/Z$ is regular, and the quotient map $\ca L\to\ca L/Z$ is an $\Bbb L_\d$-homomorphism. 

\endroster 
\endproclaim 

It is a defect of the category $\Bbb L_\d$ that, aside from some special cases (such as 
the one given by point (c) of Theorem 
E), homomorphic images of regular localities need not be regular - or even proper. One way in which to address 
this defect (but which will not be pursued here) is as follows. Define a locality $(\ca L,\D,\S)$ to be 
{\it semi-regular} if there exists a regular locality $\w{\ca L}$ and a projection (cf. I.4.4) 
$\r:\w{\ca L}\to\ca L$. Composites of projections are projections, so the category $\Bbb L_\s$ of partial 
groups whose objects are semi-regular localities is closed with respect to projections. Further, the 
Correspondence Theorem I.4.7 shows that partial normal subgroups of semi-regular localities are 
themselves semi-regular localities. What is then needed is a version of point (b) in Theorem E, in order 
to obtain the beginnings of a satisfactory category.

\vskip .2in 
\noindent 
{\bf Section 1: Localizable pairs} 
\vskip .1in 

We assume that the reader has become comfortable with the basic notions introduced in the first 
section of Part I; specifically partial groups, partial subgroups, partial normal subgroups, and homomorphisms 
of partial groups. These notions will not be reviewed here. 

From section 2 of Part I:A partial group $\ca L$ is {\it objective} if there is a set $\D$ of subgroups of 
$\ca L$ (a set of 
{\it objects}) which defines the domain $\bold D$ of the product $\Pi$ in $\ca L$, and which satisfies a 
pair of closure conditions as expressed in definition II.2.1. A (finite) {\it locality} is a finite, 
objective partial group $(\ca L,\D)$ such that $\D$ is a set of subgroups of some $S\in\D$,  
$S$ is a $p$-group for some prime $p$, and $S$ is maximal in the poset of $p$-subgroups of $\ca L$. 

A fundamental property of localities is given by I.2.6 and I.2.7. Namely, let $(\ca L,\D,S)$ be a 
locality, and let $w=(g_1,\cdots,g_n)$ be a word in the free monoid $\bold W(\ca L)$. Define $S_w$ to be 
the set of all elements $x=x_0\in S$ such that $x_0$ is conjugated to an element $x_1$ of $S$ by $g_1$, 
$x_1$ is conjugated to an element $x_2$ of $S$ by $g_2$, and so on. Then 
\roster 

\item "{(*)}" $S_w$ is a subgroup of $S$, and $w\in\bold D$ if and only if $S_w\in\D$. 

\endroster 
This result is involved in virtually every argument, and should require (and will receive) no further 
reference. In connection with (*) one says of a word $w\in\bold W(\ca L)$ that {\it $w$ is in $\bold D$ 
via $P$} if $P\in\D$ and $P\leq S_w$. A trivial consequence of (*) is: 
\roster 

\item "{(**)}" If $w\in\bold D$ then $S_w\leq S_{\Pi(w)}$. 

\endroster 

Part II introduced the relationship between a locality $(\ca L,\D,S)$ and its fusion system 
$\ca F=\ca F_S(\ca L)$. Thus, $\ca F$ is the fusion system on $S$ whose isomorphisms are 
compositions of restrictions of conjugation maps $c_g:S_g\to S_{g\i}$ for $g\in\ca L$. 

Most readers will have at least some familiarity with general fusion systems, but sections 1 and 6 of 
Part II provide all that will be needed here. The notion of ``saturation" appears in the proof of 
Therem II.6.1, but will play no further role. As mentioned in the introduction, 
the notion of $\ca F$-closed set of subgroups of $S$ ($\ca F$ a fusion system on $S$), and the 
definitions of the sets $\ca F^{cr}$, $\ca F^c$, $\ca F^q$, and $\ca F^s$ appear early in section II.2. 
Theorem II.6.2 (due to Henke [He2]) establishes the sequence of inclusions 
$$ 
\ca F^{cr}\sub\ca F^c\sub\ca F^q\sub\ca F^s, 
$$ 
and establishes that $\ca F^s$ is $\ca F$-closed.

The locality $(\ca L,\D,S)$ on $\ca F$ is defined to be proper (cf. II.2.4) if $\ca F^{cr}\sub\D$, and 
if each of the groups $N_{\ca L}(P)$ for 
$P\in\D$ is of characteristic $p$. Lemma II.2.8 provides a ``dictionary", translating back and forth 
between fusion-theoretic conditions on a given $P\in\D$ and properties of $N_{\ca L}(P)$. Other basic 
notions from Parts I and II will be recalled as needed. 

\vskip .1in 
For any locality $(\ca L,\D,S)$, define 
$$ 
O_p(\ca L)=\bigcap\{S_w\mid w\in\bold W(\ca L)\}. 
$$ 
Equivalently (by I.2.14) $O_p(\ca L)$ is the largest subgroup of $S$ which is a partial normal subgroup 
of $\ca L$. Then $O_p(\ca L)\leq O_p(\ca F)$, where $\ca F$ is the fusion system of $\ca L$. In the case 
that $\ca L$ is proper one has the important equality $O_p(\ca L)=O_p(\ca F)$ (cf. II.2.3).

\definition {Definition 1.1} Let $(\ca L,\D,S)$ be a locality and let $\ca H\leq\ca L$ be a partial 
subgroup of $\ca L$. Set $T=S\cap\ca H$, let $\G$ be a set of subgroups of $T$, and set 
$$ 
\ca H_\G=\{h\in\ca H\mid S_h\cap T\in\G\}. 
$$ 
Let $\ca E:=\ca F_T(\ca H,\G)$ be the fusion system on $T$ generated by the set of conjugation maps 
$c_h:S_h\cap T\to T$ for $h\in\ca H_\G$. Then $(\ca H,\G)$ is a {\it localizable pair} in $\ca L$ if: 
\roster 

\item "{(1)}" $\G$ is $\ca E$-closed. 

\item "{(2)}" $T$ is maximal in the poset of $p$-subgroups of $\ca H$. 

\item "{(3)}" $w\in\bold W(\ca H)$ and $S_w\cap T\in\G$ $\implies$ $w\in\bold D(\ca L)$. 

\endroster 	
\enddefinition

\proclaim {Lemma 1.2} Let $(\ca L,\D,S)$ be a locality, and let $(\ca H,\G)$ be a localizable pair 
in $\ca L$. Set $T=S\cap\ca H$ and $\ca E=\ca F_T(\ca H,\G)$. Further, set 
$$ 
\ca H_\G=\{h\in\ca H\mid S_h\cap T\in\G\}, 
$$ 
and 
$$ 
\bold D(\ca H_\G)=\{w\in\bold W(\ca H_\G)\mid S_w\cap T\in\G\}. 
$$  
Then $\ca H_\G$ is a partial group with respect to the restriction of the product $\Pi:\bold D(\ca L)\to\ca L$ 
to a mapping $\bold D(\ca H_\G)\to\ca H_\G$, and with respect to the restriction to $\ca H_\G$ of the 
inversion in $\ca L$. Moreover, $(\ca H_\G,\G,T)$ is then a locality on $\ca E$.  
\endproclaim 

\demo {Proof} We must verify that the conditions (1) through (4) in the definition (I.1.1) of a partial 
group are satisfied by $\ca H_\G$ and the given product and inversion. Thus $\ca H_\G$ is the set of words 
of length 1 in $\bold D(\ca H_\G)$, and it is plain that if the word $u\circ v$ is in 
$\bold D(\ca H_\G)$ then so 
are $u$ and $v$. Thus the condition (1) holds. If $w\in\bold D(\ca H_\G)$ then $w\in\bold D$, and 
$\Pi(w)$ is in $\ca H$ since $\ca H$ is a partial subgroup of $\ca L$. Then $\Pi(w)\in\ca H_\G$ since 
$S_w\leq S_{\Pi(w)}$, and thus $\Pi$ restricts to a mapping $\bold D(\ca H_\G)\to\ca H_\G$. The conditions 
1.1(2) and (3) are then inherited from $\Pi$. For any $h\in\ca H_\G$ we have $S_{h\i}\cap T=(S_h\cap T)^h$, 
so $\ca H_\G$ is closed under the inversion in $\ca L$. The required condition (4), that 
$\Pi(w\i\circ w)=\1$ for $w\in\bold D(\ca H_\G)$, is then immediate. Thus $\ca H_\G$ is a partial group. 

The inclusion $\ca H_\G\to\ca L$ is a homomorphism since the product in $\ca H_\G$ is the restriction of 
the product $\Pi$. Moreover it is for this same reason, and because $\G$ is $\ca E$-closed, 
that $(\ca H_\G,\G)$ is an objective partial group. 
Notice that any subgroup of $\ca H_\G$ is also a subgroup of $\ca H$. Then 1.1(2) implies that $T$ is 
maximal in the poset of $p$-subgroups of $\ca H_\G$, and thus $(\ca H_\G,\G,T)$ is a locality. By  
the definition of $\ca E$ we have $\ca E=\ca F_T(\ca H_\G)$.  
\qed 
\enddemo 

\definition {Remark} We have already, in II.2.11, encountered one sort of localizable pair. Namely, if $\G$ 
is an $\ca F$-closed subset of $\D$ then $(\ca L,\G)$ is a localizable pair, and $(\ca L_\G,\G,S)$ is the 
restriction $\ca L\mid_\G$ of $\ca L$ to $\G$. 
\enddefinition 

\proclaim {Lemma 1.3} Let the hypothesis and notation be as in 1.2, and let $\ca F_T(\ca H)$ be the fusion 
system on $T$ generated by the conjugation maps $c_h:S_h\cap T\to T$ for all $h\in\ca H$.  
In order for $(\ca H_\G,\G,T)$ to be a 
locality on $\ca F_T(\ca H)$ it suffices that $\G$ contain $\ca E^{cr}$, and that $\ca F_T(\ca H)$ be 
$\G$-generated. 
\endproclaim 

\demo {Proof} Set $\ca E=\ca F_T(\ca H_\G)$ and $\ca E_1=\ca F_T(\ca H)$. Then $\ca E$ is a fusion subsystem 
of $\ca E_1$. Let 
$\ca E_2$ be the fusion subsystem of $\ca E$ generated by the set of all $\ca E$-automorphisms of members 
$U\in(\ca E_1)^{cr}$ such that $U$ is fully normalized in $\ca E_1$ and such that $U=O_p(N_{\ca E_1}(U))$. 
If $(\ca E_1)^{cr}\sub\G$ then $\ca E_2$ is a subsystem of $\ca E$, and if $\ca E_1$ is $\G$-generated 
then $\ca E_1=\ca E_2$. Thus $\ca E=\ca E_1$ if the two stated conditions are fulfilled. 
\qed 
\enddemo

\proclaim {Lemma 1.4} Let $(\ca L,\D,S)$ be a proper locality on $\ca F$, and let $V\leq S$ be 
fully normalized in $\ca F$. Set 
$$ 
\G=\{X\in\D\mid V\norm X\}, \ \ \text{and}\ \ \S=\{Y\leq C_S(V)\mid VY\in\G\}. 
$$ 
\roster 

\item "{(a)}" If $\ca F^c\sub\D$ then $N_{\ca F}(V)^{cr}\sub\G$ and $C_{\ca F}(V)^{cr}\sub\S$. 

\item "{(b)}" If $N_{\ca F}(V)^{cr}\sub\G$ then $(N_{\ca L}(V),\G)$ is a localizable pair and 
$(N_{\ca L}(V)_\G,\G,N_S(V))$ is a proper locality on $N_{\ca F}(V)$. 

\item "{(c)}" If $C_{\ca F}(V)^{cr}\sub\S$ then $(C_{\ca L}(V),\S)$ is a localizable pair and 
$(C_{\ca L}(V)_\S,\S,C_S(V))$ is a proper locality on $C_{\ca F}(V)$. 

\endroster 
\endproclaim 

\demo {Proof} Set $\ca H=N_{\ca L}(V)$ and $\ca K=C_{\ca L}(V)$. Then $\ca H$ is a partial subgroup of 
$\ca L$, and $\ca K\norm\ca H$. Set $\ca E=N_{\ca F}(V)$ and $\ca D=C_{\ca F}(V)$. Then $\ca E$ and 
$\ca D$ are $(cr)$-generated by II.6.1. Then $\ca H_\G$ is a locality on $\ca E$ if $\G\sub\D$, and 
$\ca K_\S$ is a locality on $\ca D$ if $X\S\sub\D$, by 1.3. 

It remains to show that the locality $\ca H_\G$ is proper if $\ca E^{cr}\sub\G$, and that $\ca K_\S$ is 
proper if $\ca D^{cr}\sub\S$. Under these hypotheses on $\G$ and $\S$, it is only necessary to show that 
the normalizers of objects in $\ca H_\G$ or in $\ca K_\S$ are of characteristic $p$. As $\G\sub\D$, and since 
$$ 
N_{\ca H}(X)=N_{N_{\ca L}(X)}(V)
$$ 
for $X\in\G$, it follows from II.2.7(b) that $N_{\ca H}(X)$ is of characteristic $p$. For $Y\in\S$ we have
$VY\in\D$, and so 
$$ 
N_{\ca K}(Y)=C_{N_{\ca H}(VY)}(V)\norm N_{\ca H}(VY).  
$$ 
Then $N_{\ca K}(Y)$ is of characteristic $p$ by II.2.7(a), and the proof is complete. 
\qed 
\enddemo

\proclaim {Corollary 1.5} Let $(\ca L,\D,S)$ be a proper locality on $\ca F$, and let $T\leq S$ be 
strongly closed in $\ca F$. Set $\ca L_T=N_{\ca L}(T)$ and $\ca C_T=C_{\ca L}(T)$. 
\roster 

\item "{(a)}" $(\ca L_T,\D,S)$ is a proper locality on $N_{\ca F}(T)$.  

\item "{(b)}" Set $\S=\{V\leq C_S(T)\mid VT\in\D\}$, and assume that $C_{\ca F}(T)^{cr}\sub\S$. 
Then $(\ca C_T,\S,C_S(T))$ is a proper locality on $C_{\ca F}(T)$. Moreover, the condition 
$C_{\ca F}(T)^{cr}\sub\S$ is fulfilled if $\ca F^c\sub\D$. 

\item "{(c)}" Let $(\ca L^+,\D^+,S)$ be an expansion of $\ca L$ and set $\ca L_T^+=N_{\ca L^+}(T)$. 
Then $(\ca L_T^+,\D^+,S)$ is an expansion of $(\ca L_T,\D,S)$. 

\endroster 
\endproclaim 

\demo {Proof} As $T$ is strongly closed in $\ca F$ one observes that 
$$ 
N_{\ca F}(T)^{cr}=\{P\in\ca F^{cr}\mid T\leq P\},  
$$ 
and hence $N_{\ca F}(T)^{cr}\sub\D$. Then 1.4(a) yields (a). Point (b) is immediate from 1.4(b), and 
point (c) is immediate from (a). 
\qed 
\enddemo

\vskip .2in 
\noindent 
{\bf Section 2: The basic setup} 
\vskip .1in

Most of the results to be proved in sections 2 through 5 will be concerned with a single partial normal 
subgroup $\ca N$ of $\ca L$. (The only exceptions are 2.10 and 2.11, which concern pairs of partial 
normal subgroups.) The following notation will remain fixed.

\definition {2.1 (Basic setup)} $(\ca L,\D,S)$ is a proper locality on $\ca F$, and $\ca N\norm\ca L$ 
is a partial normal subgroup of $\ca L$. Set $T=S\cap\ca N$, and let $\ca E$ be the fusion system 
$\ca F_T(\ca N)$ on $T$, generated by the conjugation maps $c_g:S_g\cap T\to T$ for $g\in\ca N$. 
Further, set $\ca L_T=N_{\ca L}(T)$ and $\ca C_T=C_{\ca L}(T)$. 
\enddefinition 

Recall from I.3.1 that $T$ is strongly closed in $\ca F$, and $T$ is maximal in the poset of $p$-subgroups 
of $\ca N$. There is no reason to suppose that $\ca E$ should be inductive (see II.1.11), 
or that $\ca E$ should remain invariant under the process $\ca N\maps\ca N^+$ of expansion given 
by Theorem II.A2. Indeed, this non-rigidity of $\ca E$ relative to $\ca F$ will be the source of most 
of the technical difficulties that will be encountered.

\proclaim {Lemma 2.2} Assume the setup of 2.1, set $\w T=C_S(T)T$, and set $H=N_{\ca L}(\w T)$. Then 
the following hold. 
\roster 

\item "{(a)}" $H$ is a subgroup of $\ca L$, and $O_p(H)=O_p(N_{\ca F}(\w T))\in\ca F^{cr}$. 

\item "{(b)}" $\w T$ is strongly closed in $\ca F$. 

\item "{(c)}" $\ca F_S(H)=N_{\ca F}(\w T)$. 

\item "{(d)}" $\ca L_T=\ca C_T H$. 

\endroster 
\endproclaim 

\demo {Proof} Set $\w T=C_S(T)T$, and let $\w{\ca N}$ be the partial subgroup $\<\ca N,\ca C_T\>$ 
of $\ca L$ generated by $\ca N$ and $\ca C_T$. Then $\w{\ca N}\norm\ca L$ and $\w T=S\cap\w{\ca N}$, by 
I.5.5. Then (b) follows from I.3.1(a), with $\w{\ca N}$ in the role of $\ca N$. Point (d) is immediate 
from the Frattini Lemma (I.3.11). 

Clearly $\w T$ is $\ca F$-centric. Set $Q=O_p(N_{\ca F}(\w T))$. Then $Q\norm S$, so 
$\w T\norm O_p(N_{\ca F}(Q))$, and hence $Q=O_p(N_{\ca F}(Q))$. Thus 
$$ 
Q\in\ca F^{cr}, \tag*
$$ 
so $Q\in\D$, and $N_{\ca L}(Q)$ is a subgroup of $\ca L$. 

By 1.5(a), $(H,\D,S)$ is a proper locality on $N_{\ca F}(\w T)$. This yields (c), and II.2.3 shows that 
$O_p(H)=Q$. Then $H$ is the subgroup $N_{N_{\ca L}(Q)}(\w T)$ of $\ca L$. This result, together with (*), 
completes the proof of (a). 
\qed 
\enddemo

\proclaim {Lemma 2.3} Assume the setup of 2.1. Let $\ca H$ be a partial subgroup of $\ca L_T$ having 
the property that $(h\i,x,h)\in\bold D$ for all $x\in\ca N$ and all $h\in\ca H$, and let 
$\l:\ca H\to Aut(T)$ be the homomorphism which sends $h\in\ca H$ to conjugation by $h$. Then 
$Im(\l)\sub Aut(\ca E)$. 
\endproclaim 

\demo {Proof}  Let $\phi:U\to U'$ be an $\ca E$-isomorphism. By definition, there is a sequence 
$$ 
U=U_0 @>\phi_1>>\cdots@>\phi_n>>U_n=U'
$$ 
of $\ca E$-isomorphisms, such that each $\phi_i:U_{i-1}\to U_i$ is the restriction of a conjugation map 
$c_{x_i}:T\cap S_{g_i}\to T$, with $x_i\in\ca N$. Let $h\in\ca H$. By hypothesis, $(h\i,x_i,h)\in\bold D$ 
for all $i$, and then $(x_i)^h\in\ca N$ as $\ca N\norm\ca L$. As $\ca H\leq N_{\ca L}(T)$ we may define 
$V_i:=(U_i)^h$, and then define $\psi_i:V_{i_1}\to V_i$ to be given by conjugation by $(x_i)^h$. Then 
each $\psi_i$ is an $\ca E$-isomorphism, and the composite $\psi:=\psi_1\circ\cdots\circ\psi_n$ is 
given by $c_{h\i}\circ\phi\circ c_h$ as an $\ca E$-isomorphism $U^h\to(U')^h$. This shows that 
the conjugation map $c_h:T\to T$ is $\ca E$-fusion-preserving. That is, $c_h\in Aut(\ca E)$. 
The verification that the map $h\maps c_h$ is a homomorphism $\ca H\to Aut(\ca E)$ is a straightforward 
application of I.2.3(c). 
\qed 
\enddemo 

\proclaim {Corollary 2.4} Suppose that $x^h$ is defined for all $x\in\ca N$ and all 
$h\in H:=N_{\ca L}(C_S(T)T)$. Let $\G$ be a set of subgroups of subgroups of $T$ which is both 
$\ca E$-invariant and $Aut(\ca E)$-invariant. Then $\G$ is $\ca F$-invariant. 
\endproclaim 

\demo {Proof} Let $g\in\ca L$, and let $U\in\G$ with $U\leq S_g\cap T$. Let $\phi:U\to U^g$ be the conjugation 
map. Define $\ca L_T$ and $\ca C_T$ as in 1.4. By the Splitting Lemma (I.3.12) we may write $g=xf$ where 
$x\in\ca N$, $f\in\ca L_T$, and $S_g=S_{(x,f)}$. As $\ca C_T\norm\ca L_T$ we may also write $f=yh$ where 
$y\in\ca C_T$, $h\in H$, and $S_f=S_{(y,h)}$. Then $S_g=S_{(x,y,h)}$, and $\phi$ is then the composition of 
the $\ca E$-isomorphism $c_x$ followed by $c_h$. Then $U^g\in\G$ by 2.3. Then $\G$ is $\ca F$-invariant 
since $\ca F$ is generated by the conjugation maps $c_g$ with $g\in\ca L$. 
\qed 
\enddemo 

\proclaim {Lemma 2.5} Assume the setup of 2.1, and let $Q\leq T$ be a subgroup of $T$ such that $Q$ is fully 
normalized in $\ca E$. Suppose that $O_p(\ca L)Q\in\D$, and set $K=N_{\ca L}(Q)$. Then $K$ is a normal 
subgroup of the group $N_{\ca L}(Q)$, $N_T(Q)$ is a Sylow subgroup of $K$, and 
$N_{\ca E}(Q)=\ca F_{N_T(Q)}(K)$. 
\endproclaim 

\demo {Proof} Set $M=N_{\ca L}(O_p(\ca L)Q)$. As $O_p(\ca L)Q\in\D$, $M$ is a subgroup of $\ca L$. Then 
$K=\ca N\cap N_M(Q)=\ca N\cap N_{\ca L}(Q)$ is a normal subgroup of the group $N_{\ca L}(Q)$ by I.1.8. 

Set $X=N_T(V)$ and let $Y$ be a Sylow $p$-subgroup of $K$ containing $X$. By I.2.11 there exists $g\in\ca L$ 
with $Y^g\leq S$, and then $Y^g\leq T$ since $T$ is strongly closed in $\ca F$. Employ the splitting lemma 
(I.3.12) to obtain $g=fh$ with $f\in\ca N$, $h\in N_{\ca L}(T)$, and with $S_g=S_{(f,h)}$. Then 
$Y^f\leq T$, and so $Y^f\leq N_T(V^f)$. As $V$ is fully normalized in $\ca E$, we conclude that $X^f=Y^f$, 
and thus $X\in Syl_p(K)$. 

Next, let $A$ and $B$ be subgroups of $X$ containing $V$, such that $A$ and $B$ are conjugate in 
$N_{\ca E}(V)$. Let $\g:A\to B$ be an $N_{\ca E}(V)$-isomorphism. By the definition of 
$\ca E=\ca F_T(\ca N)$, this means that there exists $w=(f_1,\cdots f_n)\in\bold W(\ca N)$ such that 
$A\leq S_w$, and such that $\g$ is given by composing the conjugation maps $c_{f_i}$. As 
$O_p(\ca L)V\leq S_w$ we have $w\in\bold D$, and then $\Pi(w)\in\ca N$. 
Then $\Pi(w)\in K$, and thus $N_{\ca E}(Q)=\ca F_X(K)$. 
\qed 
\enddemo

Recall from II.1.8 that $\ca E^{cr}$ is defined to be the set of all $U\in\ca E^c$ such that there exists 
an $\ca E$-conjugate $V$ of $U$ such that $V$ is fully normalized in $\ca E$, and such that 
$V=O_p(N_{\ca E}(V))$. It has already been remarked that this definition does not quite agree with 
the usual one (see [AKO], for example. The following lemma provides the justification for this discrepancy. 
Namely, points (b) and (c) of the lemma establish a ``descent" from $\ca F^{cr}$ to $\ca E^{cr}$ 
which relies on our definition, and which would otherwise be lacking.

\proclaim {Lemma 2.6} Assume the setup of 2.1. 
\roster 

\item "{(a)}" Let $U\leq T$ be a subgroup of $T$. Assume that $U$ is fully normalized in $\ca F$ 
and that $C_T(U)\leq U$. Then $U^{\ca F}\sub\ca E^c$.  

\item "{(b)}" Let $Q\in\ca F^{cr}$ and set $U=Q\cap T$. Suppose that $U$ is fully normalized in $\ca F$. 
Then $U^{\ca F}\sub\ca E^c$, $U$ is fully normalized in $\ca E$, and $U=O_p(N_{\ca E}(U))$. 

\item "{(c)}" Let $\ca M\norm\ca L$ be a partial normal subgroup of $\ca L$ containing $\ca N$. Set 
$R=S\cap\ca M$, set $\ca D=\ca F_R(\ca M)$, and suppose that $\ca D$ is inductive. Let $Q\in\ca D^{cr}$, 
and set $U=Q\cap T$. Assume that $O_p(\ca L)Q\in\D$ and that $U$ is fully normalized in $\ca D$. 
Then $U^{\ca D}\sub\ca E^c$, $U$ is fully normalized in $\ca E$, and $U=O_p(N_{\ca E}(U))$. 

\endroster 
\endproclaim 

\demo {Proof} Let $U$ be chosen as in (a), and let $U'$ be an $\ca F$-conjugate of $U$. As $\ca F$ is 
inductive by II.6.1, there exists an $\ca F$-homomorphism $\phi:N_S(U')\to N_S(U)$ with $U'\phi=U$. As $T$ 
is strongly closed in $\ca F$ we obtain $C_T(U')\phi\leq C_T(U)$. As $C_T(U)\leq U$ it follows that 
$C_T(U')\leq U'$. Since $\ca E$-conjugates of $\ca F$-conjugates of $U$ are $\ca F$-conjugates of $U$, 
we obtain $U^{\ca F}\sub\ca E^c$, and thus (a) holds.   

Point (b) is the special case of (c) where $\ca D=\ca F$, so it remains only to prove (c). Let $\ca M$, $R$, 
$\ca D$, $Q$, and $U$ be as stated in (c). Then $U$ is fully normalized in $\ca E$ by II.1.17. 
Set $K=N_{\ca M}(Q)$. By 2.5, with $\ca M$ in the role of $\ca N$, we have $N_R(Q)\in Syl_p(K)$ and 
$N_{\ca D}(Q)=\ca F_{N_R(Q)}(K)$. Also by 2.5, $K$ is a normal subgroup of the group $N_{\ca L}(Q)$, 
where $N_{\ca L}(Q)$ is a local subgroup of $N_{\ca L}(O_p(\ca L)Q)$. Then $K$ is of characteristic $p$ 
by II.2.7, and hence $O_p(K)=O_p(N_{\ca D}(Q))$. As $Q\in\ca D^{cr}$, where $\ca D$ is inductive, it 
follows from II.1.14 that $Q=O_p(K)$. Set $D=N_{C_T(U)}(Q)$. Then  
$$ 
[Q,D]\leq C_{Q\cap T}(V)=Z(U), 
$$ 
and thus $D$ centralizes the chain $Q\geq U\geq 1$ of normal subgroups of $K$. Then $D\leq Q$ by II.2.7(c), 
and thus $C_T(U)\leq Q$. Then $C_T(U)\leq Q\cap T=U$. A straightforward variation on the proof of (a) 
then yields $U^{\ca D}\sub\ca E^c$. 

Set $X=O_p(N_{\ca E}(U))$. Then 
$U\leq X\leq T$. By 1.5 the elements of $S$ act as automorphisms of $\ca E$ by conjugation, so $Q$ acts 
on $N_{\ca E}(X)$, and thus $X$ is $Q$-invariant. By definition, each $\ca E$-automorphism of $U$ 
extends to an $\ca E$-automorphism of $X$, and this translates into the following statement. 
\roster 

\item "{(*)}" Let $\b$ be an $\ca E$-automorphism of $U$, let $\bar\b$ be an extension of $\b$ to an 
$\ca E$-automorphism of $X$, and let $x\in X$. For any $a\in X$, let $c_a$ be the automorphism of $U$ 
given by conjugation by $a$. Then $\b\i\circ c_a\circ\b=c_{a\bar\b}$. In particular, $Aut_X(U)$ 
is a normal subgroup of $Aut_{\ca E}(U)$. 

\endroster 
Set $K_0=K\cap\ca N$ and set $A=N_X(Q)$. Then $K_0$ is a normal subgroup of $K$, while (*) yields 
$Aut_A(U)\norm Aut_{K_0}(U)$. As $U\in\ca E^c$, $C_{K_0}(U)$ is the direct product 
of $Z(U)$ with a normal $p'$-subgroup of $K_0$. Here $O_{p'}(K_0)=1$ as $K$ is of characteristic 
$p$, and so $C_{K_0}(U)=Z(U)$. There is then a natural isomorphism of $Aut_{K_0}(U)$ with $K_0/Z(U)$, 
from which it follows that $A/Z(U)\norm K_0/Z(U)$. Then $A\norm K_0$, so $A\leq O_p(K_0)\leq O_p(K)=Q$. 
Then $X\leq Q$, and since also $X\leq T$ we arrive at $X\leq U$. Thus $X=U$, completing the proof of (c). 
\qed 
\enddemo 

For any partial subgroup $\ca H\leq\ca L$ define $Z(\ca H)$ to be the set of all $z\in\ca H$ such that, 
for all $h\in\ca H$, $h^z$ is defined and is equal to $h$. 

\proclaim {Lemma 2.7} Assume the setup of 2.1. Then $Z(\ca N)=C_T(\ca N)$. 
\endproclaim 

\demo {Proof} As $t^z$ is defined and is equal to $t$ for all $t\in T$ and all $z\in Z(\ca N)$, we 
have $Z(\ca N)\leq N_{\ca N}(T)$. Set $H=N_{\ca L}(C_S(T)T)$. As $\ca L$ is proper, I.3.5 yields 
$N_{\ca N}(T)\leq H$. As $H$ is a subgroup of $\ca L$ by 2.2(a), $Z(\ca N)$ is then a normal abelian subgroup 
of $H$. As $H$ is of characteristic $p$, II.2.7 implies that $Z(\ca N)$ is of characteristic $p$. 
Thus $Z(\ca N)$ is a $p$-subgroup of $N_{\ca L}(T)$. As $T$ is a maximal $p$-subgroup of $\ca N$ by 
I.3.1(a), the lemma follows. 
\qed 
\enddemo

For any fusion system $\ca F$ on a $p$-group $S$, and any non-empty set $\G$ of subgroups of $S$, define 
the {\it $\ca F$-closure} of $\G$ to be the smallest $\ca F$-closed set of subgroups of $S$ containing $\G$. 
Thus, the $\ca F$-closure of $\G$ is the set of subgroups $X$ of $S$ such that $X$ contains an 
$\ca F$-conjugate of a member of $\G$. 
\vskip .1in

For any pair $\G$ and $\S$ of non-empty sets of subgroups of $S$, write $\G\S$ for the set of all products 
$XY$ with $X\in\G$ and $Y\in\S$. If $\G=\{X\}$ is a singleton we may write $X\S$ for $\{X\}\S$.

\proclaim {Lemma 2.8} Assume the setup of 2.1, and set $X=O_p(\ca L)C_S(\ca N)$. Let $\G$ be an 
$\ca E$-closed set of subgroups of $T$ containing $\ca E^{cr}$, and let $\D_0$ be the $\ca F$-closure of 
$X\G$. Assume that $\D_0\sub\D$ and that $\ca E$ is $\G$-generated. Then $(\ca N,\G)$ is a localizable 
pair, and $(\ca N_\G,\G,T)$ is a proper locality on $\ca E$. Moreover, if $O_p(\ca L)\G\sub\D$,  then the 
following hold. 
\roster 

\item "{(a)}" $\ca F^{cr}\sub\D_0\sub\D$, the restriction $\ca L_0$ of $\ca L$ to $\D_0$ is a proper 
locality $(\ca L_0,\D_0,S)$ on $\ca F$, and $\ca N_\G=\ca N\cap\ca L_0$.  

\item "{(b)}" $N_{\ca L}(T)$ is a subgroup of $\ca L_0$, and $x^g$ is defined for all $x\in\ca L_0$ and 
all $g\in N_{\ca L}(XT)$. Moreover, the mapping $g\to c_g$ which sends $g\in N_{\ca L}(XT)$ to the 
conjugation map $c_g:\ca L_0\to\ca L_0$ is a homomorphism $\l:N_{\ca L}(XT)\to Aut(\ca L_0)$. 

\item "{(c)}" For each $g\in N_{\ca L}(T)$, the restriction of $c_g$ to $T$ is an automorphism of $\ca E$, 
and one obtains in this way a homomorphism $\l_T:N_{\ca L}(T)\to Aut(\ca E)$. 

\endroster 
\endproclaim 

\demo {Proof} We check that $(\ca N,\G)$ is a localizable pair by verifying the two conditions 
in definition 1.1. The second of these - that $T$ is maximal in the poset of $p$-subgroups of $\ca N$ - 
is given by I.3.1(c). 

Let $w\in\bold W(\ca N)$ with $S_w\cap T\in\G$. Then $X(S_w\cap T)\in\D$ by hypothesis, so $S_w\in\D$ and 
$w\in\bold D$. Thus the condition 1.1(1) obtains, and $(\ca N,\G)$ is a localizable pair. By 1.2  
$(\ca N_\G,\G,T)$ is a locality. Since $\ca E^{cr}\sub\G$ and $\ca E$ is $\G$-generated by hypothesis, 1.3 
yields $\ca E=\ca F_T(\ca N_\G)$.

Let $U\in\G$ and set $P=XU$. Then $P\in\D$, so $N_{\ca L}(P)$ is a group of characteristic $p$.  
Then $N_{\ca N}(P)$ is of characteristic $p$ by II.2.7(a), and then $N_{\ca N}(U)$ is of characteristic $p$ 
by II.2.7(b). Thus $(\ca N_\G,\G,T)$ is a proper locality on $\ca E$.  

We assume for the remainder of the proof that $O_p(\ca L)\G\sub\D$. Let $R\in\ca F^{cr}$ and, via II.1.18, 
let $R'$ be an $\ca F$-conjugate of $R$ such that both $R'$ and 
$R'\cap T$ are fully normalized in $\ca F$, and then $R'\cap T\in\ca E^{cr}$ by 2.6(b). Then 
$O_p(\ca N)\leq R'$. Thus $R'\in\D_0$, and then also $R\in\D_0$. Thus $\ca F^{cr}\sub\D_0$, and the 
restriction $(\ca L_0,\D_0,T)$ of $\ca L$ is then given by II.2.11 as a proper locality on $\ca F$. 
The equality $\ca N_\G=\ca N\cap\ca L_0$ is immediate from the definitions of $\ca N_\G$ and of $\ca L_0$. 
Thus (a) holds. 

Set $H=N_{\ca L}(T)$. As $O_p(\ca L)T=XT\in\D_0$, $H$ is a subgroup of $\ca L_0$. Let $x\in\ca L_0$, let 
$h\in H$, and set $P=S_x\cap O_p(\ca L)T$. Then $P\in\D$ (by the definition of $\ca L_0$), and 
$(h\i,x,h)\in\bold D$ via $P^h$. Points (b) and (c) now follow from 2.3. 
\qed 
\enddemo

The remainder of this section involves the following variation on the setup of 2.1.

\definition {2.9 (Product setup)} $(\ca L,\D,S)$ is a proper locality on $\ca F$, and $\ca N_1$ and $\ca N_2$ 
are partial normal subgroups of $\ca L$. Set $T_i=S\cap\ca N_i$, and set $\ca E_i=\ca F_{T_i}(\ca N_i)$. 
Assume: 
$$ 
S\cap\ca N_i\leq C_S(\ca N_{3-i})\ \ (i=1,2).\tag* 
$$ 
Set $\ca M=\ca N_1\ca N_2$ (a partial normal subgroup of $\ca L$ by I.5.1, or by [He]), and set 
$R=S\cap\ca M$ and $\ca D=\ca F_R(\ca M)$. 
\enddefinition 

\proclaim {Lemma 2.10} Assume the setup of 2.9. Assume further that $\ca D^{cr}\sub\D$ and that $\ca D$ is 
$\ca D^{cr}$-generated. Set $\G=\{P\in\D\mid P\leq R\}$, and set 
$$ 
\G_i=\{PT_{3-i}\cap T_i\mid P\in\G\}, \ \ (i=1,2). 
$$ 
Let $\D_0$ be the overgroup-closure of $\G_1\G_2$ in $S$, and let $\D_0^+$ be the overgroup-closure of 
$\G$ in $S$. Then the following hold. 
\roster 

\item "{(a)}" $\D_0$ and $\D_0^+$ are $\ca F$-closed subset of $\D$, and $\ca F^{cr}\sub\D_0\sub\D_0^+$. 

\item "{(b)}" $(\ca M,\G)$ and $(\ca N_i,\G_i)$ $(i=1,2)$ are localizable pairs, $(\ca M_{\G},\G,R)$ is a 
proper locality on $\ca D$, and each $(({\ca N_i})_{\G_i},\G_i,T_i)$ is a proper locality on $\ca E_i$. 
Moreover, $({\ca N_i})_{\G_i}=\ca N_i\cap\ca M_\G$. 

\item "{(c)}" $R\in\D_0$, and the group $N_{\ca L}(R)$ acts by conjugation on all three of the localities 
in (b). Indeed, conjugation by $h\in N_{\ca L}(R)$ is an automorphism of $\ca D$ which restricts to an 
automorphism of each $\ca E_i$. 

\item "{(d)}" $\ca D^{cr}=(\ca E_1)^{cr}(\ca E_2)^{cr}$. 

\item "{(e)}" Let $Y\leq R$ be a subgroup of $R$, set $Y_i=Y\cap T_i$, and suppose that $Y=Y_1Y_2$. 
Then $Y\in\ca D^c$ if and only if each $Y_i$ is in $(\ca E_i)^c$, and $Y\in\ca D^s$ if and only if  
$Y_i\in(\ca E_i)^s$ $(i=1,2)$. 

\endroster 
\endproclaim 

\demo {Proof} By definition, $\ca D^{cr}$ is $\ca D$-invariant, so $\G$ is $\ca D$-closed. Since 
$\ca D$ is $\ca D^{cr}$-generated and $\ca D^{cr}\sub\D$, by hypothesis, we may appeal to 2.8 with 
$\ca M$ in the role of $\ca N$. Thus $(\ca M_\G,\bar\G)$ is a localizable pair, and 
$(\ca M_\G,\bar\G,R)$ is a proper locality on $\ca D$. 

We have $R=R_1R_2=T_1T_2$ by I.5.1, so $R\in\D_0$. For any $P\in\G$ set $P_i=PT_{3-i}\cap T_i$. Then 
$P\leq P_1P_2$, and thus $\G_1\G_2\sub\G$. As $\G\sub\D$ we then have $\D_0\sub\D_0^+\sub\D$. 
Clearly $\G$ is $\ca F$-invariant, so $\D_0^+$ is $\ca F$-closed. For any $P=P_1P_2\in\G_1\G_2$ and any 
$g\in\ca L$ with $P\leq S_g$ we have $P^g=(P_1)^g(P_2)^g\in\G_1\G_2$, and it follows that $\D_0$ is 
$\ca F$-closed. 
 
Set $H=N_{\ca L}(R)$. Then $H$ is a subgroup of $\ca L$ as $R\in\G_1\G_2$. Evidently $\G$ is 
$H$-invariant, so $H$ acts on $\ca M_\G$, and then $H$ acts on $\ca D$ by 2.3. Thus (b) and (c) hold 
insofar as these points refer only to $\ca M$ and not to $\ca E_i$.  

Let $A$ and $B$ be $\ca D$-conjugate subgroups of $R$, and let $\phi:A\to B$ be a $\ca D$-isomorphism. By 
definition, this means that there exists $w=(g_1,\cdots,g_n)\in\bold W(\ca M)$ such that $\phi$ may be 
written as a composition $c_w$ of conjugation maps: 
$$ 
A=A_0@>c_{g_1}>>A_1@>>>\cdots@>>>A_{n-1}@>c_{g_n}>>A_n=B. \tag*
$$ 
By I.5.2 each $g_k$ is a product $g_k=x_ky_k$ with $x_k\in\ca N_1$, $y_k\in\ca N_2$, and with 
$S_{g_k}=S_{(x_k,y_k)}$. Set 
$$
w'=(x_1,y_1,\cdots,x_k,y_k),\ \ u=(x_1,\cdots,x_k), \ \ \text{and}\ \  v=(y_1,\cdots,y_k).  
$$ 
Then $\phi=c_{w'}$, and then $\phi=c_u\circ c_v$ since $T_i\leq C_S(\ca N_{3-1})$. Thus: 
\roster 

\item "{(1)}" Each $\ca D$-isomorphism may be factored as a composition $\psi_1\circ\psi_2$, where 
$\psi_i$ is a composition of conjugation maps by elements of $\ca N_i$ ($i=1,2)$. 

\endroster 
As a consequence: 
\roster 

\item "{(2)}" Let $X_i\leq T_i$ $(i=1,2)$, and set $X=X_1X_2$. Then 
$X^{\ca D}=(X_1)^{\ca E_1})(X_2)^{\ca E_2}$. In particular, $X$ is fully normalized in $\ca D$ 
if and only if each $X_i$ is fully normalized in $\ca E_i$. 

\endroster 
Let $U_i\in(\ca E_i)^{cr}$ $(i=1,2)$, and set $X=U_1U_2$. There then exists an $\ca E_i$-conjugate $V_i$ 
of $U_i$ such that $V_i$ is fully normalized in $\ca E_i$, and such that $V_i=O_p(N_{\ca E_i}(V_i))$. Set 
$Y=V_1V_2$. Then (2) shows that $Y\in X^{\ca D}$ and that $Y$ is fully normalized in $\ca D$. 
We compute: 
$$ 
\align 
C_R(V_1V_2)&=C_R(V_1)\cap C_R(V_2)=C_{T_1}(V_1)T_2\cap C_{T_2}(V_2)T_1\leq V_1T_2\cap V_2T_1 \\ 
 &=V_1(T_2\cap V_2T_1)=V_1V_2(T_1\cap T_2)\leq V_1V_2, 
\endalign 
$$ 
since $T_1\cap T_2\leq Z(T_1)\leq V_1$. Thus $V_1V_2\in\ca D^c$. In fact, we have shown: 
\roster 

\item "{(3)}" $(\ca E_1)^{c}(\ca E_2)^{c}\sub\ca D^c$. 

\endroster 

Let $\phi_1:A\to B$ be an $N_{\ca E_1}(V_1)$-isomorphism, where $A$ and $B$ are subgroups of $N_{T_1}(V_1)$ 
containing $V_1$. Then $\phi_1$ is a $\ca D$-isomorphism, and evidently $\phi$ extends to a 
$\ca D$-isomorphism 
$\g:AV_2\to BV_2$ whose restriction to $V_2$ is the identity map. Set $P=O_p(N_{\ca D}(Y))$. As $\g$ is in 
fact a $N_{\ca D}(Y)$-isomorphism, $\g$ extends to a $\ca D$-isomorphism $\psi:AP\to BP$ which leaves $P$ 
invariant. 
 
Write $\psi=\psi_1\circ\psi_2$ as in (1). Thus,
$$ 
AP@>\psi_1>> CQ@>\psi_2>> BP, 
$$ 
where $C=A\psi_1$ and $Q=P\psi_1=P\psi_2\i$. Then  
$$ 
Q\leq APT_1\cap PT_2=P(T_1\cap T_2)=P,  
$$ 
since 
$$ 
T_1\cap T_2\leq Z(T_1)\cap Z(T_2)\leq V_1\cap V_2\leq V_1V_2\leq P. 
$$ 
This shows that $\psi_1:AP\to BP$ is an extension of $\phi_1$. There is an obvious further extension of 
$\psi_1$ to a $\ca D$-isomorphism $\l_1:APT_2\to BPT_2$, given in the same way as $\psi_1$ as a 
composite of conjugations by elements of $\ca N_1$. Thus $\phi_1$ extends to an $\ca E_1$-homomorphism 
$A(PT_2\cap T_1)\to B(PT_2\cap T_1)$. As $V_1=O_p(N_{\ca E_1}(V))$ we conclude that 
$V_1=PT_2\cap T_1$. That is, $V_1$ is the projection of $P$ into $T_1$ relative to the decomposition  
$R=T_1T_2$. Similarly, $V_2$ is the projection of $P$ into $T_2$, and thus $P=V_1V_2$. With (3), this 
shows: 
\roster 

\item "{(4)}" $\D^{cr}\sup(\ca E_1)^{cr}(\ca E_2)^{cr}$. 

\endroster 
Let $P$ now be an arbitrary member of $\ca D^{cr}$. In order to show that $P\in(\ca E_1)^{cr}(\ca E_2)^{cr}$, 
it suffices to consider the case where $P$ is fully normalized in $\ca D$, by (2). Set $M=N_{\ca M}(P)$. Then 
$M$ is a normal subgroup of the group $N_{\ca L}(P)$. Set $P_i=PT_{3-i}\cap T_i$, and set 
$K_i=N_{\ca N_i}(P_i)$. Then $P\leq P_1P_2\norm K_i$, and $P_1P_2\in\D$. Thus $K_i$ is a normal subgroup 
of $N_{\ca M}(P_1P_2)$. We again quote I.5.2 in order to write an arbitrary $g\in M$ as a product 
$g=x_1x_2$ with $x_i\in\ca N_i$ and with $S_g=S_{(x_1,x_2)}$. Then $(P_i)^{x_i}=P_i$, and thus 
$M\leq K_1K_2$. As $\ca M_\G$ is a proper locality, II.2.3 shows that $P=O_p(M)$. Then 
$N_{P_i}(P)\leq P$, and thus $P=P_1P_2$. This completes the proof of (d). 

In order to now show that $(\ca N_i,\G_i)$ is a localizable pair, three points have to be verified. 
First: $T_i$ is maximal in the poset of $p$-subgroups of $\ca N_i$ (I.3.1(c)). Second: $\G_i$ is an 
$\ca E_i$-closed set of subgroups of $T_i$ (true, since $\G_i$ is $\ca F$-invariant). Third: 
$w\in\bold D$ for each $w\in\bold W(\ca N_i)$ 
such that $S_w$ contains some $U\in\G_i$. This third point is a consequence of (1) and (d), since 
$S_w\geq UT_{3-i}\in\D$. So then, $(\ca N_i,\G_i)$ is a localizable pair, and then 
$((\ca N_i)_{\G_i},\G_i,T_i)$ is a locality by 1.2. 

Set $\ca K_i=\ca N_i\cap\ca M_\G$ and set $\ca D_1=\ca F_{T_1}(\ca K_1)$. Then 
$\ca K_i=(\ca N_i)_{\G_i}$, and $(\ca E_1)^{cr}=(\ca D_1)^{cr}$, by (d). Thus $(\ca E_1)^{cr}\sub\G_i$. 
As $\ca D$ is $\ca D^{cr}$-generated, it follows from (1) and (d) that $\ca E_1$ is 
$(\ca E_1)^{cr}$-generated. As $\ca D_1$ is a fusion subsystem of $\ca E_1$, we conclude that 
$\ca E_1=\ca D_1$. Since $\ca M_\G$-normalizers of members of $\G$ are of characteristic $p$, the normalizer 
in $\ca K_i$ of any member of $\G_i$ is of characteristic $p$, so $(\ca N_i)_{\G_i}$ is a proper locality on 
$\ca E_i$. Thus (b) holds. 

Evidently $H$ acts on $\ca M_\G$ by conjugation, and this action restricts to an action on 
$\ca K_i$ for each $i$. The proof of (c) is then completed by 2.8(c). 

In order to complete the proof of (a) it remains to show that $\ca F^{cr}\sub\D_0$. Let $P\in\ca F^{cr}$ and 
let $Q$ be an $\ca F$-conjugate of $P$ such that both $Q$ and $Q\cap R$ are fully normalized in $\ca F$. 
Then $Q\cap R\in\ca D^{cr}$ by 2.6(b). Then (d) shows that $Q\cap R\in\D_0$, so $Q\in\D_0$. As 
$\ca D_0$ is $\ca F$-invariant we obtain $P\in\D_0$, so (a) holds. 

It now remains to prove (e). Let $Y=Y_1Y_2$ be as stated in (e). By (2) we may proceed under the 
assumption that $Y$ is fully normalized in $\ca F$ and that each $Y_i$ is fully normalized in $\ca E_i$. 
Suppose that $Y$ is $\ca D$-centric. Then $C_{T_i}(Y_i)\leq C_S(Y)\cap T_i\leq Y\cap T_i=Y_i$, and thus 
$Y_i\in(\ca E_i)^c$. With (3), we then have the first ``if and only if" in (e). 

Set $Q_i=O_p(N_{\ca E_i}(Y_i)$ and set $Q=O_p(N_{\ca D}(Y))$. Then $Q_1Q_2=Q$, as follows from 
(1). Further, by II.1.16, $Q$ is fully centralized in $\ca D$, and $Q_i$ is fully centralized in $\ca E_i$. 
Since $Q_1Q_2\in\ca D^c$ if and only if each $Q_i$ is centric in $\ca E_i$, we conclude that 
$Y\in\ca D^s$ if and only if each $Q_i\in(\ca E_i)^s$. 
\qed 
\enddemo

\proclaim {Corollary 2.11} Assume the setup of 2.9, with $\ca N_1\ca N_2=\ca L$. Set 
$\G_i=\{PT_{3-i}\cap T_i\mid P\in\D\}$ $(i=1,2)$. Then $(\ca N_i,\G_i,T_i)$ is a proper locality 
on $\ca E_i$, $\G_1\G_2$ is an $\ca F$-closed subset of $\D$, and $\ca L$ is the same partial group as 
its restriction $\ca L\mid_{\G_1\G_2}$ to $\G_1\G_2$. 
\endproclaim 

\demo {Proof} Since $\ca F$ is $(cr)$-generated by II.2.3, and since $\ca F^{cr}\sub\D$, we may apply 
2.10 with $\ca M=\ca L$ and with $\G=\D$. Then $(\ca N_i,\G_i,T_i)$ is a proper locality on $\ca E_i$, 
by 2.10(b), and $\ca L_\G=\ca L$. Set $\D_0=\G_1\G_2$. Then $\D_0$ is $\ca F$-closed, and 
$\ca F^{cr}\sub\D_0\sub\D$ by 2.10(a). Let $g\in\ca L$ and write $g=g_1g_2$ with $g_i\in\ca N_i$ and 
with $S_g=S_{(g_1,g_2)}$ (via I.5.2). Then $S_g=(S_{g_1}\cap T_1)(S_{g_2}\cap T_2)\in\D_0$, so 
$\ca L_{\D_0}=\ca L$.  
\qed 
\enddemo

\vskip .2in
\noindent 
{\bf Section 3: Alperin-Goldschmidt variations} 
\vskip .1in 

For any finite group $G$ let $\G_p(G)$ be the graph whose vertices are the Sylow $p$-subgroups of $G$, and 
whose edges are the pairs $\{X,Y\}$ of distinct Sylow $p$-subgroups such that $X\cap Y\neq 1$. By 
Sylow's theorem, the action of $G$ on $\G_p(G)$ by conjugation is transitive on vertices. 

Let $S$ be a fixed Sylow $p$-subgroup of $G$, let $\S$ be the connected component of $\G_p(G)$ containing 
the vertex $S$, and let $H$ be the set-wise stabilizer of $\S$ in $G$. Then $g\in H$ for all $g\in G$ 
such that $p$ divides $|H\cap H^g|$, and hence $H=G$ if and only if $\G_p(G)$ is connected. We shall say 
that $G$ is {\it $p$-disconnected} if $\G_p(G)$ is disconnected. Otherwise, $G$ is {\it $p$-connected}. 

\definition {Remark} In the standard terminology, one says that a proper subgroup $X$ of $G$ such that 
$p$ divides $|X|$ and such that $g\in K$ whenever $p$ divides $|X\cap X^g|$ is {\it strongly $p$-embedded} 
in $G$. One easily deduces in that case, that $Syl_p(X)$ contains a connected component $\S$ of $\G_p(G)$, 
and that $X$ contains the set-wise stabilizer $H$ of $\S$. Thus $G$ has a strongly $p$-embedded subgroup if 
and only if $\G_p(G)$ is disconnected. 
\enddefinition

The following result is well known. Since it is not so easy to find a reference for it, we provide 
a proof. 

\proclaim {Lemma 3.1} Let $X$ be a $p'$-group and let $A$ be an elementary abelian group of order $p^2$ 
such that $A$ acts on $X$. Then $X=\<C_X(a)\mid 1\neq a\in A\>$. 
\endproclaim 

\demo {Proof} Let $G$ be the semi-direct product $XA$ formed via the action of $A$ on $X$. Let $q$ be a 
prime dividing $|X|$ and let $Y$ be a Sylow $q$-subgroup of $X$. Then $G=N_G(Y)X$, so we may choose $Y$ 
to be $A$-invariant. Then $A$ acts on the elementary abelian $q$-group $V:=Y/\Phi(Y)$. By Maschke's 
Theorem $V$ is a direct sum of irreducible $A$-submodules. By Schur's Lemma $C_A(W)$ contains a 
maximal subgroup of $A$ for each irreducible $A$-submodule $W$ of $V$. As $A$ is assumed to be 
non-cyclic we thereby obtain $V=\<C_V(a)\mid 1\neq a\in A\>$. Here $C_V(a)=C_Y(a)\Phi(Y)/\Phi(Y)$ 
by coprime action (cf. II.2.7(c)), so $Y=\<C_Y(a)\mid 1\neq a\in A\>\Phi(Y)$. It is a basic property 
of the Frattini subgroup that if $Y_0$ is a subgroup of $Y$ such that $Y=Y_0\Phi(Y)$ then $Y=Y_0$. Thus   
we have the lemma in the case that $X=Y$. We conclude that $\<C_X(a)\mid 1\neq a\in A\>$ contains a Sylow 
subgroup of $X$ for each prime divisor of $|X|$, and this completes the proof. 
\qed 
\enddemo

\proclaim {Lemma 3.2} Let $G$ be a $p$-disconnected finite group, and let $\Bbb K$ be the poset (via 
inclusion) of normal subgroups $K$ of $G$ such that $p$ divides $|K/O_p(G)|$ and such that $O_p(G)\leq K$.
Then there exists a unique minimal $K\in\Bbb K$. 
\endproclaim 

\demo {Proof} We may assume that $O_p(G)=1$. Fix a Sylow $p$-subgroup $S$ of $G$ and a strongly $p$-embedded 
subgroup $H$ of $G$ containing $S$. Suppose that there exist $K,K'$ minimal in $\Bbb K$ with $K\neq K'$, and 
set $X=K\cap K'$. Then $X$ is a $p'$-group, and $S\cap KK'$ contains an elementary abelian subgroup $A$ of 
order $p^2$. By 3.1, $X$ is then generated by its subgroups $C_X(a)$ as $a$ varies over the set of 
non-identity elements of $S\cap KK'$. Thus $X\leq H$. Since 
$[K,S\cap K']\leq X$ we have $K\leq N_G(S\cap K')X$, and so $K\leq H$. Then $G=N_G(S\cap K)K\leq H$, 
which is contrary to $H$ being a proper subgroup of $G$. Thus $\Bbb K$ has a unique minimal member $K$. 
\qed 
\enddemo

In what follows we shall refer to the group $K$ in the preceding lemma as the {\it $p$-socle} of $G$.

\definition {Notation 3.3} Let $(\ca L,\D,S)$ be a proper locality on $\ca F$, and let $T\leq S$ be 
strongly closed in $\ca F$. Denote by $\bold A(\ca F)$ the set of all $P\in\ca F^{cr}$ such that: 
\roster 

\item "{(1)}" $N_S(P)\in Syl_p(N_{\ca L}(P)$, and 

\item "{(2)}" either $P=S$ or $N_{\ca L}(P)/O_p(N_{\ca L}(P))$ is $p$-disconnected. 

\endroster 
Let $\bold A_T(\ca F)$ be the set of all $P\in\bold A(\ca F)$ such that also:
\roster 

\item "{(3)}" $P\cap T$ is fully normalized in $\ca F$. 

\endroster 
\enddefinition

Note that the condition (1) is equivalent to the statement that $P$ is fully normalized in $\ca F$, 
and that condition (2) is equivalent to the statement that either $P=S$ or $Out_{\ca F}(P)$ is 
$p$-disconnected. Thus the sets $\bold A(\ca F)$ and $\bold A_T(\ca F)$ depend only on $\ca F$ and $T$.

\definition {Definition 3.4} An element $g\in\ca L$ is {\it $\bold A_T(\ca F)$-decomposable} if there 
exists $w\in\bold D$ and a sequence $\s$ of members of $\bold A_T(\ca F)$:  
$$ 
w=(g_1,\cdots g_n),\quad \s=(P_1,\cdots,P_n), 
$$
such that the following hold. 
\roster 

\item "{(1)}" $S_g=S_w$ and $g=\Pi(w)$.  

\item "{(2)}" $P_i=S_{g_i}$ for all $i$. 

\item "{(3)}" Either $P_i=S$ or $g_i\in O^p(K_i)$, where $K_i$ is the $p$-socle of $N_{\ca L}(R_i)$. 

\endroster 
We say also that $(w,\s)$ is an {\it $\bold A_T(\ca F)$-decomposition} of $g$ if (1) through (3) hold.  
\enddefinition 

Since condition (2) in 3.4 implies that the sequence $\s$ is determined by $w$, there is some redundancy 
in the definition. For that reason we shall also speak of the $\bold A_T(\ca F)$-decomposition 
$w$ and its {\it auxiliary sequence} $\s$. 
\vskip .1in 

The following result is a version of the Alperin-Goldschmidt fusion theorem [Gold]. 

\proclaim {Theorem 3.5} Let $(\ca L,\D,S)$ be a proper locality on $\ca F$, and let $T$ be 
strongly closed in $\ca F$. Then every element of $\ca L$ is $\bold A_T(\ca F)$-decomposable. 
\endproclaim 

\demo {Proof} Set $\bold A=\bold A_T(\ca F)$. The following point is then immediate from  
definition 3.4. 
\roster 

\item "{(1)}" Let $u=(x_1,\cdots,x_k)\in\bold D$ with $S_u=S_{\Pi(u)}$. Suppose that each $x_i$ is 
$\bold A$-decomposable. Then $\Pi(u)$ is $\bold A$-decomposable. 

\endroster 
Among all $g\in\ca L$ such that $g$ is not $\bold A$-decomposable, choose $g$ with $|S_g|$ as 
large as possible, and set $P=S_g$. Then $P\neq S$, as otherwise $w=(g)$ and $\s=(S)$ provide an 
$\bold A$-decomposition for $g$. Set $P'=P^g$. 

As $T$ is strongly closed in $\ca F$, there exists an $\ca F$-conjugate $Q$ of $P$ (and hence also 
of $P'$) such that both $Q$ and $Q\cap T$ are fully normalized in $\ca F$, by II.1.15. As $P,P'\in\D$,  
there exist $a,b\in\ca L$ such that $P^a=Q$, $(P')^b=Q$, and $N_S(P)^a\leq N_S(Q)\geq N_S(P')^b$. The 
maximality of $|P|$ implies that $a$ and $b$ are $\bold A$-decomposable, and the same is then true of $a\i$ 
and $b\i$ via the inverses of the words (and the reversals of the sequences of subgroups of $S$) which 
yield $\bold A$-decomposability for $a$ and $b$. Set $g'=a\i gb$ and set $M=N_{\ca L}(Q)$. Then $g'\in M$, 
$(a,g',b\i)\in\bold D$ via $Q$, and $ag'b\i=g$. If $g'$ has an $\bold A$-decomposition then so does $g$, by 
(1). Thus we may assume that $g=g'$, whence  $P=Q=P'$. 

Set $R=N_S(Q)$. Then $R\in Syl_p(M)$ by II.2.1, so $O_p(M)\leq S$, and then $O_p(M)=Q$. Let $K$ be a normal 
subgroup of $M$ which is minimal subject to $Q\leq K$ and $Q\notin Syl_p(K)$. By the Frattini lemma we may 
write $g=fh$, where $f\in O^p(K)$ and $h\in N_M(S\cap K)$. Then $h$ is  $\bold A$-decomposable, and 
$S_{(f,h)}=S_g$, so (1) implies that $f$ has no $\bold A$-decomposition. Thus we may assume that 
$g\in O^p(K)$. If $M/Q$ is $p$-disconnected then $K$ is the $p$-socle of $M$, so $Q\in\bold A$, and 
$((g),(Q))$ is a $\bold A$-decomposition of $g$. Thus $M/Q$ is $p$-connected. 

Let $\G$ be the graph $\G_p(M/Q)$, and let $R=R_0,\cdots,R_m=R^g$ be a sequence of Sylow $p$-subgroups of 
$M$ such that 
$(R_0/Q,\cdots,R_m/Q)$ is a geodesic path in $\G$ from $R/Q$ to $R^g/Q$. We may assume that, 
among all $g\in K$ having no $\bold A$-decomposition, $g$ has been chosen so that the distance $m$ from 
$R/Q$ to $R^g/Q$ in $\G$ is as small as possible. Then $m\neq 0$ as $R\nleq Q$. Suppose $m=1$. 
Then $R\cap R^g$ properly contains $Q$, and then $(R\cap R^g)^{g\i}$ is a subgroup of $S_g$ which properly 
contains $Q$, contrary to $S_g=Q$. Thus $m\geq 2$. 

Let $d\in M$ such that $(R_{m-1})^d=R_m$, and set $h=gd\i$. Then $R^h=R_{m-1}$, and the minimality of $m$ 
implies that there exists an $\bold A$-decomposition $u$ of $h$. Since $g=hd$ we have $Q=S_{(h,d)}$, and 
then (1) implies that $d$ has no $\bold A$-decomposition. Then $S_d=Q$ by the maximality of $|S_g|$ in our 
choice of $g$. Set $e=hdh\i$. Conjugation by $h$ sends 
the pair $(R,R^e)$ to $(R_{m-1},R_m)$, so $R/Q$ is adjacent to $R^e/Q$ in $\G$. As $m>1$ there then exists 
an $\bold A$-decomposition $v$ for $e$. Set $w=u\i\circ v\circ u$. Then $\Pi(w)=h\i eh=d$, and then  
$w$ is an $\bold A$-decomposition for $d$ since $S_w\geq Q=S_a$. But we have already determined that $d$ 
has no $\bold A$-decomposition, and this contradiction completes the proof. 
\qed 
\enddemo

We shall assume the setup of 2.1 until introducing a variation on that setup in 2.9. Thus, for now, 
$(\ca L,\D,S)$ is a proper locality on $\ca F$, $\ca N\norm\ca L$ is a partial normal subgroup, 
$T=\ca N\cap S$, and $\ca E=\ca F_T(\ca N)$.

\proclaim {Lemma 3.6} Let $P\in\bold A(\ca F)$, set $M=N_{\ca L}(P)$, and set $U=P\cap T$. If $P\neq S$, 
let $K$ be the $p$-socle of $M$. 
\roster 

\item "{(a)}" If $P\in\bold A_T(\ca F)$ then $U^{\ca F}\sub\ca E^c$, and $U=O_p(N_{\ca E}(U))$. 

\item "{(b)}" Either $T\leq P$ or $C_S(U)\leq P$; and if $T\nleq P$ then $K\leq N_{\ca N}(P)P$. 

\endroster 
\endproclaim 

\demo {Proof} Point (a) is a direct application of 2.6(b). If If $P=S$ then $U=T$, and then (b) holds 
vacuously. Thus we may assume that $P\neq S$, and hence that $M$ is $p$-disconnected. Set $K=N_{\ca N}(P)$. 
Then $K\norm M$, and $[P,K]\leq P\cap K\leq U$. 

Set $D=N_{C_S(U)}(P)$ and set $X=C_S(U)P$. Then $N_X(P)=DP$. Set $H=\<D^M\>$. Then $[U,H]=1$, so $H\cap K$ 
centralizes the chain $P\geq U\geq 1$. As $M$ is of characteristic $p$, II.2.7(c) yields 
$H\cap K\leq P$. Thus $HKP/P$ is the direct product of $HP/P$ with $KP/P$. As $M$ is $p$-disconnected, 
either $HP/P$ or $KP/P$ is a $p'$-group, so either $N_T(P)\leq P$ or $D\leq P$. This yields  
the first of the statements in (b). Now suppose that $T\nleq P$. Then $N_T(P)\nleq P$, and hence 
$K\leq\<N_T(P)^M\>P\leq N_{\ca N}(P)P$. Thus (b) holds. 
\qed 
\enddemo

Set $H=N_{\ca L}(C_S(T)T)$. Then $H$ is a subgroup of $\ca L$, by 2.2(a).

\definition {Definition 3.7} Denote by $\S_T(\ca F)$ the set of all products $O_p(\ca L)UV$ such that:  
\roster 

\item "{$\cdot$}" $U=(P\cap T)^h$ for some $P\in\bold A_T(\ca F)$ and some $h\in H$, and 

\item "{$\cdot$}" $V\in C_{\ca F}(T)^{cr}$. 

\endroster 
\enddefinition 

Recall from 1.5 that $(N_{\ca L}(T),\D,S)$ is a proper locality on $N_{\ca F}(T)$, and that 
$C_{\ca L}(T)\norm N_{\ca L}(T)$. We shall often write $\ca L_T$ for the locality $N_{\ca L}(T)$, and 
$\ca C_T$ for $C_{\ca L}(T)$. 

\proclaim {Lemma 3.8} Let $X\in\S_T(\ca F)$, and write $X=O_p(\ca L)UV$ as in the preceding definition. 
Then $U=X\cap T$ and $V=C_X(T)$. Moreover, $\S_T(\ca F)$ depends only on $\ca F$ and $T$, and not on 
$\ca L$. 
\endproclaim 

\demo {Proof} Notice first of all that $O_p(\ca L)T\leq O_p(\ca L_T)$. As $\ca L_T$ is a proper locality on 
$N_{\ca F}(T)$ we have $O_p(\ca L_T)=O_p(N_{\ca F}(T))$, and then 
$$ 
O_p(\ca L)T\cap\ca C_T\leq O_p(N_{\ca F}(T))\cap\ca C_T\leq O_p(\ca F_{C_S(T)}(\ca C_T))
\leq O_p(C_{\ca F}(T)),  
$$ 
since $\ca F_{C_S(T)}(\ca C_T)$ is a subsystem of $C_{\ca F}(T)$. Then 
$$ 
C_X(T)=C_{O_p(\ca L)U}(T)V\leq O_p(C_{\ca F}(T))V=V, \tag*
$$ 
since $V\in C_{\ca F}(T)^{cr}$. 

Let $P\in\bold A_T(\ca F)$ and let $h\in H$ with $U=(P\cap T)^h$. If $T\leq P$ then $U=T=X\cap T$. So 
assume that $T\nleq P$. Then $C_S(T)\leq P$ by 3.6. As $P\in\ca F^{cr}$ we have also $O_p(\ca L)\leq P$. 
Then 
$$ 
X^{h\i}=O_p(\ca L)(P\cap T)V^{h^{\i}}\leq P, 
$$ 
so $X^{h\i}\cap T=P\cap T$, and $X\cap T=U$. 

Since $O_p(\ca L)=O_p(\ca F)$, and since $H$-conjugation on subgroups of $T$ is the same as 
$Aut_{\ca F}(C_S(T)T)$-conjugation (2.2(c)), $\ca F_T(\ca F)$ depends only on $\ca F$ and $T$.  
\qed 
\enddemo

\proclaim {Lemma 3.9} Assume $\S_T(\ca F)\sub\D$. Then: 
\roster 

\item "{(a)}" $\ca F_{C_S(T)}(\ca C_T)=C_{\ca F}(T)$. 

\item "{(b)}" The image of $H$ under the natural homomorphism $H\to Aut(C_S(T))$ is contained in 
$Aut(C_{\ca F}(T))$. 

\item "{(c)}" Every member of $\S_{C_S(T)}(N_{\ca F}(T))$ contains a member of $S_T(\ca F)$. In 
particular, $\S_{C_S(T)}(N_{\ca F}(T))\sub\D$. 

\endroster 
\endproclaim 

\demo {Proof} Set $\Psi=\{Y\leq C_S(T)\mid O_p(\ca L)TY\in\D\}$. Since $\S_T(\ca F)\sub\D$ we have 
$C_{\ca F}(T)^{cr}\sub\Psi$. We may then appeal to 1.4(c), with $O_p(\ca L)T$ in the role of $V$ and with 
$\Psi$ in the role of $\S$. Thus $(\ca C_T,\Psi)$ is a localizable pair and  
$(\ca C_T)_\Psi,\Psi,C_S(T))$ is a proper locality on $\ca C_{\ca F}(T)$. This yields (a). 

Set $\ca F_T=N_{\ca F}(T)$ and let $\D_0$ be the $\ca F_T$-closure of $\S_T(\ca F)$. Let 
$R\in(\ca F_T)^{cr}$. As $\ca F_T=\ca F_S(\ca L_T)$, $\ca F_T$ is inductive by 6.1. By II.1.15 there is 
then an $\ca F_T$-conjugate $R'$ of $R$ such that $R'\cap C_S(T)$ is fully normalized in $\ca F_T$. Then 
$R'\cap C_S(T)\in C_{\ca F}(T)^{cr}$ by 2.6(b). Since $T\leq R'$, we have thus shown that 
$(\ca F_T)^{cr}\sub\D_0$. Since $\D_0\sub\D$, we have the restriction $\ca L_0$ of $\ca L$ to 
$\D_0$ (II.1.10), and plainly $H$ is a subgroup of $\ca L_0$. Then 2.3 applies with $(\ca C_T)_\Psi$ in 
the role of $\ca N$, and yields the desired action of $H$ on $C_{\ca F}(T)$. That is, (b) holds. 

Set $T^*=C_S(C_S(T))$. By definition 3.8, $\S_{C_S(T)}(\ca L_T)$ is the set of products 
$O_p(\ca L_T)UV$ such that: 
\roster 

\item "{(1)}" $U=(Q\cap C_S(T))^h$ for some $Q\in\bold A_{C_S(T)}(\ca L_T)$ and some 
$h\in N_{\ca L_T}(T^*)$, and 

\item "{(2)}" $V\in C_{\ca F_T}(C_S(T))^{cr}$. 

\endroster 
Let $h$, $U$, and $V$ be as in (1) and (2). Then $h\in H:=N_{\ca L}(T)$, and (2) yields $T\leq V$. 
Set $U'=Q\cap C_S(T)$. Then $U'\in C_{\ca F}(T)^{cr}$ by 2.6(b), and then (b) yields 
$U\in C_{\ca F}(T)^{cr}$. Then $O_p(\ca L)TU\in\S_T(\ca F)$. Since $TU\leq VU$ we then 
obtain (c). 
\qed 
\enddemo

\proclaim {Lemma 3.10} Let $U\in\ca E^c$ and let $V$ be a subgroup of $C_S(T)$, such that 
$O_p(\ca L)UV\in\D$. Then $N_{\ca N}(C_S(T)U)$ and $N_{\ca C_T}(V)$ are subgroups of the group 
$N_{\ca L}(UV)$, and 
$$
[N_{\ca N}(C_S(T)U),N_{\ca C_T}(V)]=[N_{\ca N}(UV),O^p(N_{\ca C_T}(UV)]=1. 
$$
\endproclaim 

\demo {Proof} Set $R=UV$ and set $M^*=N_{\ca L}(O_p(\ca L)R)$, $M=N_{\ca L}(R)$, $K=N_{\ca N}(R)$, 
$K_0=N_{\ca N}(C_S(T)U)$, and $X=N_{\ca C_T}(R)$. As $O_p(\ca L)R\in\D$, $M^*$ is a subgroup of $\ca L$, 
and $M=N_{M^*}(R)$ is a subgroup of $\ca L$. Moreover, $M$ and $K$ are of characteristic $p$ by II.2.6. 
  
We have $[C_S(T),K_0]\leq C_S(T)\cap\ca N=Z(T)\leq U$, as $U\in\ca E^c$. Thus $K_0\leq K$, while clearly 
$N_{\ca C_T}(V)\leq X$. As $K$ is of characteristic $p$ and $U\in\ca E^c$ we have 
$C_K(U)\leq U$. Then $[K,X]\leq U$, and so $[K,O^p(X)]=1$ by coprime action. 
\qed
\enddemo

\definition {Definition 3.11} Let $g\in\ca L$, and let $w=(a,g_1,\cdots,g_n)\in\bold D$. Then $w$ is an 
{\it $\ca N$-decomposition} of $g$ if: 
\roster 

\item "{(1)}" $S_w=S_g$, 

\item "{(2)}" $\Pi(w)=g$, 

\item "{(3)}" $a\in N_{\ca L}(T)$, and 

\item "{(4)}" there exists a sequence $(X_1,\cdots,X_n)$ of members of $\S_T(\ca F)\cap\D$ such that 
$C_S(T)\leq X_i$ and such that  
$$ 
g_i\in O^p(N_{\ca N}(X_i))\cap O^{p'}(N_{\ca N}(X_i)) \ \ (1\leq i\leq n).     
$$ 
\endroster 
\enddefinition  

A word of caution regarding the above definition: Even though $C_S(T)\leq X_i$ for all $i$, it cannot be 
concluded that $C_S(T)\leq S_w$. Indeed, there are examples to the contrary. The point is that $C_S(T)$ 
need not be invariant under $g_i$. 

\vskip .1in 
The following result may be thought of as a refinement of the splitting Lemma (I.3.12). 

\proclaim {Theorem 3.12} Assume $\S_T(\ca F)\sub\D$. Then each $g\in\ca L$ has an $\ca N$-decomposition.
Moreover, if $g\in\ca N$, and $w=(a,g_1,\cdots,g_n)$ is an $\ca N$-decomposition of $g$, then 
$a\in N_{\ca N}(C_S(T)T)$. 
\endproclaim 

\demo {Proof} Set $\w\S=\{X\in\S_T(\ca F)\mid C_S(T)\leq X\}$. Let $g\in\ca L$, and let $\Phi$ be the set 
of all words $v=(x_1,\cdots,x_n)\in\bold W(\ca L)$ such that $S_w=S_g$, $\Pi(w)=g$, and having the 
following property.  
\roster 

\item "{(*)}" For each $i$, one of the following holds. 

\itemitem {\rm{(i)}} $x_i\in H$, and $S_{x_i}\in\bold A(N_{\ca F}(T))$. 

\itemitem {\rm{(ii)}} $C_S(T)\nleq S_{x_i}\in\bold A(N_{\ca F}(T))$, and $x_i\in O^p(K_i)$ where $K_i$ is 
the $p$-socle of $N_{\ca L}(S_{x_i})$. 

\itemitem {\rm{(iii)}} $T\nleq S_{x_i}$, and $x_i\in O^p(N_{\ca N}(X_i))\cap O^{p'}(N_{\ca N}(X_i))$ 
for some $X_i\in\w\S$.  

\endroster 
We shall see, first of all, that $\Phi$ contains the set of all $\bold A_T(\ca F)$-decompositions of $g$. 
Let $u=(x_1,\cdots,x_n)$ be an arbitrary such. By 3.6, (i) or (ii) holds for any index $j$ such 
that $S_{x_j}\in\bold A(\ca L_T)$. Let then $i$ be an index such that $S_{x_i}\notin\bold A(\ca L_T)$. 
Set $x=x_i$, $P=S_{x_i}$, and let $K$ be the $p$-socle of $N_{\ca L}(P)$. Then $K=O^{p'}(K)$ by definition, 
and $K\leq N_{\ca N}(P)P$ by 3.6(b). Set $D=(P\cap T)C_S(T)$. Then $D\in\w\S$, and $D\leq P$ 
(again by 3.6(b)). Then $D\norm K$ since $[P,N_{\ca N}(P)]\leq P\cap T$. Since $x\in O^p(K)$ by 
definition 3.4, we obtain (iii), and thus $u\in\Phi$. 

Now let $v$ be an arbitrary member of $\Phi$, and suppose that there is segment $(c,b)$ of $v$ such 
that $T\nleq S_c$ and $T\leq S_b$. Let $X\in\w\S$ such that $c\in N_{\ca L}(X)$ (and such that $c$ and $X$ 
in the role of $x_i$ and $X_i$ satisfy also the stronger condition given by (iii) in (*)). 

Set $Q=S_b$, and suppose first that $C_S(T)\nleq Q$. Then $b\in O^p(L)$,  
where $L$ is the $p$-socle of $N_{\ca L}(Q)$. Set $E=N_{C_S(T)}(Q)$. Then $E\nleq Q$, so 2.2 yields 
$L\leq\<E^{N_{\ca L}(Q)}\>$, and hence $L\leq C_{\ca L}(T)$ as $T\norm L$. Let $X\in\w\S$ such that 
$c\in N_{\ca L}(X)$ (and such that $c$ and $X$ in the role of $x_i$ and $X_i$ satisfy also the stronger 
condition given by (iii)). We have $C_Q(T)\in C_{\ca F}(T)^{cr}$ by 2.6(b), so $[K,O^p(L)]=1$ by 3.10. 
Thus $c$ and $b$ commute. 
Since $c\in\ca N$ and $b\in N_{\ca L}(T)$ we have $S_{(c,b)}=S_{(b,c^b)}$ by I.3.2(a), and so 
$S_{(c,b)}=S_{(b,c)}$. Write $v=v_1\circ(c,b)\circ v_2$, and set $v'=v_1\circ(b,c)\circ v_2$.
Then $S_v=S_{v'}$ (so $v'\in\bold D$), and $\Pi(v')=\Pi(v)$. The condition (*) 
is evidently in place for $v'$, so $v'\in\Phi$. 

Suppose next that $C_S(T)\leq Q$. Then $C_S(T)T\norm N_{\ca L}(Q)$, and so $b\in H$. We have 
$X=(X\cap T)C_S(T)$ by the definition of $\w\S_T(\ca F)$, and $X\in\D$ by hypothesis. Then 
$(b,b\i,c,b)\in\bold D$ 
via $X$, and $S_{(c,b)}=S_{(b,c^b)}$ as before. Conjugation by $b$ is an isomorphism from 
$N_{\ca L}(X)$ to $N_{\ca L}(X^b)$ (I.2.3(b)), so $c^b\in O^p(N_{\ca L}(X^b))$ and 
$c^b\in O^{p'}(N_{\ca L}(X^b))$. Moreover we have $X^b\in\w\S$, since $\w\S$ is $H$-invariant. With $v$ 
defined as in the preceding case, set $v'=v_1\circ(b,c^b)\circ v_2$. We conclude that $v'\in\Phi$. 

It follows from the preceding analysis that there exists $v\in\Phi$ such that $v$ can be 
written as $v=v_1\circ v_2$, where (i) or (ii) holds for each entry of $v_1$, and where (iii) holds 
for each entry of $v_2$. As $v_1$ is a prefix of a member of $\bold D$ we have $v_1\in\bold D$. Set 
set $a=\Pi(v_1)$. Then $a\in\ca L_T$, and set $w=(a)\circ v_2$. Then $\Pi(w)=g$, and since 
$S_g\leq S_w\leq S_v=S_g$ we obtain $S_g=S_w$, and $w$ is an $\ca N$-decomposition of $g$.  

Suppose finally that $g\in\ca N$. Then $a\in\ca N$ by cancellation (I.1.4(e)), and so $a\in N_{\ca N}(T)$. 
As $\ca L$ is proper, it follows from I.3.5 that $N_{\ca N}(T)\leq H$. Thus $a\in H$, 
and the proof is complete.  
\qed 
\enddemo 

\proclaim {Corollary 3.13} Assume that $\S_T(\ca F)\sub\D$, and let $\G$ be the set of subgroups 
$U$ of $T$ of the form $(P\cap T)^h$, with $P\in\bold A_T(\ca F)$ and with $h\in H$. Then  
$\G\sub\ca E^c$, and $\ca E$ is $\G$-generated. 
\endproclaim 

\demo {Proof} The existence of an $\ca N$-decomposition for $g\in\ca N$ shows that the conjugation map 
$c_g:S_g\cap T\to (S_g\cap T)^g$ is a composition of $\ca E$-homomorphisms between members of $\G$. Thus 
$\ca E$ is $\G$-generated. That $\G$ is a subset of $\ca E^c$ is given by 2.6(b). 
\qed 
\enddemo

The following corollary to 3.12 answers a question that was left hanging from Part II. Namely, in 
Theorem A (the union of Theorems A1 and A2) one has invariance of the poset of partial normal subgroups 
of $\ca L$ under expansion of objects from $\D$ to $\D^+$, but nothing is said about what becomes of the 
fusion systems of the various partial normal subgroups of $\ca L$ under this process. One is told only (in 
the setup of Theorem A2) that if $\ca N\norm\ca L$ with $T=S\cap\ca N$, then also $T=S\cap\ca N^+$.

\proclaim {Corollary 3.14} Assume that $\S_T(\ca F)\sub\D$, and 
let $\D^+$ be an $\ca F$-closed set of subgroups of $S$ such that  
$\ca F^{cr}\sub\D\sub\D^+\sub\ca F^s$. Let $(\ca L^+,\D^+,S)$ be the unique (in the sense of Theorem 
I.A1) proper locality on $\ca F$ which contains $\ca L$ and such that 
the inclusion map $\ca L\to\ca L^+$ is a homomorphism of partial groups. Let $\ca N^+$ be the partial normal 
subgroup of $\ca L^+$ (as in Theorem II.A2) generated as a partial subgroup of $\ca L^+$ by the set of 
all $\ca L^+$-conjugates of elements of $\ca N$. Then: 
\roster 

\item "{(a)}" $\ca F_T(\ca N)=\ca F_T(\ca N^+)$, and  

\item "{(b)}" $\ca N^+$ is the partial subgroup $\<\ca N\>$ of $\ca L^+$ generated by $\ca N$. 

\endroster 
\endproclaim 

\demo {Proof} By 3.12, $\ca E$ is generated by the set of conjugation 
maps $c_g:S_g\cap T\to T$ such that $g\in\ca N$ and such that $S_g$ contains a member of $\w\S_T(\ca F)$. 
The same is true of $\ca E^+$, so (a) holds. 

Let $g\in\ca N^+$, and let $w$ be an $\ca N^+$-decomposition of $g$. Then $w$ is a sequence of elements of 
$\ca N$, and so $\ca N^+\leq\<\ca N\>$. The reverse inclusion is immediate (as $\ca N^+$ is a partial 
subgroup of $\ca L^+$). Thus (b) holds. 
\qed 
\enddemo

\proclaim {Lemma 3.15} Suppose that $O_p(\ca L)\ca E^c\sub\D$. Then:  
\roster 

\item "{(a)}" $\S_T(\ca F)\sub\D$. 

\item "{(b)}" The image of the natural homomorphism $\ca L_T\to Aut(T)$ is contained in $Aut(\ca E)$. 

\item "{(c)}" $(\ca N,\ca E^c)$ is a localizable pair, and $((\ca N)_{\ca E^c},\ca E^c,T)$ is a proper 
locality on $\ca E$. 

\endroster 
\endproclaim 

\demo {Proof} Let $P\in\bold A_T(\ca F)$, and set $U=P\cap T$. Then $U^H\sub\ca E^c$ by 2.6(b), and then  
$O_p(\ca L)U^H C_{\ca F}(T)^{cr}\sub\D$ by hypothesis. This establishes (a). Set $\G=\ca E^c$. Then 3.12 
shows that $\ca E$ is $\G$-generated.  

Let $V\in\D$ with $P\leq T$, set $M=N_{\ca L}(V)$, and set $K=N_{\ca N}(V)$. Then $K$ is of 
characteristic $p$ by II.2.6(a), and then $U\in\ca E^c$ if and only if $C_K(V)=Z(V)$. This shows that 
$\ca E^c$ is $\ca F$-invariant. 

As $\G\sub\D$, and since $T$ is maximal in the poset of $p$-subgroups of $\ca N$, it is immediate that 
$(\ca N,\G)$ is a localizable pair. Then $(\ca N_\G,\G,T)$ is a locality on a fusion subsystem $\ca E_0$ 
of $\ca E$, by 1.2. As $\ca E$ is $\G$-generated we get $\ca E_0=\ca E$, and then (c) follows from 1.3.  
Point (b) is given by 2.3. 
\qed 
\enddemo 

We end this section with an application of 3.5.

\proclaim {Theorem 3.16} Let $(\ca L,\D,S)$ be a proper locality on $\ca F$, and let $\Bbb K$ be the set of 
all partial normal subgroups $\ca K\norm\ca L$ such that $\ca L/\ca K$ is an abelian group. Set 
$[\ca L,\ca L]=\bigcap\Bbb K$. Then: 
\roster 

\item "{(a)}" $[\ca L,\ca L]=\<[N_{\ca L}(P),N_{\ca L}(P)]\mid P\in\bold A(\ca F)\>\in\Bbb K$. 

\item "{(b)}" $C_S([\ca L,\ca L])\norm\ca L$. 

\item "{(c)}" Let $\D^+$ be an $\ca F$-closed subset of $\ca F^s$ containing $\D$, and let $\ca L^+$ be 
the expansion of $\ca L$ to $\D^+$ given by Theorem II.A1, and let $[\ca L,\ca L]^+$ be the expansion 
of $[\ca L,\ca L]$ to a partial normal subgroup of $\ca L^+$ given by Theorem II.A2. Then 
$[\ca L,\ca L]^+=[\ca L^+,\ca L^+]$. 

\endroster 
\endproclaim 

\demo {Proof} For brevity, write $\ca L'$ for $[\ca L,\ca L]$, and for any group $G$ is write $G'$ for 
$[G,G]$. By I.1.9 $\ca L'$ is the union of its subsets $\ca L_k'$ $(k\geq 0)$, where 
$\ca L_0'=\bigcup\{N_{\ca L}(P)'\mid P\in\D\}$, and where $\ca L_{k+1}'$ is the set of all 
products $\Pi(w)$ with $w\in\bold W(\ca L_k')\cap\bold D$. 

Let $Y$ be the set of all $g\in\ca L$ for which there exists $x\in\ca L'$ such that $x^g$ is defined 
and $x^g\notin\ca L'$. Assume $Y\neq\nset$ and choose $g\in Y$ so that $|S_g|$ as large as possible, 
and then so that the minimal length of an $\bold A(\ca L)$-decomposition for $g$ is as small as possible.  
Suppose first that $S_g=S$, and let $k$ be the least index for which there exists $x\in\ca L_k'$ and 
$g\in N_{\ca L}(S)$ with $x^g\notin\ca L'_k$. Then $k>0$ since, by I.2.3(b), conjugation by $g$ induces 
an isomorphism $N_{\ca L}(P)\to N_{\ca L}(P^g)$ for any $P\in\D$. Thus $x=\Pi(w)$ for some 
$w\in\bold W(\ca L_{k-1}')\cap\bold D$. Set $w=(x_1,\cdots,x_n)$ and set 
$w^g=(g\i,x_1,g,\cdots,g\i,x_n,g)$. Then $w^g\in\bold D$ via $(S_w)^g$, and so 
$x^g=\Pi(w^g)=x_1^g\cdots x_n^g\in\ca L_k'$. This contradiction (to $g\in Y$) shows that $S_g\neq S$. 

Let $v$ be an $\bold A(\ca L)$-decomposition for $g$ of minimal length. Then $S_v=S_g$, and hence 
$v\i\circ(x)\circ v\in\bold D$ via $S_{(g\i,x,g)}$. The minimality of the length of $v$ then 
implies that $v=(g)$ is of length 1, and thus $g\in N_{\ca L}(Q)$ where $Q=S_g\in\bold A(\ca L)$. 
Set $H=N_{\ca L}(Q)$ and $R=N_S(Q)$. Then $g=fh$ for some $f\in H'$ and some $h\in N_H(R)$, and 
we have $Q=S_{(f,h)}$. Set $u=(h\i,f\i,x,f,h)$. Then $S_u=S_{(g\i,x,g)}\leq Q$, and 
$x^g=\Pi(u)=(x^f)^h$. Here $x^f\in\ca L'$ since $\ca L'$ is a partial group, and then 
$(x^f)^h\in\ca L'$ since the maximality of $S_g$ in the choice of $g$ yields $h\notin Y$. 
Thus $Y=\nset$, which is to say that $\ca L'\norm\ca L$. 

The set $\bar\D$ of objects of the quotient locality $\bar{\ca L}=\ca L/\ca L'$ is the set of all 
$\bar P=P/(S\cap\ca L')$ such that $S\cap\ca L'\leq P\in\D$. For any such $P$ we have 
$N_{\bar{\ca L}}(\bar P)\cong N_{\ca L}(P)/N_{\ca L'}(P)$, and thus $N_{\bar{\ca L}}(\bar P)$ is 
abelian. Set $\Theta(\bar P)=O_{p'}(N_{\bar{\ca L}}(\bar P))$. Then $\Theta(\bar P)$ centralizes 
$N_{\bar S}(\bar P)$, and a straightforward argument by induction shows that 
$\Theta(\bar P)\leq\Theta(\bar S)$ for all $\bar P\in\bar\D$. Thus $\bar{\ca L}=N_{\bar{\ca L}}(\bar S)$ 
is the group $\bar S\times\Theta(\bar S)$. Since $\bar S\norm\bar{\ca L}$, the ``correspondence  theorem" 
I.4.5 yields $\ca L'S\norm\ca L$, and then $\ca L=\ca L' N_{\ca L}(S)$ by the Frattini Lemma (I.3.11). 
Since $C_S(\ca L')\norm N_{\ca L}(S)$ and $C_S(\ca L')\norm\ca L'S$, it follows from I.3.13 that 
$C_S(\ca L')\norm\ca L$. Thus (a) holds. 

Now let $\ca N\norm\ca L$ be a partial normal subgroup of $\bar{\ca L}$ such that $\ca L/\ca N$ is an 
abelian group, and let $\r:\ca L\to\ca L/\ca N$ be the canonical projection (see I.4.4). Then 
$\r$ maps subgroups of $\ca L$ to subgroups of $\ca L/\ca N$ by I.4.2; hence $N_{\ca L}(P)'\leq\ca N$ 
for all $P\in\D$. Thus $\ca L'\leq\ca N$, completing the proof of (b). 

Let $\D^+$ and $\ca L^+$ be given as in (c). For any partial normal subgroup $\ca N\norm\ca L$ let 
$\ca N^+$ be the unique partial normal subgroup of $\ca L^+$ which intersects $\ca L$ in $\ca N$, as given 
by Theorem II.A2. In particular, taking $\ca N=(\ca L^+)'\cap\ca L$, we have $\ca N^+=(\ca L^+)'$. 
Let $\r^+:\ca L^+\to\ca L^+/(\ca L^+)'$ be the canonical projection. Then $\ca N$ is the kernel of the 
restriction $\r$ of $\r^+$ to $\ca L$. Since $Im(\r^+)$ is an abelian group, it follows that 
$\ca L'\leq\ca N$, and so $(\ca L')^+\leq(\ca L^+)'$. On the other hand, we have 
$\ca L/\ca L'\cong\ca L^+/(\ca L')^+$ by II.5.6(a). As $\ca L/\ca L'$ is an abelian group we thereby obtain 
$(\ca L^+)'\leq(\ca L')^+$, which completes the proof of (c). 
\qed 
\enddemo

\vskip .2in 
\noindent 
{\bf Section 4: $\S_T(\ca F)$} 
\vskip .1in 

The preceding section indicates that the set $\S_T(\ca F)$ introduced in 3.7 plays an important role, if 
it so happens that $\S_T(\ca F)\sub\D$. But the only indication that has been given so far, as to when this 
condition is met is the one given by lemma 3.15: in which there is the rather 
strong assumption that $O_p(\ca L)\ca E^c\sub\D$. 
One aim of this section is to show that one indeed has $\S_T(\ca F)\sub\D$ provided only that $\D$ is 
large enough. For example, it will suffice that $\D$ be as large as possible, i.e. that $\D$ be 
the set $\ca F^s$ of $\ca F$-subcentric subgroups of $S$. 

\vskip .1in
We continue the setup of 2.1. Thus $(\ca L,\D,S)$ is a proper locality on $\ca F$, $\ca N\norm\ca L$ is a 
fixed partial normal subgroup of $\ca L$, $T=S\cap\ca N$, and $\ca E=\ca F_T(\ca N)$. We also continue the 
notation: $H=N_{\ca L}(C_S(T)T)$ (a subgroup of $\ca L$ by 2.2), $\ca L_T=N_{\ca L}(T)$, and 
$\ca C_T=C_{\ca L}(T)$. Recall from 1.5 that $(\ca L_T,\D,S)$ is a proper locality on $N_{\ca F}(T)$, 
and that $\ca C_T\norm\ca L_T$.

\proclaim {Lemma 4.1} Assume that $S=C_S(T)T$, and assume that $\ca F^c\sub\D$.  
\roster 

\item "{(a)}" $UV\in\ca F^c$ for all $(U,V)\in\ca E^c\times C_{\ca F}(T)^c$. 

\item "{(b)}" Let $V$ be a subgroup of $C_S(T)$ containing $Z(T)$, and with $TV\in\ca F^c$. Then 
$V\in C_{\ca F}(T)^c$. 

\endroster 
\endproclaim 

\demo {Proof} First, let $U\leq T$ be a subgroup of $T$, and let $U'\in U^{\ca F}$. Thus there exists 
$w=(g_1,\cdots,g_n)\in\bold W(\ca L)$, and a sequence $(U_0,\cdots,U_n)$ of $\ca F$-conjugates of $U$ such 
that $U_0=U$, $U_n=U'$, and $U_k=(U_{k-1})^{g_k}$ for all $k$ with $1\leq k\leq n$. Write $g_k=x_ky_k$ where 
$x_k\in\ca N$ and $y_k\in C_{\ca L}(T)$. We have $N_{\ca N}(T)\leq N_{\ca L}(C_S(T)T)$ by I.3.5, and 
$S_{(x_k,y_k)}=S_{g_k}$ by I.3.10 and I.3.12. It follows that $U_k=(U_{k-1})^{x_k}$, and thus: 
\roster 

\item "{(1)}" $U^{\ca F}=U^{\ca E}$.  

\endroster 
Assume now that (a) is false, and among all $(U,V)\in\ca E^c\times C_{\ca F}(T)^c$ with 
$UV\notin\ca F^c$ choose $(U,V)$ so that $|U||V|$ is as large as possible. Set 
$R=UV$. Notice that $N_T(U)\leq N_S(R)$ and that $N_{C_S(T)}(V)\leq N_S(R)$. As $S=C_S(T)T$  
we have $N_S(R)=N_T(U)N_{C_S(T)}(V)$, so the maximality of $|U||V|$ yields $N_S(R)\in\ca F^c$. 

By II.1.15 there exists an $\ca F$-conjugate $R'$ of $R$ such that both $R'$ and $R'\cap T$ are fully 
normalized in $\ca F$. There is then an $\ca F$-homomorphism $\l:N_S(R)\to N_S(R')$ such that $R'=R\l$. 
As $N_S(R)\in\D$, $\l$ is given by conjugation by some $g\in\ca L$. We may employ I.3.12 in order to 
write $g=fy$, with $f\in N_{\ca L}(T)$, $y\in\ca N$, and with $S_g=S_{(f,y)}$. Set $U_1=U^f$ and $V_1=V^f$. 
Then $U_1\in U^{\ca E}$ by (1), and so $U_1\in\ca E^c$. Applying I.3.12 also to $\ca L_T$ and its partial 
normal subgroup $\ca C_T$, we obtain $f=hx$ where 
$h\in H$, $x\in\ca C_T$, and $S_f=S_{(h,x)}$. Here $H=N_{\ca L}(S)$ since $S=C_S(T)T$, so conjugation by 
$h$ is defined on all of $\ca L$, and then 2.3 shows that conjugation by $h$ preserves $C_{\ca F}(T)^c$. 
Thus $V_1\in C_{\ca F}(T)^c$, and we may therefore assume that $U=U_1$ and $V=V_1$, and that $g=y$. 
Then I.3.1(b) yields: 
\roster 

\item "{(2)}" $V^g\leq VT$. 

\endroster 
Set $U'=U^y$ and $V'=V^y$. We now compute: 
$$ 
\align 
C_S(R')&=C_S(U'V')=C_S(V')\cap C_S(U')=C_S(V')\cap C_{C_S(T)T}(U')\\ 
&=C_S(V')\cap C_S(T)Z(U')\quad\text{(as $U'$ is centric in $\ca E$)}\\ 
&=C_{C_S(T)}(V')Z(U')\quad\text{(as $[U',V']\cong[U,V]=1$)}.  \\
&=C_{C_S(T)}(V'T)Z(U')=C_{C_S(T)}(VT)Z(U') 
\quad\text{(by (2))}\\ 
&=Z(V)Z(U'). 
\endalign 
$$ 
Thus, in order to show that $R'$, and hence $R$, is centric in $\ca F$ it suffices to show that 
$Z(V)\leq R'$. Since $C_S(R)=Z(V)Z(U)$, it then suffices to show that $|Z(V)Z(U')|=|Z(V)Z(U)|$. As 
$V\leq C_S(T)$, and since $U,U'\in\ca E^c$ we have:  
$$ 
Z(V)\cap Z(U')=Z(V)\cap U'=Z(V)\cap T=Z(V)\cap U=Z(V)\cap Z(U). 
$$ 
Thus $|Z(V)\cap Z(U')|=|Z(V)\cap Z(U)|$, and hence $|Z(V)Z(U')|=|Z(V)Z(U)|$, as required. Thus $R\in\ca F^c$,  
and this contradiction completes the proof of (a). 

Now let $V$ be a subgroup of $C_S(T)$ containing $Z(T)$, and with $TV\in\ca F^c$, and suppose that 
$V\notin C_{\ca F}(T)^c$. Let $V'$ be a $C_{\ca F}(T)$-conjugate of $V$, with $V'$ fully centralized in 
$C_{\ca F}(T)$. Then $C_{C_S(T)}(V')\nleq V'$ by II.1.10. That is, we have $C_S(TV')\nleq V'$. On the 
other hand, since $TV'\cap T=T$ and $TV'\cap C_S(T)=Z(T)V'=V'$, we have 
$$ 
C_S(TV')\cap TV'=Z(TV')=Z(T)Z(V')\leq V'. 
$$ 
This shows that $C_S(TV')\nleq TV'$, and so $TV'\notin\ca F^c$. But $C_{\ca F}(T)=\ca F_{C_S(T)}(\ca C_T)$ 
by 1.5(b), and thus $TV'$ is an $\ca F$-conjugate of $TV$ via the same sequence of conjugation maps by 
elements of $\ca C_T$ that sends $V$ to $V'$. As $TV\in\ca F^c$ and $TV'\notin\ca F^c$, we have the  
contradiction which proves (b). 
\qed 
\enddemo

\proclaim {Corollary 4.2} If $S=C_S(T)T$ and $\ca F^c\sub\D$, then $\S_T(\ca F)\sub\D$. 
\endproclaim 

\demo {Proof} Let $X=O_p(\ca L)UV\in\S_T(\ca F)$, with $U$ and $V$ given as in definition 3.7. Then 
$U\in\ca E^c$ by 2.6(b), while $V\in C_{\ca F}(T)^c$ by definition. The preceding lemma then yields 
$UV\in\ca F^c$, and hence $X\in\D$. 
\qed 
\enddemo 

\proclaim {Lemma 4.3} Assume $\S_T(\ca F)\sub\D$ and assume $T\in\ca F^q$. 
\roster 

\item "{(a)}" If $\D=\ca F^s$ then $(\ca N,\ca E^c)$ is a localizable pair, and 
$(\ca N_{\ca E^c},\ca E^c,T)$ is a proper locality on $\ca E$. 

\item "{(b)}" If $\D=\ca F^s$ then the image of the natural homomorphism $\l:\ca L_T\to Aut(T)$ is 
contained in $Aut(\ca E)$. 

\item "{(c)}" $\ca E^c\sub\ca F^q$, and $\ca E^s\sub\ca F^s$. 

\item "{(d)}" $\ca E$ is $(cr)$-generated. 

\endroster  
\endproclaim 

\demo {Proof} If points (c) and (d) hold in the case that $\D=\ca F^s$ then they hold in general, by 3.14. 
Thus, we may assume throughout that $\D=\ca F^s$. 

Set $\ca U^c=\{U\in\ca E^c\mid U^{\ca F}\sub\ca E^c\}$, and set 
$\ca U^s=\{U\in\ca E^s\mid U^{\ca F}\sub\ca E^s\}$. The strategy of the proof is to first establish 
all three parts of the lemma under the following assumption.  
\roster 

\item "{(1)}" $\ca U^c=\ca E^c\sub\D$. 

\endroster 
Assume (1). Then 3.15 yields (a) and (b), and (d) then follows from (a) and from II.2.10. 
Further, (b) implies that $\ca U^s=\ca E^s$. Pick 
$U\in\ca E^s$. In order to show that $U\in\ca F^s$ we may assume that $U$ is fully normalized in $\ca F$. 
Set $X=O_p(N_{\ca E}(U))$. Then $X\in\ca E^c$, and so $X\in\ca F^s$ (and $X\in\D$) by (1). Further,  
II.1.16 implies that $X$ is fully centralized in $\ca F$. Set $M=N_{\ca L}(X)$ and $K=C_{\ca L}(X)$. 
Then $K$ is of characteristic $p$, and $X=O_p(K)$ by II.2.3. 

Set $Y=O_p(N_{\ca F}(U))$. Then $X\norm Y$, and then $Y=O_p(N_M(U))$. As $X\leq Y$ we may compute: 
$$ 
C_{\ca L}(Y)=C_M(Y)=C_{N_M(U)}(Y). 
$$ 
As $N_M(U)$ is of characteristic $p$ by II.2.6(b), we conclude that $C_{\ca L}(Y)\leq Y$. Then 
$Y\in\ca F^c$ by II.2.8(b), and $U\in\ca F^s$. Thus $\ca E^s\sub\ca F^s$. 

It now suffices to prove the following stronger version of (1).   
\roster 

\item "{(2a)}" $\ca U^c\sub\ca F^q$, and   

\item "{(2b)}" $\ca U^c=\ca E^c$.   

\endroster  

Assume that (2a) is false, and choose $U\in\ca U^c-\ca F^q$ with $|U|$ as large as possible. As $U\in\ca U^c$ 
we may assume that $U$ is fully normalized in $\ca F$. Set $\G=\{X\in\D\mid U\norm X\}$ and 
$\S=\{Y\leq C_S(U)\mid UY\in\D\}$. By 1.4 (and since $\ca F^c\sub\D$) both $(\ca N_{\ca L}(U),\G)$ and 
$(C_{\ca L}(U),\S)$ are localizable pairs. Write $\ca L_U$ for $N_{\ca L}(U)_\G$, and $\ca C_U$ for 
$C_{\ca L}(U)_\S$. Then, by 1.4, $(\ca L_U,\G,N_S(U))$ is a proper locality on $N_{\ca F}(U)$, 
$(\ca C_U,\S,C_S(U))$ is a proper locality on $C_{\ca F}(U)$, and Evidently $\ca C_U=C_{\ca L_U}(U)$, 
and so $\ca C_U\norm\ca L_U$. 

Set $B=N_T(U)$. Thus $U$ is properly contained in $B$, and so $B\in\ca F^q$.  Let $V\in C_{\ca F}(U)^{cr}$ 
with $V$ fully normalized in $C_{\ca F}(U)$, and set $D=N_{\ca C_U}(P)$. Then $V\in\D$ by 1.4(a), so $D$ is 
a subgroup of $\ca C_U$, and $D=N_{\ca C_U}(UV)$. Set $K=\ca N\cap D$. Then $B\leq K\norm D$, and then   
$$ 
[B,D]\leq C_K(U)=Z(U)\times O_{p'}(K), \tag*  
$$
since $B\in\ca E^c$. As $\ca L_U$ is a proper locality, $N_{\ca L_U}(UV)$ is of characteristic $p$, 
so also $D$ and $K$ are of characteristic $p$ by II.2.6. Thus $O_{p'}(K)=1$, and so (*) yields 
$[B,D]\leq Z(U)$. Then $[B,O^p(D)]=1$ by II.2.6(c). As $B\in\ca F^q$ we conclude that $O^p(D)=1$. 
As $V\in C_{\ca F}(V)^{cr}$ it follows that $V=C_S(U)=\ca C_U$. Then $C_{\ca F}(U)$ is the trivial 
fusion system on $C_S(U)$, and hence $U\in\ca F^q$. This completes the proof of (2a). 

Finally, assume (2b) to be false, and among all $U\in\ca E^c-\ca U^c$, choose $U$ so that $|U|$ is as large 
as possible. Then $U\neq T$, so $U$ is a proper subgroup of $N_T(U)$, and then $N_T(U)\in\ca F^q$ by (2a). 
Thus $N_T(U)\in\D$, and we may then argue - in a perhaps familiar way - as follows. Let $U'$ be a fully 
normalized $\ca F$-conjugate of $U$, and let $\phi:N_S(U)\to S$ be an $\ca F$-homomorphism such that 
$U\phi=U'$. Then $\phi$ is given by conjugation by an element $g\in\ca L$ (as $N_T(U)\in\D$), and $g=xh$ for 
some $x\in\ca N$ and $h\in\ca L_T$ with $S_g=S_{(x,h)}$ (I.2.12). As $\ca E^c$ is both $\ca E$-invariant 
and $\ca L_T$-invariant, we then have $U'\in\ca E^c$, and a similar argument then show that 
$(U')^{\ca F}\sub\ca E^c$. Thus $U\in\ca U^c$ after all. This completes the proof of (2b), and thereby  
proves the lemma. 
\qed 
\enddemo

\definition {Remark} It appears that point (c) of the preceding lemma cannot be improved upon in any 
obvious way. For example, $\ca E^{cr}$ need not be contained in $\ca F^c$, even if $T\in\ca F^c$. 
For example, let $G$ be a semi-direct product $V\rtimes H$ where $H=GL_3(2)$ and where $V$ is 
elementary abelian of order 16, with $C_V(H)=1$. Then $G$ may be viewed as a proper 
locality whose objects are the overgroups of $V$ in a Sylow $2$-subgroup $S$ of $G$. Let $N$ be the 
subgroup of $G$ of index 2 in $G$, and set $T=S\cap N$. Then $T$ is centric in $\ca F_S(G)$, and 
$V\cap T$ is centric radical in $\ca F_T(\ca N)$;  but $V\cap T$ is not centric in $\ca F_S(G)$. 
\enddefinition

Our aim now is to use 4.1 and 4.3 to show that $\S_T(\ca F)\sub\ca F^q$ if $\D=\ca F^s$.  
At the same time, we want to obtain information about the special case where $C_S(\ca N)$ is 
contained in $\ca N$. The next result (a corollary of 3.11) prepares the way for these goals.

\proclaim {Lemma 4.4} Assume $\S_T(\ca F)\sub\D$, and let $\w\S_T(\ca F)$ be the set of all 
$X\in\S_T(\ca F)$ such that $C_S(T)\leq X$. Then  
$$
C_S(\ca N)=\bigcap\{C_S(N_{\ca N}(X))\mid X\in\w\S_T(\ca F)\}. \tag*
$$ 
\endproclaim 

\demo {Proof} Let $R$ be the right-hand intersection in (*). Then $C_S(\ca N)\leq R$, so it remains to 
prove the reverse inclusion. Let $f\in\ca N$ and let $g\in R$. As $\S_T(\ca F)\sub\D$ by hypothesis, 
there exists an $\ca N$-decomposition $w=(a,f_1,\cdots,f_n)$ of $f$, and $a\in\ca N\cap H$, by 3.12. Set 
$w'=(g\i,a,g,g\i,f_1,g,\cdots,g\i,f_n,g)$. Then $w'\in\bold D$ via $(S_w)^g$. By definition 3.11, each 
$f_i$ normalizes a member of $\w\S_T(\ca F)$, so $(f_i)^g=f_i$. As $C_S(T)T\in\w\S_T(\ca F)$, also 
$a^g=a$. Thus $f^g=\Pi(w')=\Pi(w)=g$, and so $g\in C_S(\ca N)$. Thus $R\leq C_S(\ca N)$, and the reverse 
inclusion is obvious. 
\qed 
\enddemo 

\proclaim {Proposition 4.5} Assume $\S_T(\ca F)\sub\D$, and assume that $C_S(\ca N)\leq\ca N$. Set 
$\G=\ca E^s$. 
\roster 

\item "{(a)}" $\ca E^c\sub\ca F^q$, and $\S_T(\ca F)\sub\ca F^q$. 

\item "{(b)}" If $\D=\ca F^s$, then $(\ca N,\ca E^s)$ is a localizable pair, and $(\ca N_\G,\G,T)$ is a 
proper locality on $\ca E$. 

\item "{(c)}" If $\D=\ca E^s$ then $\ca L_T$ is a subgroup of $\ca L$, and the image of $\ca L_T$ in 
$Aut(T)$ is a subgroup of $Aut(\ca N_\G)$. 

\item "{(d)}" If $\D=\ca F^s$ then $O_p(\ca E)=O_p(\ca N_\G)=O_p(\ca N)\norm\ca L$. 

\endroster  
\endproclaim 

\demo {Proof} It follows from 3.14 that if (a) holds if it holds under the assumption that 
$\D=\ca F^s$. We may therefore assume throughout that $\D=\ca F^s$. 

The proof will have two parts. In the first we assume: 
\roster 

\item "{(*)}" $\S_T(\ca F)\sub\ca F^q$. 

\endroster 
Given (*), we will show that (a) through (d) hold. Once that has been achieved, we will then be able 
to show that the hypothesis (*) is redundant, and to thereby complete the proof.   

Assume (*). Define $w\S_T(\ca F)$ as in 4.4, let $X\in\w\S_T(\ca F)$, and set $U=X\cap T$. Then 
$X=O_p(\ca L)C_S(T)U$. Let $V\in C_{\ca F}(T)^{cr}$ with $V$ fully normalized in $C_{\ca F}(T)$. Then 
$O_p(\ca L)TV\in\S_T(\ca F)$, so $O_p(\ca L)TV\in\D$. Then $N_{C_{\ca F}(T)}(V)$ is the fusion system of 
$N_{\ca C_T}(TV)$ 
over $C_{C_S(T)}(V)$, by II.2.2. As $\ca L$ is proper, $N_{\ca L}(TV)$ is of characteristic $p$. 
Since $N_{\ca L}(TV)=N_{\ca L_T}(V)$, and since $\ca C_T\norm\ca L_T$, the group $K:=O^p(N_{\ca C_T}(V))$ 
is of characteristic $p$ by II.2.6(a). Set $V_0=V\cap K$. Then $V_0\leq C_S(\ca N)$, by 3.9 and 4.4. 
As $C_S(\ca N)\leq T$ by hypothesis, we then have $V_0\leq Z(K)$. Then $[V,K]=1$ by II.2.6(c), and 
then $K=1$ since $V$ is centric in $C_{\ca F}(T)$. As $V$ is also 
radical in $C_{\ca F}(T)$ we conclude that $V=C_S(T)$. Thus $C_{\ca F}(T)$ is the trivial fusion 
system on $C_S(T)$. That is, we have $T\in\ca F^q$. Point (a) is then given by 4.3(c) and by (*). Points  
(b) and (c) are given by the relevant points in 4.3, in conjunction with expansion from 
$\ca E^c$ to $\ca E^s$ via Theorem II.A. 
 
Next, it follows from (b) and II.2.3 that $O_p(\ca E)=O_p(\ca N_\G)$, and from (c) that $O_p(\ca E)$ is 
$H$-invariant. As $\ca N_\G=\ca N\cap\ca L_0$ we have $\ca N_\G\norm\ca L_0$. 
As $\ca L_0=\ca N\ca L_T$ we obtain $O_p(\ca E)\norm\ca L_0$, and then $O_p(\ca E)\norm\ca L$ by 
Theorem II.A2. This establishes (d). It now remains to remove the hypothesis (*). 

Let $\ca M$ be the partial subgroup $\<\ca C_T\ca N\>$ of 
$\ca L$. Then $\ca M\norm\ca L$ and $S\cap\ca M=C_S(T)T$, by I.5.5. Set $\w T=C_S(T)T$ and set 
$\ca D=\ca F_{\w T}(\ca M)$. Then $\w T\in\ca F^c$, so 4.3(a) implies that $(\ca M,\ca D^c)$ is a 
localizable pair and that $\ca M_{\ca D^c}$ is a proper locality on $\ca D$. 
Then $\S_T(\ca F)\sub\ca D^c$ by 4.1(a). with $\ca M$ in the role of $\ca L$. Since  
$\ca D^c\sub\ca F^q$ by 4.3(c) we thereby obtain $\S_T(\ca L)\sub\ca F^q$. Thus (*) holds, 
and the proof is complete.  
\qed 
\enddemo

\proclaim {Corollary 4.6} Assume the setup of 2.1. Then $\S_T(\ca F)\sub\ca F^q$. 
\endproclaim 

\demo {Proof} By 3.8 the definition of $\S_T(\ca F)$ depends only on $\ca F$ and on the strongly 
closed subgroup $T$, and not on $\D$. By Theorem II.A we may therefore assume that $\D=\ca F^s$.  
In this case the proof is given by repeating - verbatim - the final paragraph in the proof of 4.5. 
\qed 
\enddemo

\proclaim {Corollary 4.7} Assume the setup of 2.1, and assume that $T$ is abelian. Then $T\norm\ca N$, and 
$\ca N$ is a subgroup of $N_{\ca L}(C_S(T)T)$. Moreover, if $T\leq Z(\ca N)$ then $\ca N=T$. 
\endproclaim 

\demo {Proof} As $\S_T(\ca F)\ca F^q\sub\ca F^s$ there exists an $\ca F$-closed set $\D^+$ of subgroups 
of $S$ such that $\D\cup\S_T(\ca F)\sub\D^+\sub\ca F^s$. Let $\ca L^+$ be the expansion of $\ca L$ to 
$\D^+$ via Theorem II.A1, and let $\ca N^+$ be the partial normal subgroup of $\ca L^+$ corresponding to 
$\ca N$ via Theorem II.A2. Let $g\in\ca N$. Then $g\in\ca N^+$, so $g$ has an $\ca N^+$-decomposition 
$w=(a,g_1,\cdots,g_n)$. For each index $i$ (if any) with $1\leq i\leq n$ there exists 
$X_i\in\S_T(\ca F)$ such that $T\nleq X_i$, $N_T(X_i)\in Syl_p(N_{\ca N}(X_i))$, and 
$g_i\in O^p(N_{\ca N}(X_i))$. As $N_T(X_i)$ is abelian, these conditions imply that 
$O^p(N_{\ca N}(X_i))=1$, and so $w=(a)$. Thus $g\in N_{\ca N}(C_S(T)T$ by 3.12. 

Now $\ca N$ is a subgroup of the group $H:=N_{\ca L}(C_S(T)T)$. Assume that $T\leq Z(\ca N)$. Then 
$\ca N=T\times O_{p'}(\ca N)$, and then $O_{p'}(\ca N)=1$ since $H$ is of 
characteristic $p$. Thus $\ca N=T$ in this case.   
\qed 
\enddemo

Recall (cf. 2.7) that $Z(\ca N)$ is the set of all $z\in\ca N$ such that $g^z$ is defined and is equal to 
$f$ for all $g\in\ca N$.

\proclaim {Lemma 4.8} Suppose $\S_T(\ca F)\sub\D$. Then $C_S(\ca N)\norm H$, and 
$Z(\ca N)=C_T(\ca N)\norm\ca L$. 
\endproclaim 

\demo {Proof} Define $\w\S_T(\ca F)$ as in 4.4. Thus, each $X\in\w\S$ is of the form $O_p(\ca L)C_S(T)U$
 where $U=(P\cap T)^h$ for some $P\in\bold A_T(\ca F)$ and some $h\in H$. As $C_S(T)\norm H$, 
$\w\S_T(\ca F)$ is $H$-invariant. Since 4.4(*) can be expressed as 
$$
C_S(\ca N)=C_S(T)\cap(\bigcap\{C_S(N_{\ca N}(X))\mid X\in\w\S_T(\ca F)\}), 
$$ 
it follows that $C_S(\ca N)\norm H$. As $Z(\ca N)=C_T(\ca N)$ by 2.7, we obtain 
$$ 
Z(\ca N)=T\cap C_S(\ca N)\norm H. 
$$
As $C_T(\ca N)\norm\ca C_T$, we may apply I.3.13 with $\ca L_T$ and $\ca C_T$ in the roles of $\ca L$ 
and $\ca N$ to obtain $Z(\ca N)\norm\ca L_T$. Then apply I.3.13 to $\ca L$ and $\ca N$ to obtain 
$Z(\ca N)\norm\ca L$. 
\qed 
\enddemo

\vskip .2in 
\noindent 
{\bf Section 5. $\ca N^\perp$} 
\vskip .1in 

We continue the setup and the notation of the preceding sections. Thus $(\ca L,\D,S)$ is a proper locality 
on $\ca F$, $\ca N$ is a partial normal subgroup of $\ca L$, $T=S\cap\ca N$, and $\ca E=\ca F_T(\ca N)$. 
Further, we have the abbreviations $\ca L_T=N_{\ca L}(T)$, $\ca C_T=C_{\ca L}(T)$, and $H=N_{\ca L}(C_S(T)T)$. 
The collection $\S_T(\ca F)$ of subgroups of $S$, defined in 3.7, will continue to play a key role. 

In this section we will be taking a roundabout path to the structure of $\ca L$ and $\ca N$ by way of the 
structure of $\ca L_T$ and $\ca C_T$. Recall from 2.2 that $H$ is a subgroup of $\ca L$, and from 1.5 that 
$(\ca L_T,\D,S)$ is a proper locality on $N_{\ca F}(T)$, that $\ca C_T\norm\ca L_T$, and that 
$C_{\ca F}(T)=\ca F_{C_S(T)}(\ca C_T)$ if $\S_T(\ca F)\sub\D$.

\proclaim {Theorem 5.1} Assume $O_p(\ca L)\ca E^{cr}\sub\D$, and assume that $\ca E$ is $(cr)$-generated. 
Let $\l$ be a mapping which assigns to each $U\in\ca E^{cr}$ a group $\l(U)$, with  
$$ 
\l(U)\norm N_{\ca L}(U)\ \ \text{and}\ \ \l(U)\leq N_{\ca N}(U). 
$$ 
Set $\ca M(\l)=\<\l(U)\mid U\in\ca E^{cr}\>$. Then $\l(U)^h$ and $\l(U^h)$ are defined for all 
$h\in\ca L_T$, and $M(\l)\norm\ca L$ if the following two conditions hold for all $U\in\ca E^{cr}$.       
\roster 

\item "{(1)}" $O^p(N_{\ca N}(U))\cap O^{p'}(N_{\ca N}(U))\leq\l(U)$, and 

\item "{(2)}" $\l(U)^h=\l(U^h)$ for all $h\in\ca L_T$. 

\endroster 
\endproclaim 

\demo {Proof} Let $\G$ be the overgroup closure of $\ca E^{cr}$ in $T$. Then $O_p(\ca L)\G\sub\D$, and 
by 2.8 shows that $(\ca N,\G)$ is a localizable pair, that $(\ca N_\G,\G,T)$ is a proper locality on $\ca E$, 
$\ca L_T$ is a subgroup of $\ca L$ which acts on $\ca N_\G$ by conjugation, and the image of the natural 
homomorphism $\ca L_T\to Aut(T)$ is contained in $Aut(\ca E)$. In particular, $\l(U)^h$ and $\l(U^h)$ are 
defined for $U\in\ca E^{cr}$ and $h\in\ca L_T$. 

Let $\D_0$ be the overgroup closure of $O_p(\ca L)\G$ in $S$. Then $\ca F^{cr}\sub\D_0\sub\D$ by 2.6(b), 
the restriction $(\ca L_0,\D_0,S)$ of $\ca L$ to $\D_0$ is a proper locality on $\ca F$, and
$\ca N_\G=\ca N\cap\ca L_0$ by 2.8(a). We shall write $\ca N_0$ for $\ca N_\G$. Let $\ca M_0$ be 
the partial subgroup of $\ca L_0$ generated by $\{\l(V)\mid V\in\ca E^{cr}\}$. 

Set $\ca U=\{(P\cap T)^h\mid P\in\bold A_T(\ca L),\ h\in N_{\ca L}(C_S(T)T\}$. Then $\ca U\sub\ca E^{cr}$ 
by 2.6(b), and so $\S_T(\ca F)\sub\D_0$. Then every element of $\ca L_0$ has an $\ca N_0$-decomposition, 
by 3.12. Let $f\in\ca M_0$ and let $g\in\ca L_0$ such that $f^g$ is defined in $\ca L_0$. Let 
$u=(a,f_1,\cdots,f_m)$ be an $\ca N$-decomposition of $f$, let $v=(b,g_1,\cdots,g_n)$ be an 
$\ca N$-decomposition of $g$, and set $w=v\i\circ u\circ v$. Thus 
$$ 
w=(g_n\i,\cdots,g_1\i,b\i,a,f_1,\cdots,f_m,b,g_1,\cdots,g_n), 
$$ 
and $S_w=S_{(g\i,f,g)}\in\D_0$. By (1) and definition 3.11, there exist elements $U_i$ and $V_j$ of 
$\ca U$ (and hence of $\ca E^{cr}$) such that $f_i\in\l(U_i)$ and $g_j\in\l(V_j)$. Further, we have 
$a,b\in N_{\ca L}(T)$. As $f\in\ca M_0$, and $f=a(f_1\cdots f_m)$ 
with $f_1\cdots f_m\in\ca M_0$, the Dedekind lemma (I.1.10) implies that $a\in\ca M_0$, and so 
$a\in N_{\ca M_0}(T)$. We may employ the Frattini Calculus (I.3.4), to obtain a word 
$$ 
w'=(b\i,a,b,g_n',\cdots,g_1',f_1',\cdots,f_m',g_1,\cdots,g_n) \tag* 
$$ 
such that $S_w=S_{w'}$, $\Pi(w)=\Pi(w')$, and (by (2)) with $g_j'$ and $f_i'$ in $\ca M_0$ for all 
$i$ and $j$. 

Set $Y_0=\bigcup\{\l(V)\mid V\in\ca E^{cr}$, and for each $k>0$ set 
$$ 
Y_k=\{\Pi(\w u)\mid \w u\in\bold W(Y_{k-1})\cap\bold D(\ca L_0)\}. 
$$ 
Then $\ca M_0$ is (by I.1.9) the union of the sets $Y_k$. We now show by induction on $k$ that 
each $Y_k$ is invariant under conjugation by $\ca L_T$. By (2), $Y_0$ is $\ca L_T$-invariant. 
Suppose that $Y_{k-1}$ is $\ca L_T$-invariant, let $\w u=(y_1,\cdots,y_d)\in Y_k$, and let 
$h\in\ca L_T$. Then $S_{\w u}\cap O_p(\ca L)T\in\D_0$, so the word 
$$ 
\w u^h=(h\i,y_1,h,\cdots,h\i,y_d,h) 
$$ 
is in $\bold D(\ca L_0)$. Since each of $(y_1)^h,\cdots,(y_d)^h$ is in $Y_{k-1}$, and since 
$\Pi(\w u^h)=\Pi(\w u)^h$ (by $\bold D(\ca L_0)$-associativity), the induction is complete. Thus 
$\ca L_T$ acts on $\ca M_0$ by conjugation. 
In particular, $N_{\ca M_0}(T)$ is a normal subgroup of $\ca L_T$, and so $a^b\in\ca M_0$. 
Now (*) yields $f^g=a^bf'$ for some $f'\in\ca M_0$, and so $f^g\in\ca M_0$. Thus $\ca M_0\norm\ca L_0$. 

Set $\ca M=\ca M(\l)$, and let $\ca M^+$ be the partial subgroup of $\ca L$ generated by the 
set $\ca M^{\ca L}$ of all $\ca L$-conjugates of elements of $\ca M$. Thus: 
$$ 
\ca M_0\sub\ca M\leq\ca M^+. 
$$
Since $\D_0$ is contained in $\D$, we have $\S_T(\ca L)\sub\D$, and so each element of $\ca L$ has 
an $\ca N$-decomposition. Let $f$ now be an arbitrary element of $\ca M^+$, and let $u=(a,f_1,\cdots,f_m)$ 
be an $\ca N$-decomposition (in $\ca L$) of $f$. As we have already seen, (1) yields $f_i\in\ca M_0$ for all 
$i$. Then $f_1\cdots f_m\in\ca M$, and $a\in N_{\ca M^+}(T)$. By Theorem II.A2, we have 
$\ca M^+\norm\ca L$, and $\ca L_0\cap\ca M^+=\ca M_0$. Then 
$$ 
N_{\ca M^+}(T)=N_{\ca L}(T)\cap\ca M^+=N_{\ca L}(T)\cap\ca M_0, 
$$ 
as $N_{\ca L}(T)\leq\ca L_0$. This shows that $a\in\ca M_0$, and that $f=\Pi(u)\in\ca M$. Thus 
$\ca M=\ca M^+$, and so $\ca M\norm\ca L$. 
\qed 
\enddemo

\definition {Remark} The hypotheses of the preceding proposition are fulfilled trivially if $\ca N=\ca L$, 
and by the mapping $U\maps O^p(N_{\ca N}(U))$ (and similarly for $U\maps O^{p'}(N_{\ca N}(U))$). 
\enddefinition

By II.7.2 there is a smallest partial normal subgroup $\ca K:=O^p_{\ca L}(\ca N)$ of $\ca L$ 
such that $\ca KT=\ca N$, and a smallest partial normal subgroup $\ca K'=O^{p'}_{\ca L}(\ca N)$ of $\ca L$ 
such that $T\leq\ca K$. In the following arguments $\ca K$ will play a decisive role, and it will be 
useful to know that the defining properties of $\ca K$ are invariant under change of objects. Thus, by  
II.7.3, if $\D\sub\D^+\sub\ca F^s$ and $\D^+$ is $\ca F$-closed, then 
$$ 
O^p_{\ca L^+}(\ca N^+)=O^p_{\ca L}(\ca N)^+ 
$$ 
(and the analogous property holds for $O^{p'}_{\ca L}(\ca N)$).

\proclaim {Corollary 5.2} Assume $O_p(\ca L)\ca E^{cr}\sub\D$, and define mappings $\l$ and $\l'$ on 
$\ca E^{cr}$ by 
$$ 
\l(U)=O^p(N_{\ca N}(U))\quad\text{and}\quad \l'(U)=O^{p'}(N_{\ca N}(U)). 
$$ 
Then $O^p_{\ca L}(\ca N)=\ca M(\l)$, and $O^{p'}_{\ca L}(\ca N)=\ca M(\l')$. 
\endproclaim 

\demo {Proof} Let $\mu$ be either of the mappings $\l$ or $\l'$. Then $\mu(U)$ is a subgroup of $N_{\ca N}(U)$ 
which is normal in $N_{\ca L}(O_p(\ca L)U)$. The conditions (1) and (2) of 5.1 are immediate from 
I.2.3(b), so 5.1 yields $\ca M(\mu)\norm\ca L$. 

Let $\ca N_\G$ and $\ca L_0$ be defined as in the first two paragraphs of the proof of 5.1, and set 
$\ca M_0=\ca M(\mu)\cap\ca N_\G$. Then $\ca M_0\norm\ca L_0$, and $\ca M_0$ is the partial subgroup of 
$\ca L_0$ generated (in $\ca L_0$) by $\{\mu(U)\mid U\in\ca E^c\}$. As $\ca N_\G$ is a locality, we have 
the quotient locality $\bar{\ca N_\G}/\ca M_0$, and the canonical projection 
$\r:\ca N_\G\to\bar{\ca N_G}$. As $\ca N_\G$ is a proper locality on $\ca E$, each element $g\in\ca N_\G$ 
has an $\bold A_T(\ca E)$-decomposition $w=(g_1,\cdots,g_n)$ (cf. 3.4). By definition, either    
$g_i\in O^p(N_{\ca N}(R_i))\cap O^{p'}(N_{\ca N}(R_i))$ for some $R_i\in\ca E^{cr}$ 
(in which case $g_i\in Ker(\r)$), or $g_i\in N_{\ca N}(T)$. It follows that $\bar{\ca N_G}$ is a $p$-group 
if $\mu=\l$, and a $p'$-group if $\mu=\l'$. 

We provide the remaining details for the case $\mu=\l$. As $T\r$ is a maximal $p$-subgroup of $\bar{\ca N_G}$ 
we obtain $\ca N_\G=\ca M_0T$. This shows that $O^p_{\ca L_0}(\ca N_\G)\leq\ca M_0$. On the other hand, the 
image of $\ca N_\G$ in $\ca L_0/O^p_{\ca L_0}(\ca N_\G)$ is a $p$-group, so 
$O^p_{\ca L_0}(\ca N_\G)$ contains each of the groups $O^p(N_{\ca N}(U))$ for $U\in\G$. Thus 
$O^p_{\ca L_0}(\ca N_\G)=\ca M_0$. We now employ Theorem A in order to view $\ca L$ as an expansion 
$(\ca L_0)^+$ of $\ca L$, and to provide a correspondence between the partial normal subgroups of $\ca L_0$ 
and the partial normal subgroups of $\ca L$. Then lemma II.7.3 yields 
$$ 
O^p_{\ca L}(\ca N)=O^p_{\ca L_0}(\ca N_\G)^+=\ca M_0(\l)^+.  
$$  
Here $\ca M_0(\l)^+=\ca M(\l)$ since $\ca M_0(\l)=\ca M(\l)\cap\ca L_0$, and this completes the proof 
in the case that $\mu=\l$. The proof for $\mu=\l'$ is essentially the same. 
\qed 
\enddemo

\definition {Notation 5.3} For elements $x,y\in\ca L$ write $[x,y]=\1$ if $S_{(x,y)}\cap S_{(y,x)}\in\D$ 
and $xy=yx$. For non-empty subsets $X$ and $Y$ of $\ca L$, write $[X,Y]=\1$ if $[x,y]=\1$ for all 
$x\in X$ and all $y\in Y$. 
\enddefinition

\proclaim {Lemma 5.4} Assume $\S_T(\ca F)\sub\D$. Then the following hold. 
\roster 

\item "{(a)}" $O_p(\ca L_T)C_{\ca F}(T)^{cr}\sub\D$. In particular, the hypothesis of 5.2 is fulfilled, 
with $\ca L_T$ in the role of $\ca L$ and with $\ca C_T$ in the role of $\ca N$. 

\item "{(b)}" $O^p_{\ca L_T}(\ca C_T)=\<O^p(N_{\ca C_T}(V)\mid V\in C_{\ca F}(T)^{cr}\>$. 

\item "{(c)}" Let $X\in\S_T(\ca F)$ and let $V\in C_{\ca F}(T)^{cr}$. Then 
$[N_{\ca N}(X),O^p(N_{\ca C_T}(V))]=\1$. 

\item "{(d)}" $O^p_{\ca L_T}(\ca C_T)C_S(\ca N)\norm\ca L_T$. 

\endroster 
\endproclaim 

\demo {Proof} Point (a) follows from the observation that $O_p(\ca L)T C_{\ca F}(T)^{cr}\sub\S_T(\ca F)$   
and that $T\leq O_p(\ca L_T)$. As $C_{\ca F}(T)=\ca F_{C_S(T)}(\ca C_T)$ by 1.5, we may apply 5.2 with  
$\ca L_T$ and $\ca C_T$ in the roles of $\ca L$ and $\ca N$, and thereby obtain (b). Point (c) is given 
by 3.10. As $C_S(\ca N)\norm\ca L_T$ by 4.8, and $O^p_{\ca L_T}(\ca C_T)\norm\ca L_T$ by definition II.7.1,  
we have also point (d). 
\qed 
\enddemo

\definition {Definition 5.5} For any proper locality $(\ca L,\D,S)$ on $\ca F$, fix an expansion 
$(\ca L^s,\ca F^s,S)$ of $\ca L$ to a proper locality on $\ca F$ whose set of objects is the set 
$\ca F^s$ of $\ca F$-subcentric subgroups of $S$. We shall refer to such a locality $\ca L^s$ as a 
{\it subcentric closure} of $\ca L$. For any $\ca N\norm\ca L$, let $\ca N^s$ be the partial normal 
subgroup of $\ca L^s$ whose intersection with $\ca L$ is $\ca N$. Set 
$$ 
\ca N^\perp=O^p_{\ca L_T}(\ca C_T)C_S(\ca N^s),  
$$ 
where $T=S\cap\ca N$, $\ca L_T=N_{\ca L}(T)$, and $\ca C_T=C_{\ca L}(T)$. When no confusion is likely 
we shall write 
$$ 
T^\perp:=S\cap\ca N^\perp, 
$$ 
even though $T^\perp$ depends on $\ca N$ rather than on $T$. The uniqueness of $\ca L^s$ (up to a unique 
isomorphism which restricts to the identity map on $\ca L$), as given by Theorem II.A1, shows that 
$C_S(\ca N^s)$ does not depend on the choice of subcentric closure of $\ca L$, and thus 
$\ca N^\perp$ and $T^\perp$ are well-defined. 
\enddefinition

\proclaim {Lemma 5.6} $\D$ may be chosen so that $\S_T(\ca F)\sub\D$, and so that: 
\roster 

\item "{(*)}" For each $X\in\S_T(\ca F)$ such that $C_S(T)\leq X$, each $f\in N_{\ca N}(X)$, 
and each $w\in\bold W(\ca C_T)\cap\bold D$, we have $S_{(f)\circ w}\cap C_S(T)T\in\D$. 

\endroster 
\endproclaim 

\demo {Proof} Assume to begin with that $\D=\ca F^s$. Set $\ca M=\<\ca C_T\ca N\>$, 
$\w T=C_S(T)T$, $\ca D=\ca F_{\w T}(\ca M)$, and $\G=\ca D^c$. Then $\ca M\norm\ca L$ and $S\cap\ca M=\w T$ 
by I.5.5. As $C_S(\ca M)\leq C_S(T)\leq\ca M$ we may appeal to 4.5 with $\ca M$ in the role of $\ca N$, and 
conclude that $\G\sub\D$, $(\ca M,\G)$ is a localizable pair, $(\ca M_\G,\G,\w T)$ is a proper locality 
on $\ca D$, and that there is a natural conjugation action of $H$ on $\ca D$. Then 4.1 applies to 
$\ca M_\G$ in the role of $\ca L$, and yields $\S_T(\ca F)\sub\G$. Let $\D_0$ be the overgroup closure of 
$\G$ in $S$. As $\ca L=\ca MH$ by the 
Frattini Lemma, the action of $H$ on $\ca D$ implies that $\D_0$ is $\ca F$-closed. 

Note that, since $\ca M_\G$ is a proper locality on $\ca D$, $\ca D$ is $\G$-generated. Then by 2.8(a),  
with $\ca M$ in the role of $\ca N$, we have $\ca F^{cr}\sub\D_0$, and the restriction $\ca L_0$ of 
$\ca L$ to $\ca L_0$ is a proper locality on $\ca F$, with $\ca M_\G=\ca M\cap\ca L_0$. 
Set $\ca N_0=\ca N\cap\ca L_0$. Thus $\ca N_0\norm\ca M_\G$, and it now suffices to 
show that (*) holds with $\ca L_0$, $\ca N_0$ and $C_{\ca L_0}(T)$ in place of $\ca L$, $\ca N$ and 
$\ca C_T$. 

Let $X\in\S_T(\ca F)$ with $C_S(T)\leq X$, let each $f\in N_{\ca N_0}(X)$, and let 
$w\in\bold W(\ca C_{\ca L_0}(T))\cap\bold D(\ca L_0)$. Set $U=X\cap T$ and $V=S_w\cap C_S(T)$. We now 
apply 4.1 to the locality $\ca M_\G$ in the role of $\ca L$ (and with $\w T$ in the role of $S$). Since 
$S_w\cap\w T=TV\in\ca D^c$, and $Z(T)\leq V$, it follows from 4.1(b) that $V\in C_{\ca F}(T)^c$. Since 
$U\in\ca E^c$ by 2.6(b), and since $\ca E=\ca F_T(\ca N_0)$ by 3.14, it follows from 4.1(a) that 
$UV\in\ca D^c$. Thus $UV\in\D_0$. Since $f\i\in N_{\ca N}(X)$, and since $UV\leq UC_S(T)\leq X$, we 
have $(UV)^{f\i}\leq S$. One then observes $(f)\circ w\in\bold D(\ca L_0)$ via $(UV)^{f\i}$. Thus 
(*) holds, and the proof is complete. 
\qed 
\enddemo

\proclaim {Theorem 5.7} Let $(\ca L,\D,S)$ be a proper locality and let $\ca N\norm\ca L$ be a partial 
normal subgroup. 
\roster 

\item "{(a)}" Both $O^p_{\ca L_T}(\ca C_T)$ and $\ca N^\perp$ are partial normal subgroups of $\ca L$. 

\item "{(b)}" $S\cap\ca N^\perp=T^\perp\leq C_S(\ca N)$, and if $\S_T(\ca F)\sub\D$ then $T^\perp=C_S(\ca N)$. 

\item "{(c)}" Let $f\in\ca N$ and $g\in\ca N^\perp$, and suppose that either $(f,g)\in\bold D$ or 
$(g,f)\in\bold D$. Then $[f,g]=\1$. 

\item "{(d)}" Let $(\ca L^+,\D^+,S)$ be an expansion of $\ca L$, and for any $\ca K\norm\ca L$ let 
$\ca K^+$ be the unique partial normal subgroup of $\ca L^+$ (given by Theorem II.A2) whose intersection 
with $\ca L$ is equal to $\ca K$. Then $(\ca N^+)^\perp=(\ca N^\perp)^+$. 

\endroster 
\endproclaim 

\demo {Proof} By 5.6 there is a choice of $\D$ such that:  
\roster 

\item "{(*)}" $\S_T(\ca F)\sub\D$ and, 

\item "{(**)}" For each $X\in\S_T(\ca F)$ such that $C_S(T)\leq X$, each $f\in N_{\ca N}(X)$, 
and each $w\in\bold W(\ca C_T)\cap\bold D$, we have $S_{(f)\circ w}\in\D$. 

\endroster 
Set $\ca K=O^p_{\ca L_T}(\ca C_T)$, and assume (*). Then, by 5.4(b), $\ca K$ is generated (as a partial 
subgroup of $\ca L$) by the union $Y_0$ of the groups $O^p(N_{\ca C_T}(P))$ taken over all 
$P\in C_{\ca F}(T)^{cr}$. Recursively define $Y_n$ for $n>0$ to be the set of all $\Pi(w)$ with 
$w\in\bold W(Y_{n-1})\cap\bold D$. Then $\ca K$ is the union of the sets $Y_n$, by I.1.9. Now let 
$X\in\S_T(\ca F)$ with $C_S(T)\leq X$. Then $[N_{\ca N}(X),Y_0]=\1$ by 5.4(c). Let $n$ be any index such 
that $[N_{\ca N}(X),Y_n]=\1$, let $f\in N_{\ca N}(X)$, and let $w=(g_1,\cdots,g_k)\in\bold W(Y_n)\cap\bold D$. 
Assume now that (**) holds, so that $(f)\circ w\in\bold D$ via some $Q\leq C_S(T)T$. Set 
$$ 
w'=(f,g_1,f\i,f,g_2,\cdots,f\i,f,g_k,f\i,f). 
$$ 
Then $w'\in\bold D$ via $Q$. Since $g_i^{f\i}=g_i$ for all $i$ we obtain  
$$ 
fg=\Pi(f\circ w)=\Pi(w')=\Pi(w\circ f)=gf, 
$$ 
by $\bold D$-associativity (I.1.4(b)). Thus $[N_{\ca N}(X),Y_{n+1}]=\1$, and induction yields the 
following result.  
\roster 

\item "{(1)}" If (*) and (**) hold then $[N_{\ca N}(X),\ca K]=\1$ for all $X\in\S_T(\ca F)$ such 
that $C_S(T)\leq X$.  

\endroster 
Continue to assume (*) and (**), let $f$ be an arbitrary element of $\ca N$, and let $u=(a,f_1,\cdots,f_m)$ 
be an $\ca N$-decomposition of $f$. Then $a\in N_{\ca N}(\w T)$ by 3.12, and so $a=tf_0$ where $t\in T$ and 
where $f_0\in O^p(N_{\ca N}(\w T))$. Set $u'=(t,f_0,\cdots,f_m)$, let $s\in S\cap\ca K$, and set  
$$ 
u''=(s,t,s\i,s,f_0,s\i,\cdots,s,f_m,s\i). 
$$ 
Then $u''\circ(s)\in\bold D$ via $(S_f)^{s\i}$. Since $s\in\ca K\leq\ca C_T$, and since $(f_i)^s=f_i$ 
for all $i$ with $0\leq i\leq m$ by (1), we obtain
$$ 
fs=\Pi(u\circ(s))=\Pi(u''\circ(s))=\Pi((s)\circ u)=sf. 
$$ 
Thus $S\cap\ca K\leq C_S(\ca N)$. By definition 5.5 we then have 
$S\cap\ca N^\perp=T^\perp=(S\cap\ca K)C_S(\ca N^s)$. Here  
$C_S(\ca N)=C_S(\ca N^s)$ by (*) and 4.4, so we have shown:  
\roster 

\item "{(2)}" If (*) and (**) hold then $T^\perp=C_S(\ca N)$. 

\endroster 
Let $(\ca L^+,\D^+,S)$ be an expansion of $\ca L$ and let $\ca N^+$ be the unique partial normal 
subgroup of $\ca L^+$ whose intersection with $\ca L$ is equal to $\ca N$. Then $C_S(\ca N)=C_S(\ca N^+)$ 
by 4.4. Thus $T^\perp=C_S(\ca N^+)$. Set $\ca L_T^+=N_{\ca L^+}(T)$. Then $(\ca L_T^+,\D^+,S)$ is a proper 
locality on $N_{\ca F}(T)$ by 1.5(a), so $(\ca L_T^+,\D^+,S)$ is an expansion of the locality 
$(\ca L_T,\D,S)$, and $(\ca N^\perp)^+\norm\ca L_T^+$. Observe that  
$$ 
(\ca N^\perp)^+\ca L_T=(\ca N^\perp)^+\ca L=\ca N^\perp. 
$$ 
Thus $(\ca N^\perp)^+$ is the partial normal subgroup of the locality $\ca L_T^+$ which corresponds to 
$\ca N^\perp$ via Theorem II.A2. Recall that we have $\ca N^\perp=\ca KT^\perp$, where 
$\ca K=O^p_{\ca L_T}(\ca C_T)T^\perp$. Let $\ca K^+$ be the partial normal subgroup of $\ca L_T$ 
corresponding to $\ca K$. Then $\ca K^+=O^p_{\ca L_T^+}(C_{\ca L^+}(T))$ by II.7.3, and then 
$\ca K^+T^\perp\norm\ca L_T^+$ by 5.4(d). Then, since $(\ca K^+T^\perp)\cap\ca L=\ca N^\perp$, 
we conclude that $(\ca N^\perp)^+=\ca K^+T^\perp$. Then $(\ca N^\perp)^+=(\ca N)^+)^\perp$ by  
definition 5.5.  

Now let $(\ca L',\D',S)$ be the restriction of $\ca L^+$ to a locality (necessarily proper, by II.2.11) on 
an $\ca F$-closed subset $\D'$ of $\ca F^s$ containing $\ca F^{cr}$. For any partial subgroup $\ca H$ of 
$\ca L^+$ set $\ca H'=\ca H\cap\ca L'$. By a straightforward exercise with definition II.7.1 one has 
$\ca K'=O^p_{(\ca L^+_T)'}(C_{\ca L^+}(T)')$, and hence $(\ca N^\perp)'=(\ca N')^\perp$. Since any 
proper locality on $\ca F$ can be obtained (up to isomorphism) from $\ca L$ by a procedure of first 
expanding and then restricting, by II.A1, we thereby obtain (d). By 4.4 and (2), we obtain also point (b). 
In particular (and without recourse to (*) or (**)): 
\roster 

\item "{(3)}" $T^\perp\leq C_S(\ca N)$. 

\endroster 
Assume that (*) holds. Then $\S_{C_S(T)}(\ca C_T)\sub\D$ by 3.9(c), and so every element of $\ca L_T$ 
has a $\ca C_T$-decomposition by 3.12. Let $g\in\ca N^\perp$, and let $v=(b,g_1,\cdots,g_n)$ be a 
$\ca C_T$-decomposition of $g$. Write $b=sg_0$ where $s\in C_S(T)$ and $g_0\in O^p(N_{\ca K}(T^\perp))$. 
Set $v'=(s,g_0,\cdots,g_n)$. Then $S_v=S_{v'}$ and $g=\Pi(v')$. As $g\in\ca N^\perp$ and  
$g_0\cdots g_n\in\ca K\leq\ca N^\perp$, it follows that $s\in S\cap\ca N^\perp$. Thus $s\in C_S(\ca N)$ 
by (3). Let $f\in\ca N$, and let $u=(f_1,\cdots,f_m)$ by an $\ca N$-decomposition of $f$. 
Suppose that $(f,g)\in\bold D$. Then $u\circ v'\in\bold D$, and $fg=\Pi(u\circ v')$.  
Since $s$ commutes with each $f_i$ we have 
$$ 
fg=s\Pi(f_1,\cdots,f_m,g_0,\cdots,g_n). 
$$ 
Since each $[f_i,g_j]=\1$ for all $i$ and $j$, by 5.4(c), we have $S_{(f_i,g_j)}=S_{(g_j,f_i)}$ and 
$P^{f_ig_j}=P^{g_jf_i}$ for all $P\in\D$ with $P\leq S_{(f_i,g_j)}$. From this fact, a straightforward 
argument by induction on $m+n$ yields: 
$$ 
fg=s\Pi(g_1,\cdots,g_n,f_0,\cdots,f_m), 
$$ 
and then $fg=gf$. Evidently the assumption $(f,g)\in\bold D$ can be replaced by the symmetric assumption 
$(g,f)\in\bold D$, so: 
\roster 

\item "{(4)}" If $\S_T(\ca F)\sub\D$ then (c) holds. 

\endroster 
Since $\ca L=\ca N\ca L_T$, and since $\ca N^\perp\norm\ca L_T$, one may employ (4) and the splitting 
lemma (I.3.12) to obtain $\ca N^\perp\norm\ca L$ in the case that $\S_T(\ca F)\sub\D$. But then 
$\ca N^\perp\norm\ca L$ in general, by (d). Similarly, since $\ca K\leq N^\perp$ and $\ca K\norm\ca L_T$, 
we obtain $\ca K\norm\ca L$. Thus (a) holds, and the proof is complete. 
\qed 
\enddemo

\vskip .2in 
\noindent 
{\bf Section 6: $\del(\ca F)$, $\d(\ca F)$, and $F^*(\ca L)$} 
\vskip .1in 

As always, $(\ca L,\D,S)$ is a proper locality on $\ca F$. We are now in position to consider all of the 
partial normal subgroups $\ca N$ of $\ca L$ simultaneously. And as always, there is the difficulty in the 
back-ground, that the fusion system $\ca E$ of $\ca N$ need not be stable with respect to the processes 
of expansion (given by Theorem II.A) and restriction (given by II.2.11). Since $S\cap\ca N$ is stable with 
respect to these two processes, one step towards addressing the difficulty is with a simple notational 
device, as follows. 

For $\ca N\norm\ca L$, and $T=S\cap\ca N$, define $\ca E^s(\ca N)$ to be the fusion system 
on $T$ of the form $\ca F_T(\ca N^s)$; where $\ca N^s$ is the partial normal subgroup of an expansion 
$(\ca L^s,\ca F^s,S)$ of $\ca L$ (as in Theorem II.A1) whose intersection with $\ca L$ is 
$\ca N$ (as in Theorem II.A2).By 3.14, $\ca E^s(\ca N)$ is stable with respect to 
restriction from $\ca F^s$ to $\ca F$-closed sets $\D$ containing $\S_T(\ca F)$.

\definition {Definition 6.1} A partial normal subgroup $\ca M\norm\ca L$ will be said to be {\it large} 
in $\ca L$ if $S\cap\ca M^\perp\leq\ca M$. Let $\Bbb M:=\Bbb M(\ca L)$ be the set of all large partial normal 
subgroups of $\ca L$. Define $\del(\ca F)$ to be the overgroup closure in $S$ of 
$\bigcup\{O_p(\ca L)\ca E^s(\ca M)^{cr}\mid \ca M\in\Bbb M\}$. 
\enddefinition

\proclaim {Lemma 6.2} 
\roster 

\item "{(a)}" $\del(\ca F)$ depends only on $\ca F$ (and not on the choice of 
a proper locality on $\ca F$). 

\item "{(b)}" $\del(\ca F)$ is $\ca F$-closed, and $\ca F^{cr}\sub\del(\ca F)\sub\ca F^s$. 

\item "{(c)}" Let $\ca M\in\Bbb M$ and set $\ca D=\ca E^s(\ca M)$. Then $\ca D$ is $(cr)$-generated,  
and $\ca D^{cr}\sub\ca F^q$. 

\endroster 
\endproclaim 

\demo {Proof} By Theorem II.A1, any proper locality $\ca L'$ on $\ca F$ is isomorphic to a ``version" of 
$\ca L$ obtained by first expanding $\D$ to $\ca F^s$, and then restricting to an $\ca F$-closed subset 
$\D'$ of $\ca F^s$ containing $\ca F^{cr}$. By Theorem II.A2 there is a canonical isomorphism 
$\ca N\maps\ca N'$ from the poset of partial normal subgroups of $\ca L$ to that of $\ca L'$, and 
5.5(c) implies that $(\ca N^\perp)'=(\ca N')^\perp$. Then $S\cap\ca N^\perp=S\cap(\ca N')^\perp$ (by 
II.A2). As $O_p(\ca L)=O_p(\ca L')$ by II.2.3, we obtain (a). 

Let $\ca M\in\Bbb M$, set $R=S\cap\ca M$, and set $\ca D=\ca E^s(\ca M)$. Let $\ca L^+$ be the expansion 
of $\ca L$ to $\D^+:=\ca F^s$, and let $\ca M^+$ be the partial normal subgroup of $\ca L^+$ whose 
intersection with $\ca L$ is $\ca M$. Then $\ca D=\ca F_R(\ca M^+)$ by definition. We have  
$$ 
C_S(\ca M^+)=S\cap(\ca M^+)^\perp=S\cap(\ca M^\perp)+=S\cap\ca M 
$$ 
by 5.5, so $C_S(\ca M^+)\leq\ca M^+$. 

Let $\G$ be the overgroup closure of $\ca D^{cr}$ in $R$, and let $\D_0$ be the overgroup closure of 
$O_p(\ca L)\G$ in $S$. Then $\G\sub\ca F^q$ by 4.3(c), so $\D_0\sub\D^+$, and thus $\del(\ca F)\sub\ca F^s$. 
We have $\ca F^{cr}\sub\del(\ca F)$ as a consequence of 2.6(b). As $R\in\G\sub\ca F^q$, $\ca D$ is 
$(cr)$-generated by 4.3(d), and thus (c) holds. Then 2.8 applies with $\ca M$ in the role of $\ca N$, and 
2.8(c) shows that $\G$ is invariant under the action of $N_{\ca L}(R)$. By the splitting lemma (I.3.12) $\G$ 
is then $\ca F$-invariant, so $\D_0$ is $\ca F$-invariant. As $\del(\ca F)$ is the union of the various sets 
$\D_0$ taken over all $\ca M\in\Bbb M$, the proof of (b) is complete. 
\qed 
\enddemo

\proclaim {Lemma 6.3} Let $\ca N\norm\ca L$ be a partial normal subgroup of $\ca L$, set $T=S\cap\ca N$,  
and $\ca E=\ca F_T(\ca N)$. Further, set $T^\perp=S\cap\ca N^\perp$ and 
$\ca E^\perp=\ca F_{T^\perp}(\ca N^\perp)$. Then $\S_T(\ca F)\sub\del(\ca F)$, and if 
$\del(\ca F)\sub\D$ then:
\roster 

\item "{(a)}" $O_p(\ca L)\ca E^{cr}(\ca E^\perp)^{cr}\sub\del(\ca F)$. 

\item "{(b)}" $\ca E$ is $(cr)$-generated. 

\endroster 
\endproclaim 

\demo {Proof} We begin by proving (a) and (b) under the assumption that $\S_{S\cap\ca K}(\ca F)\sub\D$ for 
all $\ca K\norm\ca L$. Thus, set $\ca M=\ca N\ca N^\perp$ and $R=TT^\perp$. As $\S_T(\ca F)\sub\D$ by 
assumption, 5.7(b) yields $T^\perp=C_S(\ca N)$, and then I.5.1 yields $\ca M\norm\ca L$ and $R=S\cap\ca M$. 
Set $\ca D=\ca F_R(\ca M)$. As $\S_R(\ca M)\sub\D$ we have $\ca D=\ca E^s(\ca M)$ by 3.14. Note that 
$$ 
S\cap\ca M^\perp=C_S(\ca M)\leq C_S(\ca N)=T^\perp\leq\ca M, 
$$ 
so that $\ca M\in\Bbb M$. Then $\ca D$ is $(cr)$-generated and $\ca D^{cr}\sub\D$ by 6.2(c). 
We may now appeal to 2.10, with $\ca N_1=\ca N$ and $\ca N_2=\ca N^\perp$, to obtain 
$$ 
\ca D^{cr}=\ca E^{cr}(\ca E^\perp)^{cr} 
$$ 
via 2.10(d). This yields (a). By 2.10(c) $\ca E$ is the fusion system of a proper locality, and then (b) 
follows from II.2.10. 

Assume now that $\D=\ca F^s$. Then $\S_{S\cap\ca K}(\ca F)\sub\D$ for all $\ca K\norm\ca L$, and we may 
may make further use of 2.10. Namely, 2.10(c) yields: 
\roster 

\item "{(*)}" Both $\ca E^{cr}$ and $(\ca E^\perp)^{cr}$ are $N_{\ca L}(R)$-invariant. 

\endroster 
Let $V\in C_{\ca F}(T)^{cr}$. By II.1.15 there exists an $N_{\ca F}(T)$-conjugate $V'$ of $V$ such that 
$X':=V'\cap T^\perp$ is fully normalized in $N_{\ca F}(T)$. Then $X'$ is fully normalized in $C_{\ca F}(T)$ 
by II.1.17. Set $\ca L_T=N_{\ca L}(T)$ and   $\ca C_T=C_{\ca L}(T)$. 

Note that $C_{\ca F}(T)^{cr}$ is $N_{\ca L}(R)$-invariant, by an application of 2.10(c) to 
the partial normal subgroup $T\ca C_T$ of the locality $\ca L_T$. Thus $V'\in C_{\ca F}(T)^{cr}$. 
As $C_{\ca F}(T)=\ca F_{C_S(T)}(C_{\ca L}(T))$ by 1.4, and since $C_{\ca F}(T)$ is 
inductive by II.6.1, the hypothesis of 2.6(c) is fulfilled with $(\ca L_T,\ca C_T,\ca N^\perp)$ in the 
place of $(\ca L,\ca M,\ca N)$. We therefore conclude that $X'\in(\ca E^\perp)^{cr}$. 
Set $X=V\cap T^\perp$. Then $X$ is an $N_{\ca F}(T)$-conjugate of $X'$. As $N_{\ca F}(T)=\ca F_S(\ca L_T)$, 
and $\ca L_T=\ca E^\perp N_{\ca L_T}(R)$ by the Frattini Lemma, it follows from (*) that 
$X\in(\ca E^\perp)^{cr}$. Thus: 
\roster 

\item "{(**)}" $V\cap T^\perp\in(\ca E^\perp)^{cr}$ for all $V\in C_{\ca F}(T)^{cr}$. 

\endroster 
Now let $Y\in\S_T(\ca F)$, and write $Y=O_p(\ca L)UV$ as in 3.7. Notice that since 
$T^\perp=S\cap\ca E^\perp$ is weakly closed in $\ca F$ we have $N_{\ca L}(C_S(T)T)\leq N_{\ca L}(R)$. 
As $\ca E^{cr}$ is $N_{\ca L}(R)$-invariant by (*), we obtain $U\in\ca E^{cr}$ from 2.6(b). Then 
$U(V\cap T^\perp)\in\ca D^{cr}$ by 2.10(d), and thus $Y\in\del(\ca F)$. Having thus shown that 
$\S_T(\ca F)\sub\del(\ca F)$, the proof of the lemma is complete.  
\qed 
\enddemo 

Our next aim is to show that the intersection of the set of large partial normal subgroups of $\ca L$ 
containing $O_p(\ca L)$ is itself large. Lemmas 6.4 through 6.7 will achieve that result. 

Recall that for any partial subgroup $\ca H\leq\ca L$, $Z(\ca H)$ is defined to be the set of all 
$z\in\ca H$ such that, for all $h\in\ca H$, $h^z$ is defined and is equal to $h$.

\proclaim {Lemma 6.4} Assume $\del(\ca F)\sub\D$, let $\ca N$ be a partial normal subgroup of $\ca L$, and 
set $T=S\cap\ca N$. Then $T^\perp=C_S(\ca N)$, and $Z(\ca N)=C_T(\ca N)$. Moreover, the following are 
equivalent. 
\roster 

\item "{(1)}" $T^\perp=Z(\ca N)$. 

\item "{(2)}" $T^\perp\leq\ca N$. 

\item "{(3)}" $\ca N^\perp=Z(\ca N)$. 

\endroster 
\endproclaim 

\demo {Proof} We have $\S_T(\ca F)\sub\D$ by 6.3. Then $T^\perp=C_S(\ca N)$ by definition, 
and $Z(\ca N)=C_T(\ca N)$ by 4.8. The implication (3)$\implies$(1) is then immediate. Since   
the implication (1)$\implies$(2) is trivial, it remains only to show that (2)$\implies$(3). 

Assume (2). Then $T^\perp=Z(\ca N)$, and so $\ca E^\perp$ is the trivial fusion system on 
an abelian $p$-group. As $\ca E^\perp$ is the fusion system of $\ca N^\perp$, where 
$\ca N^\perp=O^p_{\ca L_T}(\ca C_T)T^\perp$, it follows that $\ca N^\perp=T^\perp$, and thus (3) holds.  
\qed 
\enddemo 

Next, recall from II.7.2 that for any $\ca N\norm\ca L$, $O^p_{\ca L}(\ca N)$ is the smallest partial normal 
subgroup $\ca K\norm\ca L$ such that $\ca N=(S\cap\ca N)\ca K$. 

\proclaim {Lemma 6.5} Assume $\del(\ca F)\sub\D$, and let $\ca M$ and $\ca N$ be partial normal subgroups of 
$\ca L$. Suppose that $\ca M\cap\ca N\leq Z(\ca N)$. Then $O^p_{\ca L}(\ca N)\leq\ca M^\perp$. 
\endproclaim 

\demo {Proof} Set $R=S\cap\ca M$, $T=S\cap\ca N$, and $\ca E=\ca F_T(\ca N)$. Then $[R,T]\leq\ca M\cap\ca N$, 
so $[R,T]\leq Z(\ca N)$ by hypothesis. Let $U\in\ca E^{cr}$. Then $R$ normalizes $U$, and then $R$ normalizes 
$UT^\perp$ since $T^\perp:=S\cap\ca N^\perp\norm S$. Moreover, we have $UT^\perp\in\del(\ca F)$ by 6.3(a), 
so $UT^\perp\in\D$, and thus both $R$ and $N_{\ca N}(U)$ are subgroups of the group $N_{\ca L}(UT^\perp)$. 
As $Z(\ca N)\leq O_p(\ca N)$ by 2.7, we obtain 
$$ 
[R,N_{\ca N}(U)]\leq Z(\ca N)\leq U, 
$$ 
and hence $[R,O^p(N_{\ca N}(U))]=\1$ by II.2.7(c). 

Set $\ca L_R=N_{\ca L}(R)$ and $\ca C_R=C_{\ca L}(R)$. Then $\ca C_R=O^p_{\ca L_R}(\ca C_R)C_S(R)$.  
The canonical projection $\ca L_R\to\ca L_R/O^p_{\ca L_R}(\ca C_R)$ maps $O^p(N_{\ca N}(U))$ to a 
$p$-group, and so $O^p(N_{\ca N}(U))\leq O^p_{\ca L_R}(\ca C_R)$. Thus $O^p(N_{\ca N}(U))\leq\ca M^\perp$ 
for all $U\in\ca F^{cr}$, and the lemma now follows from 5.2. 
\qed 
\enddemo 

\proclaim {Lemma 6.6} Assume $\del(\ca F)\sub\D$, and let $\ca M$ and $\ca N$ be partial normal subgroups of 
$\ca L$. Assume that $C_S(\ca M)\leq\ca M$ and that $\ca M\cap\ca N^\perp\leq Z(\ca N)$. Then 
$\ca N^\perp=C_S(\ca N)$. 
\endproclaim 

\demo {Proof} Since $Z(\ca N)\leq Z(\ca N^\perp)$, the hypothesis that $\ca M\cap\ca N^\perp$ be 
contained in $Z(\ca N)$ enables an application of the preceding lemma with $\ca N^\perp$ in the role of 
$\ca N$. Thus $O^p_{\ca L}(\ca N^\perp)\leq\ca M^\perp$. Since $C_S(\ca M)\leq\ca M$, 6.4 yields 
$\ca M^\perp=Z(\ca M)=C_S(\ca M)$, and thus $O^p_{\ca L}(\ca N^\perp)$ is a $p$-group. Since 
$S\cap\ca N^\perp=C_S(\ca N)$ by 6.6, the lemma follows. 
\qed 
\enddemo

\proclaim {Lemma 6.7} Assume $\del(\ca F)\sub\D$, and let $\ca M$ and $\ca N$ be partial normal subgroups 
of $\ca L$. Assume that $O_p(\ca L)C_S(\ca M)\leq\ca M$ and that $O_p(\ca L)C_S(\ca N)\leq\ca N$. Then 
$O_p(\ca L)C_S(\ca M\cap\ca N)\leq\ca M\cap\ca N$. 
\endproclaim 

\demo {Proof} Set $\ca K=\ca M\cap\ca N$. Obviously $O_p(\ca L)\leq\ca K$, so the problem is to show 
that $C_S(\ca K)\leq\ca K$. Observe: 
$$ 
\ca M\cap(\ca N\cap\ca K^\perp)=\ca K\cap\ca K^\perp=Z(\ca K)\leq Z(\ca K^\perp). 
$$ 
Since $Z(\ca K)\leq\ca N$ we then have $Z(\ca K)\leq Z(\ca N\cap\ca K^\perp)$, and we may apply 6.5 with 
$\ca N\cap\ca K^\perp$ in the role of $\ca N$. Thus $O^p_{\ca L}(\ca N\cap\ca K^\perp)\leq\ca M^\perp$. 
As $C_S(\ca M)\leq\ca M$, $\ca M^\perp=C_S(\ca M)$ by 6.4.  
Thus $\ca N\cap\ca K^\perp$ is a $p$-group, and indeed a normal $p$-subgroup of $\ca L$. As 
$O_p(\ca L)\leq\ca K$, we obtain $\ca N\cap\ca K^\perp\leq\ca K$, and so $\ca N\cap\ca K^\perp\leq Z(\ca K)$. 
Now the hypothesis of 6.6 is satisfied, with $\ca N$ in the role of $\ca M$, and with $\ca K$ in the 
role of $\ca N$. We conclude that $\ca K^\perp=C_S(\ca K)$. As $\ca K^\perp\norm\ca L$, and 
$O_p(\ca L)\leq\ca K$, we obtain $C_S(\ca K)\leq\ca K$, as required. 
\qed 
\enddemo

\definition {Definition 6.8} If $\del(\ca F)\sub\D$, define $F^*(\ca L)$ to be the intersection of the large 
partial normal subgroups of $\ca L$ containing $O_p(\ca L)$: 
$$ 
\ca F^*(\ca L)=\bigcap\{\ca N\mid O_p(\ca L)C_S(\ca N)\leq\ca N\norm\ca L\}.  
$$ 
More generally, define $F^*(\ca L)$ to be $F^*(\ca L^+)\cap\ca L$, where $(\ca L^+,\D^+,S)$ is 
the expansion of $\ca L$ to a proper locality on $\ca F$ with $\D^+=\D\cup\del(\ca F)$. 
\enddefinition

\proclaim {Corollary 6.9} If $\del(\ca F)\sub\D$ then $F^*(\ca L)$ is a large partial normal subgroup of 
$\ca L$ containing $O_p(\ca L)$, and is the unique smallest such. Moreover, and in general $($i.e. 
whether or not $\del(\ca F)\sub\D)$, we have $F^*(\ca L^+)=F^*(\ca L)^+$ for any expansion 
$(\ca L^+,\D^+,S)$ of $\ca L$ to a proper locality on $\ca F$. 
\qed 
\endproclaim 

\demo {Proof} The first assertion is immediate from 6.7. Now drop the assumption that $\del(\ca F)$ is 
contained in $\D$, and let $(\ca L^+,\D^+,S)$ be an expansion of $\ca L$ to some $\D^+$ with 
$\ca D\cup\del(\ca F)\sub\D^+\sub\ca F^s$. Let $(\ca L^s,\ca F^s,\S)$ be the expansion of $\ca L$ to 
$\ca F^s$. Let $\bold X^+$ be the set of large 
partial normal subgroups $\ca N^+$ of $\ca L^+$ containing $O_p(\ca F)$, and 
similarly define $\bold X^s$ relative to $\ca L^s$. 

By 6.3 we have $\S_T(\ca F)\sub\del(\ca F)$ for all $T=S\cap\ca N$ with $\ca N\norm\ca L$. Then by 3.14, 
the correspondence (Theorem II.A2) between the set of partial normal subgroups of $\ca L^+$ and the 
set of partial normal subgroups of $\ca L^s$ restricts to a bijection $\bold X^+\to\bold X^s$. Then 
$$ 
F^*(\ca L^+)=\bigcap\bold X^+=\ca L^+\cap(\bigcap\bold X^s)=\ca L^+\cap F^*(\ca L^s).   
$$ 
This shows that $F^*(\ca L^+)$ is determined by $F^*(\ca L^s)$. Then 
$$ 
F^*(\ca L)=\ca L\cap F^*(\ca L^s)=\ca L\cap\F^*(\ca L^+).   
$$ 
The correspondence given by Theorem II.A2 then yields $F^*(\ca L^+)=F^*(\ca L^+)$. 
\qed 
\enddemo

Set $F^*(\ca F)=\ca F_{S\cap F^*(\ca L)}(F^*(\ca L))$, and  
define $\d(\ca F)$ to be the overgroup-closure of $F^*(\ca F)^s$ in $S$.

\proclaim {Lemma 6.10} $\d(\ca F)$ depends only on $\ca F$ (and not on the choice of proper locality 
on $\ca F$). Moreover, $\d(\ca F)$ is $\ca F$-closed, and 
$$ 
\ca F^{cr}\sub\del(\ca F)\sub\d(\ca F)\sub\ca F^s. 
$$
\endproclaim 

\demo {Proof} As $\del(\ca F)$ depends only on $\ca F$ by 6.2(a), and since $\ca L^s$ depends only on 
$\ca F$ by Theorem II.A1, it follows that $\d(\ca F)$ depends only on $\ca F$. Now let $\ca M$ be an 
arbitrary large partial normal subgroup of $\ca L$ containing $O_p(\ca L)$, set $R=S\cap\ca M$, and 
set $\ca D=\ca F_R(\ca M)$. As $\S_R(\ca F)\sub\del(\ca F)$ by 6.3, it follows from 3.14 that $\ca D$ is 
independent of the choice of $\D$, provided only that $\del(\ca F)\sub\D\sub\ca F^s$.   
Taking $\D=\ca F^s$, we obtain $\ca D^c\sub\ca F^q$ by 4.5(a), so $R\in\ca F^q$, and then 
$\ca D^s\sub\ca F^s$ by 4.3(c). In particular, by taking $\ca M=F^*(\ca L)$, we obtain 
$F^*(\ca F)^s\sub\ca F^s$. 

As $\ca D^{cr}\sub\D$, it follows from 2.8 that $\ca D^s$ is $\ca F$-invariant, and 
hence the overgroup-closure of $\ca D^s$ in $S$ is $\ca F$-closed. By specializing again to 
$\ca M=F^*(\ca L)$, we conclude that $\d(\ca F)$ is $\ca F$-closed. As $\ca F^{cr}\sub\del(\ca F)$ by 
6.2(b), it now only remains to show that $\del(\ca F)\sub\d(\ca F)$. 

Let $P\in\ca D^{cr}$. Then there exists $Q\in P^{\ca F}$ such that $P\cap F^*(\ca L)$ is fully normalized 
in $\ca F$, by II.1.15. We note that $\ca D^{cr}$ is $\ca F$-invariant (again by 2.8), 
so $Q\in\ca D^{cr}$, and then $Q\cap F^*(\ca L)\in F^*(\ca F)^{cr}$ by 2.6(c). Thus $Q\in\d(\ca F)$, 
and then $P\in\d(\ca F)$ as $\d(\ca F)$ is $\ca F$-invariant. Thus $\del(\ca F)\sub\d(\ca F)$, and 
the proof is complete. 
\qed 
\enddemo

\vskip .2in 
\noindent 
{\bf Section 7: Regular localities} 
\vskip .1in 

In the preceding section we 
produced a mapping $\d$ which, to each saturated fusion system $\ca F$ on a $p$-group $S$, assigns 
the overgroup closure in $S$ of $\ca F^*(\ca F)^s$. 
By 6.8(a) and Theorem II.A1 there is a unique (up to unique 
isomorphism) proper locality $\ca L$ on $\ca F$ whose set of of objects is $\d(\ca F)$.

\definition {Definition 7.1} A locality $(\ca L,\D,S)$ on $\ca F$ is {\it regular} provided that 
$\ca L$ is proper, and $\D=\d(\ca F)$. 
\enddefinition

\proclaim {Proposition 7.2} Let $(\ca L,\D,S)$ be a regular locality, let $Z\leq Z(\ca L)$ be a subgroup 
of $Z(\ca L)$, and let $\r:\ca L\to\ca L/Z$ be the canonical projection. Then $\ca L/Z$ is regular, and 
$F^*(\ca L)\r =\ca F^*(\ca L/Z)$. 
\endproclaim 

\demo {Proof} Set $\bar{\ca L}=\ca L/Z$, and write $\bar X$ for the image under $\r$ of any subset of 
$\ca L$, or any set of subgroups of $\ca L$. Set $\ca F=\ca F_S(\ca L)$ and set 
$\bar{\ca F}=\ca F_{\bar S}(\bar{\ca L})$. The restriction $\s:S\to\bar S$ of $\r$ to $S$ is 
fusion-preserving, and we shall denote also by $\s$ the associated homomorphism $\ca F\to\bar{\ca F}$ of 
fusion systems. Then $\s$ maps the set $Hom(\ca F)$ of $\ca F$-homomorphisms onto $Hom(\bar{\ca F})$, so 
II.1.19 applies and yields the following information concerning a subgroup $\bar X$ of $\bar S$ and 
its preimage $X$ in $S$.   
\roster 

\item "{(1)}" $X$ is fully normalized in $\ca F$ if and only if $\bar X$ is fully normalized in 
$\bar{\ca F}$.  

\item "{(2)}" $\bar{O_p(N_{\ca F}(X))}\leq O_p(N_{\bar F}(\bar X))$. 

\item "{(3)}" If $\bar X\in\bar{\ca F}^c$ then $X\in\ca F^c$, and if $\bar X\in\bar{\ca F}^{cr}$ then 
$X\in\ca F^{cr}$. 

\endroster 
As a consequence of (3) we have $\bar{\ca F}^{cr}\leq\bar\D$. For any $P\in\D$ the group $N_{\ca L}(P)/Z$ 
is of characteristic $p$, by II.2.6(c), so $\bar{\ca L}$ is a proper locality. We next show: 
\roster 

\item "{(4)}" $\bar{\ca F^s}=\bar{\ca F}^s$. 

\endroster  
In one direction: let $Z\leq P\leq S$ with $P$ fully normalized in $\ca F$ and with $\bar P\in\bar{\ca F}^s$. 
Let $Q$ be the pre-image in $S$ of $O_p(N_{\bar{\ca L}}(\bar P))$. Then $Q=O_p(N_{\ca L}(P))$, and 
$Q\in\ca F^c$ by (3). Thus $\bar{\ca F^s}\sup\bar{\ca F}^s$. In the other direction: let 
$P\in\ca F^s$. Then $ZP\in\ca F^s$, as is any $\ca F$-conjugate of $P$. In order to show that 
$\bar P\in\bar{\ca F}^s$ we may therefore assume that $Z\leq P$, that $P$ is fully normalized in $\ca F$, 
and (by II.1.18) that $Q:=O_p(N_{\ca F}(P))$ is fully normalized in $\ca F$. Let $D$ be the pre-image of 
$C_{\bar S}(\bar Q)$ in $S$. Thus $[Q,D]\leq Z$. Let $(\ca L^s,\ca F^s,S)$ be the expansion of $\ca L$  
to a proper locality on $\ca F$ whose set of objects is $\ca F^s$, and set $M=N_{\ca L^s}(P)$. 
Thus $M$ is a subgroup of $\ca L^s$ of characteristic $p$. Here $\ca L^s$ is generated by $\ca L$ as a 
partial group, by point (c) in Theorem II.A1, so $Z\leq Z(\ca L^s)$. Thus $Z\leq Z(M)$, and $D$ centralizes 
the chain $Q\geq Z\geq 1$ of normal subgroups of $M$. Now II.2.6(c) shows that $D=Q$, and thus $\bar Q$ is 
centric in $\bar{\ca F}$. Thus (4) holds. We next show: 
\roster 

\item "{(5)}" Let $\ca M\norm\ca L$ be a partial normal subgroup of $\ca L$ containing $O_p(\ca L)$. Then 
$\ca M$ is large in $\ca L$ if and only in $\ca M/Z$ is large in $\ca L/Z$.  

\endroster 
Evidently $C_S(\ca M)/Z\leq C_{\bar S}(\bar{\ca M})$, so if $\ca M/Z$ is large in $\ca L/Z$ then $\ca M$ 
is large in $\ca L$. Conversely, assume that $\ca M$ is large in $\ca L$, set $R=S\cap\ca M$, set 
$\bar{\ca M}=\ca M/Z$, and let $\ca N$ be the pre-image in $\ca L$ of $\bar{\ca M}^\perp$. Then 
$\ca N\norm\ca L$, and $[\ca N,R]\leq Z$. As $S\cap F^*(\ca L)\in\D$ and $F^*(\ca L)\leq\ca M$ we have 
$R\in\D$, and thus $\ca N$ is a subgroup of the group $N_{\ca L}(R)$. Then $[O^p(\ca N),R]=1$. 
By 6.3 and 6.10 we have $\S_R(\ca F)\sub\D$, so $O^p(\ca N)\leq\ca M^\perp$ by 6.5. Thus 
$S\cap O^p(\ca N)\leq Z(\ca N)$, and $O^p(\ca N)$ is a $p'$-group, normal in $N_{\ca L}(R)$. As 
$N_{\ca L}(R)$ is of characteristic $p$ we conclude that $\ca N$ is a $p$-group, so 
$\ca N\leq O_p(\ca L)$, and thus $\ca N\leq\ca M$. This shows that $\bar{\ca M}$ is large in 
$\bar{\ca L}$, and completes the proof of (5).  

It is immediate from (5) that $F^*(\ca L)/Z=F^*(\ca L/Z)$. Then (4) yields $\bar\D=\d(\bar{\ca F})$ and 
completes the proof.  
\qed 
\enddemo

\proclaim {Lemma 7.3} Let $P\leq S$ be a subgroup of $S$, and suppose that $O_p(\ca L)P\in\d(\ca F)$. Then 
$P\in\d(\ca F)$. 
\endproclaim 

\demo {Proof} Set $\ca M=F^*(\ca L)$, $R=S\cap\ca M$, and $\ca D=\ca F_R(\ca M)$. As 
$O_p(\ca L)P\in\d(\ca F)$ we have $R\cap O_p(\ca L)P\in\ca D^s$. As $O_p(\ca L)\leq R$ we have 
$R\cap O_p(\ca L)P=O_p(\ca L)(R\cap P)$, and then $R\cap P\in\ca D^s$ by II.6.2(c). That is, $P\in\d(\ca F)$. 
\qed 
\enddemo

For the remainder of this section $(\ca L,\D,S)$ will be a regular locality on $\ca F$. The main goal 
is to show that every partial normal subgroup of $\ca L$ is itself a regular locality.

\proclaim {Lemma 7.4} Let $(\ca L,\D,S)$ be a regular locality on $\ca F$, and let $\ca N\norm\ca L$ be 
a partial normal subgroup of $\ca L$ such that $O_p(\ca L)C_S(\ca N)\leq\ca N$. Adopt the notation of 2.1, 
and set $\G=\{U\leq T\mid U\in\D\}$. Then $(\ca N,\G,T)$ is a regular locality on $\ca E$, and 
$F^*(\ca N)=F^*(\ca L)$. 
\endproclaim 

\demo {Proof} The hypothesis that $O_p(\ca L)C_S(\ca N)$ be contained in $\ca N$ is equivalent to   
$F^*(\ca L)\leq\ca N$. Then $\D$ is the overgroup closure in $S$ of $\G$, and so $S_w\cap T\in\G$ 
for each $w\in\bold W(\ca N)\cap\bold D$. Then $(\ca N,\G,T)$ is a locality by 1.2, and indeed a 
locality on $\ca E$ by the definition $\ca E$ as $\ca F_T(\ca N)$. Here $\ca E^{cr}\sub\del(\ca F)$, so
$\ca E^{cr}\sub\D$ by 6.10; and $\ca N$ is then a proper locality 
by II.2.6(a). We note that $N_{\ca L}(T)$ acts on $\ca L$ (and hence on $\ca N$, and also on $\ca E$) by 
conjugation, since $P\cap T\in\D$ for all $P\in\D$. Then $O_p(\ca N)=O_p(\ca L)$ by the Frattini Lemma 
(I.3.12).  

Set $\G^+=\ca E^s$, and let $\D^+$ be the overgroup closure of $\G^+$ in $S$. As $T\in\ca F^q$ by 4.5(a), 
we have $\G^+\sub\ca F^s$ by 4.3(c), and then $\D^+\sub\ca F^s$ by II.6.2(a). As $N_{\ca L}(T)$ 
acts on $\ca E$, it follows that $\D^+$ is $\ca F$-closed. Let $(\ca L^+,\D^+,S)$ be the expansion of 
$\ca L$ to a proper locality on $\ca F$, via Theorem II.A1. For any $\ca K\norm\ca L$ let 
$\ca K^+\norm\ca L^+$ be the partial normal subgroup of $\ca L^+$, given by Theorem II.A2, whose 
intersection with $\ca L$ is $\ca K$. Notice that $(\ca N^+,\G^+,T)$ is a proper locality on $\ca E$, 
and so $\ca N^+$ is an expansion of $\ca N$. As $N_{\ca L}(T)=N_{\ca L^+}(T)$ permutes the set 
$\Bbb M(\ca N^+)$ of large partial normal subgroups $\ca K$ of $\ca N^+$ containing $O_p(\ca N)$, it 
follows from the Frattini Lemma that $F^*(\ca N^+)\norm\ca L^+$. 

Set $\ca K^+=F^*(\ca N^+)$. Applying 6.4 with $\ca N^+$ and $\ca K^+$ in the roles of $\ca L$ and $\ca N$, 
we obtain $T\cap(\ca K^+)^\perp\leq Z(\ca K^+)$. Then 6.5 yields 
$$ 
O^p_{\ca L^+}((\ca K^+)^\perp)\leq(\ca N^+)^\perp. 
$$ 
But $(\ca N^+)^\perp=(\ca N^\perp)^+$ by 5.7(d), and $\ca N^\perp=Z(\ca N)$ by 6.4. Thus 
$(\ca K^+)^\perp$ is a normal $p$-subgroup of $\ca L^+$, and so $(\ca K^+)^\perp\leq O_p(\ca L^+)$. 
As $O_p(\ca L^+)=O_p(\ca N^+)=O_p(\ca K^+)$, we conclude that $\ca K^+$ is a large partial normal 
subgroup of $\ca L^+$, and hence $F^*(\ca L^+)\leq\ca K^+$. The reverse inclusion holds since 
$F^*(\ca L)^+$ is a large partial normal subgroup of $\ca N^+$ containing $O_p(\ca N^+)$, and since 
$\ca K^+$ is by definition the intersection of all such. Moreover, Theorem II.A2 then yields: 
$$
\ca K^+\cap\ca L=F^*(\ca L)^+\cap\ca L=F^*(\ca L),  
$$ 
which in turn yields $F^*(\ca E)=\ca F^*(\ca F)$. Then $\G=\d(\ca E)$, and the proof is complete. 
\qed 
\enddemo

\proclaim {Lemma 7.5} Let $(\ca L,\D,S)$ be a regular locality on $\ca F$, and let $\ca N\norm\ca L$ 
be a partial normal subgroup of $\ca L$ such that $\ca L=O_p(\ca L)\ca N$. Adopt the notation of 2.1,  
and set $\G=\{U\leq T\mid U\in\D\}$. Then $(\ca N,\G,T)$ is a regular locality on $\ca E$, and 
$F^*(\ca N)=F^*(\ca L)\cap\ca N$. 
\endproclaim 

\demo {Proof} Set $\ca K=F^*(\ca L)\cap\ca N$. Then 
$$ 
F^*(\ca L)=F^*(\ca L)\cap\ca L=F^*(\ca L)\cap O_p(\ca L)\ca N=O_p(\ca L)\ca K, 
$$ 
by the Dedekind Lemma. Then 7.3 implies that $P\cap T\in\G$ for each $P\in\D$ such that $O_p(\ca L)\leq P$.  
Then $S_w\cap T\in\G$ for each $w\in\bold W(\ca N)\cap\bold D$, and so $(\ca N,\G,T)$ is a locality on 
$\ca E$. 

As in the proof of 7.4, it will be necessary to consider ``versions" of $\ca L$ and of $\ca N$ with sets 
of objects other than $\D$ and $\G$. Let $\D_0$ be the overgroup closure of $\G$ in $S$. Thus $P\in\D_0$ 
if and only if $P\cap T\in\G$, so $O_p(\ca L)\D\sub\D_0\sub\D$, and the restriction of $\ca L$ to $\D_0$ is 
then equal (as a partial group) to $\ca L$. Set $\G^+=\ca E^s$. As $C_{\ca L}(S)$ is a $p$-group, we have 
$O^p(C_{\ca L}(T))=1$, and so $T\in\ca F^q$ by II.2.8(c). Then $\G^+\sub\ca F^s$ by 4.3(c). Let 
$\D_0^+$ be the overgroup closure of $\G^+$ in $S$. Then $\D_0^+$ is $\ca F$-closed, and 
$\D_0\sub\D_0^+\sub\ca F^s$. Let $(\ca L^+,\D_0^+,S)$ be the $\D_0^+$-expansion of $\ca L$ to a proper 
locality on $\ca F$, and for each partial normal subgroup $\ca H$ of $\ca L$ let $\ca H^+$ be the 
partial normal subgroup of $\ca L^+$ such that $\ca H=\ca H^+\cap\ca L$. As 
$O_p(\ca L)\ca N^+\cap\ca L=\ca L$ we then have $\ca L^+=O_p(\ca L)\ca N^+$. 

Let $P\in\ca E^{cr}$ and set $Q=O_p(\ca L)P$. Then $Q\in\D_0^+$, and 
$$ 
N_{\ca L^+}(Q)=O_p(\ca L)N_{\ca N^+}(Q)=O_p(\ca L)N_{\ca N^+}(P).  
$$  
Since $P\in\ca E^{cr}$ we have $P=O_p(N_{\ca N^+}(P))$, and then $Q=O_p(N_{\ca L^+}(Q))$. Then 
$Q\in\ca F^{cr}$ by II.2.8(a), so $Q\in\D_0$. Then $P=Q\cap T\in\G$, and thus $\ca E^{cr}\sub\G$. 
We have thus shown that $(\ca N,\G,T)$ is proper. 

Set $\ca H=\ca K^\perp\cap\ca N$. As $O_p(\ca N)\leq\ca K$, 5.7(c) yields $[O_p(\ca N),\ca H]=1$. 
As $[O_p(\ca L),\ca H]\leq O_p(\ca N)$ it follows from II.2.6(c) that $[O_p(\ca L),O^p_{\ca N}(\ca H)]=1$. 
Thus $T\cap O^p_{\ca N}(\ca H)\leq C_S(F^*(\ca L))$, and 6.4 then yields  
$$ 
T\cap O^p_{\ca N}(\ca H)\leq Z(O^p_{\ca N}(\ca H)). 
$$ 
Then $O^p_{\ca N}(\ca H)$ is a $p$-group by 3.17, and so $\ca H\leq T$. As $\ca H\norm\ca L$ we get 
$\ca H\leq O_p(\ca L)$. Then $\ca H\leq\ca K$, and since  
$C_T(\ca K)\leq\ca H$ we conclude that $O_p(\ca N)C_T(\ca K)\leq\ca K$. Thus $F^*(\ca N)\leq\ca K$, and 
$O_p(\ca L)F^*(\ca N)\leq F^*(\ca L)$.  

In order to obtain the reverse inclusion, observe first of all that 
$$ 
C_S(O_p(\ca L)F^*(\ca N^+))\leq O_p(\ca L)C_S(F^*(\ca N)^+)=O_p(\ca L)C_T(F^*(\ca N^+))\leq O_p(\ca L). 
$$ 
This shows:
$$ 
F^*(\ca L^+)\leq O_p(\ca L)F^*(\ca N^+). \tag* 
$$
Since $F^*(\ca N^+)=F^*(\ca N)^+$ by 6.9, and since 
$$ 
O_p(\ca L)F^*(\ca N)^+\cap\ca L=O_p(\ca L)(F^*(\ca N)^+\cap\ca L)=O_p(\ca L)F^*(\ca N), 
$$ 
we get  
$$
(O_p(\ca L)F^*(\ca N))^+=O_p(\ca L)F^*(\ca N)^+=O_p(\ca L)F^*(\ca N^+).  
$$ 
Then (*) yields $F^*(\ca L)\leq O_p(\ca L)F^*(\ca N)$, and thus 
$F^*(\ca L)=O_p(\ca L)F^*(\ca N)$. 

We have still to show that $(\ca N,\G,T)$ is a regular locality. Thus, it remains to show that 
$\G=\d(\ca E)$. Let $X\in\d(\ca E)$ and set $U=X\cap F^*(\ca N)$. As $\ca E=\ca E^+$ by 3.14, the 
definition of $\d(\ca E)$ yields $U\in F^*(\ca E)^s$. Here $T\cap F^*(\ca E)\in\ca E^q$ by 4.5(a), 
so $U\in\ca E^s$ by 4.3(c). Then $U\in\ca F^s$ since we have already seen that $\ca E^s\sub\ca F^s$. 
As $U\leq F^*(\ca N)\leq F^*(\ca L)$ it follows that $U\in F^*(\ca F)^s$. Then $U\in\d(\ca F)=\D$ and 
$X\in\D$. Thus $X\in\G$, and we have $\d(\ca E)\sub\G$. Now let $Y\in\G$ and set $V=Y\cap F^*(\ca L)$. 
Then $V\in\G$ since $\D=\d(\ca F)$. Then $N_{\ca L}(V)$ is of characteristic $p$, so also 
$N_{F^*(\ca N)}(V)$ is of characteristic $p$ by II.2.6(a). As $V\leq F^*(\ca L)\cap\ca N=F^*(\ca N)$ 
we conclude that $V\in F^*(\ca E)^s$, so $V\in\d(\ca E)$. Then $Y\in\d(\ca E)$, and thus  
$\G=\d(\ca E)$, as required.  
\qed 
\enddemo

\proclaim {Lemma 7.6} Let $(\ca L,\D,S)$ be a regular locality on $\ca F$, and assume the ``product setup"      
of 2.9, with $\ca L=\ca N_1\ca N_2$. Then: 
\roster 

\item "{(a)}" $\ca N_i\leq(\ca N_{3-i})^\perp$ for $i=1$ and $2$. 

\item "{(b)}" Each $(\ca N_i,\d(\ca E_i),T_i)$ is a regular locality on $\ca E_i$, and 
$\d(\ca E_i)=\{U\leq T_i\mid UC_S(\ca N_i)\in\D\}$.  

\item "{(c)}" $F^*(\ca L)=F^*(\ca N_1)F^*(\ca N_2)$.  

\endroster 
\endproclaim 

\demo {Proof} We have $\ca N_2\leq C_{\ca L}(T_1)$ by 2.9. Applying II.7.4 to the locality 
$N_{\ca L}(T_1)$ then yields 
$$ 
O^p{N_{\ca L}(T_1)}(\ca N_2)\leq O^p_{N_{\ca L}(T_1)}(C_{\ca L}(T_1))\leq (\ca N_1)^\perp. \tag1
$$ 
By definition, $O^p_{N_{\ca L}(T_1)}(C_{\ca L}(T_1))$ is the intersection of the set $\Bbb K$ of partial 
normal subgroups $\ca K$ of $N_{\ca L}(T_1)$ such that: 
\roster 

\item "{(2)}" $\ca K\leq\ca N_2$ and $\ca KT_2=\ca N_2$. 

\endroster 
Evidently $\Bbb K$ contains the set of partial normal subgroups $\ca K$ of $\ca L$ such that (2) holds, 
so $O^p_{\ca L}(\ca N_2)\leq O^p{N_{\ca L}(T_1)}(\ca N_2)$. Then (1) yields 
$O^p_{\ca L}(\ca N_2)\leq(\ca N_1)^\perp$. Since $\S_{T_1}(\ca F)\sub\D$ we have 
$C_S(\ca N_1)=(\ca N_1)^\perp$ by 5.7(b), and so $T_2\leq(\ca N_1)^\perp$. Thus $\ca N_2\leq(\ca N_1)^\perp$, 
and point (a) follows. 

Set $\ca M=F^*(\ca L)$, set $\ca M_1=\ca M\cap\ca N_1$, and set $\ca K=(\ca M_1)^\perp\cap\ca N_1$. 
Then $\ca K\norm\ca L$ by 5.7(a), and 
$$ 
\ca K\cap\ca M=\ca K\cap\ca M_1\leq Z(\ca K), 
$$
so 6.5 yields $O^p_{\ca L}(\ca K)\leq\ca M^\perp$. Now 4.7 yields $\ca M^\perp=Z(\ca M)$ and 
$O^p_{\ca L}(\ca K)=\1$. Thus $\ca K$ is a normal $p$-subgroup of $\ca L$, and so $\ca K\leq O_p(\ca L)$. 
Notice that since $[\ca N_i,T_{3-i}]=\1$ we have $O_p(\ca L)=O_p(\ca N_1)O_p(\ca N_2)$. As  
$O_p(\ca L)\leq\ca M$ we then have $O_p(\ca L)=O_p(\ca M_1)O_p(\ca M_2)$, and $\ca K\leq O_p(\ca M_1)$.  
We now compute: 
$$ 
\align 
C_{T_1}(\ca M_1)&\leq C_S(\ca M_1\ca M_2)=C_S(\ca M_1)\cap C_S(\ca M_2) \\
&=C_{T_1}(\ca M_1)T_2\cap C_{T_2}(\ca M_2)T_1=Z(\ca M_1)T_2\cap Z(\ca M_2)T_1 \\
&=Z(\ca M_1)(T_2\cap Z(\ca M_2)T_1)=Z(\ca M_1)Z(\ca M_2)(T_1\cap T_2)\leq O_p(\ca L)\\
&\leq\ca M_1\ca M_2.  
\endalign 
$$ 
This shows that $C_{T_1}(\ca M_1)\leq\ca M_1$, and that $\ca M_1\ca M_2$ is a large partial normal subgroup 
of $\ca L$ containing $O_p(\ca L)$. We have shown: 
\roster 

\item "{(3)}" $F^*(\ca L)=\ca M_1\ca M_2$, and 

\item "{(4)}" $O_p(\ca N_i)C_S(\ca M_i)\leq \ca M_i$ ($i=1,2)$.  

\endroster 

Next, set $\G_i=\{PT_{3-i}\cap T_i\mid P\in\D\}$ and set $\D_0=\G_1\G_2$. Then $(\ca N_i,\G_i,T_i)$ is a 
proper locality on $\ca E_i$, $\D_0$ is an $\ca F$-closed subset of $\D$ containing $\ca F^{cr}$, and 
the restriction of $ca L$ to $\D_0$ is equal to $\ca L$ as a partial group, by 2.11 and by 2.10(a). Set 
$\G_i^+=\ca E_i^s$ and let $\D_0^+$ be the overgroup closure of $\G_1^+\G_2^+$. Then $\D_0^+$ is an 
$\ca F$-closed subset of $\ca F^s$ by 2.10(e). Form the corresponding expansions $(\ca N_i)^+$ and 
$\ca L^+$ via Theorem II.A1. Let $\ca K_i^+=F^*(\ca N_i)^+$ be the partial normal subgroup of $(\ca N_i)^+$ 
corresponding to $\ca K_i=F^*(\ca N_i)$ via Theorem II.A2. Note that $\ca N_i^+\leq(\ca N_{3-i}^+)^\perp$  
(a) and 5.7(d). Then 5.7(c) implies that $\ca K_i\norm\ca L$ and $\ca K_i^+\norm\ca L^+$. 
Let $(\ca K_1\ca K_2)^+$ be the partial normal subgroup of $\ca L^+$ whose intersection with $\ca L$ 
is $\ca K_1\ca K_2$. Then $(\ca K_1\ca K_2)^+=\ca K_1^+\ca K_2^+$ by II.5.3. We compute: 
$$ 
\align 
C_S(\ca K_1^+\ca K_2^+)&=C_S(\ca K_1^+)\cap C_S(\ca K_2^+)=C_{T_1}(\ca K_1^+)T_2
\cap C_{T_2}(\ca K_2^+)T_1 \\ 
&=Z(\ca K_1^+)T_2\cap Z(\ca K_2^+)T_1=Z(\ca K_1^+)Z(\ca K_2^+)(T_1\cap T_2)\leq\ca K_1^+\ca K_2^+,  
\endalign 
$$ 
and so $(\ca K_1\ca K_2)^+$ is a large partial normal subgroup of $\ca L^+$. As also 
$O_p(\ca F)=O_p(\ca L)\leq\ca K_1\ca K_2$, we conclude that 
$F^*(\ca L^+)\leq(\ca K_1\ca K_2)^+$. This yields $\ca M\leq\ca K_1\ca K_2$. Since also 
$\ca K_i\leq\ca M_i$ by (4), we obtain $\ca K_i=\ca M_i$. That is, (c) holds. 

Let $U_i\in F^*(\ca E_i)^s$. Then $U_i\in\d(\ca E_i)\sub(\ca E_i)^s$, and then $U_1U_2\in\ca F^s$ 
by 2.10(e). As $U_1U_2\leq F^*(\ca L)$ by (c), we obtain $\d(\ca E_1)\d(\ca E_2)\sub\d(\ca F)$. 
As $\d(\ca F)=\D$ we have $U_1(T_2\cap\ca M_2)\in\D$ in particular, and so $U_1(T_1\cap T_2)\in\G_1$. 
This shows that $Z(\ca N_1)\d(\ca E_1)\sub\G_1$.  

Let $Q_1\in\G_1$. Then $Q_1=P\cap T_1$ for some $P\in\D$ with $T_2\leq P$. Then 
$P\cap\ca M\in F^*(\ca F)^s$ as $\D=\d(\ca F)$. Here $P\cap\ca M=(P\cap\ca M_1)(T_2\cap\ca M_2)$, 
so $P\cap\ca M_1\in F^*(\ca E_1)^s$ by 2.10(e). Thus $Q_1\in\d(\ca E_1)$, and this shows that 
$\G_1\sub\d(\ca E_1)$. Since $Z(\ca N_1)\leq S_w$ for all $w\in\bold W(\ca N_1)$, the expansion of 
$(\ca N_1,\G_1,T_1)$ to a regular locality on $\d(\ca E_1)$ has the same underlying partial group 
$\ca N_1$. This completes the proof of (b), and of the lemma. 
\qed 
\enddemo

\proclaim {Theorem 7.7} Let $(\ca L,\D,S)$ be a regular locality on $\ca F$, and let $\ca N\norm\ca L$ be 
a partial normal subgroup. Set $T=S\cap\ca N$, $\ca E=\ca F_T(\ca N)$, and  
$\G=\{U\leq T\mid C_S(\ca N)U\in\D\}$. Then: 
\roster 

\item "{(a)}" $(\ca N,\G,T)$ is a regular locality on $\ca E$.  

\item "{(b)}" $F^*(\ca N)=\ca N\cap F^*(\ca L)$. 

\item "{(c)}" $\ca N^\perp=C_{\ca L}(\ca N)$. That is, $\ca N^\perp$ is the set of elements $g\in\ca L$ 
such that $x^g=x$ for all $x\in\ca N$. 

\item "{(d)}" $O_p(\ca E)=O_p(\ca N)\norm\ca L$. 

\item "{(e)}" $(g\i,x,g)\in\bold D$ for all $x\in\ca N$ and all $g\in N_{\ca L}(T)$. Moreover, the mapping 
$c_g:\ca N\to\ca N$ given by $x\maps\Pi(g\i,x,g)$ is an automorphism of $\ca N$, and the mapping 
$\eta:N_{\ca L}(T)\to Aut(\ca N)$ given by $g\maps c_g$ is a homomomorphism of partial groups with kernel 
$\ca N^\perp$. 

\item "{(f)}" The image of the natural homomorphism $N_{\ca L}(T)\to Aut(T)$ is contained in $Aut(\ca E)$. 

\endroster 
\endproclaim 

\demo {Proof} If $\ca N$ is large and $O_p(\ca L)\leq\ca N$ then (a) and (b) are given by 7.3 and 7.4. Set 
$\ca M=\ca N\ca N^\perp$. Thus (a) and (b) hold with $O_p(\ca L)\ca M$ in the role of $\ca N$. In order 
to prove (a) and (b) in general, we may then assume that $\ca L=O_p(\ca L)\ca M$. By 7.5 we then obtain (a) 
and (b) for $\ca M$ in the role of $\ca N$. In this way we reduce to the case where $\ca L=\ca M$, where 
(a) and (b) are then given by 7.6. 

In proving (c) we may now assume that $\ca L$ is the regular locality $\ca N\ca N^\perp$. Let $f\in\ca N$ and 
$g\in\ca N^\perp$, set $U=S_f\cap T$, and set $V=S_g\cap T^\perp$. Then $UV\in\D$ by 7.6(b), and so 
$(f,g)\in\bold D$. Then $[f,g]=\1$ by 5.7(c), and this yields point (c). That $O_p(\ca E)=O_p(\ca N)$ is 
now a consequence of (a) and II.2.3. Point (d) will then follow from (f) and the factorization 
$\ca L=N_{\ca L}(T)\ca N$ given by the Frattini Lemma. Thus, it remains to prove (e) and (f). 

Set $H=N_{\ca L}(TT^\perp)$ and set $\ca L_T=N_{\ca L}(T)$. Then $(\ca L_T,\D,S)$ is a locality, 
and $\ca N^\perp\norm\ca L_T$, so the Frattini Lemma yields $\ca L_T=\ca N^\perp H$. Here $TT^\perp\in\D$ 
by an application of (a) to $\ca N\ca N^\perp$ in the role of $\ca N$, so $H$ is a subgroup of $\ca L$. 
Let $x\in\ca N$ and set $P=S_x\cap T$. Then $P\in\G$. Now let $g\in\ca L_T$ and employ the Splitting 
Lemma so as to write $g=yh$ with $y\in\ca N^\perp$ and $h\in\ca H$, and with $S_g=S_{(y,h)}$. 
Set $Q=S_y\cap T^\perp$. Then $PQ\in\D$ by 2.11, and one may verify that the word $(h\i,y\i,x,y,h)$ 
is in $\bold D$ via $(PQ)^h$. Thus $(g\i,x,g)\in\bold D$, and (f) follows from 2.3. A similar argument 
shows that if $(x_1,\cdots,x_n)\in\bold D(\ca N)$ then 
$$ 
(g\i,x_1,g,\cdots,g\i,x_n,g)\in\bold D, 
$$ 
and hence that the conjugation map $c_g:\ca N\to\ca N$ is an endomorphism of $\ca N$ as a partial group. 
One checks that $c_g\circ c_{g\i}$ is the identity map on $\ca N$, so we have a mapping 
$\eta:\ca L_T\to Aut(\ca N)$ given by $g\maps c_g$. The restriction of $\eta$ to $H$ is easily seen to be 
a homomorphism (of groups), and a further exercise with the Splitting Lemma will verify that $\eta$ itself 
is a homomorphism of partial groups, completing the proof of (e). 
\qed 
\enddemo

\definition {Definition 7.8} 
A partial subgroup $\ca K$ of a locality $\ca L$ is {\it subnormal} in $\ca L$ (denoted 
$\ca K\norm\norm\ca L$) if there exists a sequence $(\ca N_0,\cdots,\ca N_k)$ of partial subgroups 
of $\ca L$ such that $\ca K=\ca N_0\norm\ca N_1\norm\cdots\norm\ca N_k=\ca L$. 
\enddefinition

\proclaim {Corollary 7.9} Let $\ca L$ be a regular locality and let $\ca K\norm\norm\ca L$ be a 
partial subnormal subgroup of $\ca L$. Then $\ca K$ is a regular locality. 
\endproclaim 

\demo {Proof} By induction on the length of a subnormal chain from $\ca K$ to $\ca L$ it suffices to show 
that $\ca K$ is regular in the case that $\ca K\norm\ca L$, in which case we are done by 7.7(a). 
\qed 
\enddemo

\proclaim {Lemma 7.10} Let $\ca L$ be a regular locality, let $\ca N\norm\ca L$ be a partial normal 
subgroup, and let $\Bbb K$ be a non-empty set of partial normal subgroups of $\ca N$. Set $T=S\cap\ca N$, 
and let $\L=N_{Aut(\ca N)}(T)$ be the set of automorphisms of $\ca N$ (as a partial group) which leave $T$ 
invariant. Suppose that $\Bbb K$ is invariant under $\L$, and set $\ca K=\bigcap\Bbb K$. Then 
$\ca K\norm\ca L$. 
\endproclaim 

\demo {Proof} We have $\ca K\norm\ca N\ca N^\perp$ by 7.7(c), and $\ca K$ is $S$-invariant since, by I.2.9, 
$S$ acts on $\ca L$ by conjugation. Thus $\ca K\norm\ca M:=O_p(\ca L)\ca N\ca N^\perp$. 
Set $H=N_{\ca L}(S\cap\ca M)$. Then $H$ is a subgroup of $N_{\ca L}(S\cap F^*(\ca L))$. 
As $\ca L$ is regular we have $\bold D=\{w\in\bold W(\ca L)\mid S_w\cap F^*(\ca L)\in\D\}$, so $H$ acts on 
$\ca L$ by conjugation. Then $O^p(\ca N)$ is $H$-invariant, and we may employ I.3.13 to conclude that 
$\ca K\norm\ca L$. 
\qed 
\enddemo 

Recall from II.7.1 that $O^p(\ca L)$ is defined to be $O^p_{\ca L}(\ca L)$ (and similarly for 
$O^{p'}(\ca L)$). Recall also the definition of $[\ca L,\ca L]$ from 3.16. 

\proclaim {Corollary 7.11} Let $\ca L$ be a regular locality and let $\ca N\norm\ca L$ be a partial 
normal subgroup. Then $O^p(\ca N)$, $O^{p'}(\ca N)$, and $[\ca N,\ca N]$ are partial normal subgroups of 
$\ca L$. Moreover, we have $O^p_{\ca L}(\ca N)=O^p(\ca N)$ and $O^{p'}_{\ca L}(\ca N)=O^{p'}(\ca N)$. 
\endproclaim 

\demo {Proof} Let $\Bbb K$ be the set of partial normal subgroups of $\ca N$ defined by any one of 
the following three conditions: (1) $\ca H\in\Bbb K$ if $\ca HT=\ca N$, (2) $\ca H\in\Bbb K$ if 
$T\leq\ca H$, and (3) $\ca H\in\Bbb K$ if $\ca N/\ca H$ is an abelian group. Set $\ca K=\bigcap\Bbb K$. 
Then $\ca K\norm\ca L$ by 7.10, and $\ca K$ is variously $O^p(\ca N)$, $O^{p'}(\ca N)$, or $[\ca N,\ca N]$. 
Now let $\Bbb K'$ be the set of all partial normal subgroups $\ca H'$ of $\ca L$ such that 
$\ca H'T=\ca N$. Then $O^p_{\ca L}(\ca N)=\bigcap\Bbb K'$ by definition, and 
$O^p_{\ca L}(\ca N)\leq O^p(\ca N)$ since $O^p(\ca N)\in\Bbb K'$. The reverse inclusion holds since 
$O^p_{\ca L}(\ca N)\in\Bbb K$ (with $\Bbb K$ as in (1)). Thus $O^p_{\ca L}(\ca N)=O^p(\ca N)$, and 
similarly $O^{p'}_{\ca L}(\ca N)=O^{p'}(\ca N)$. 
\qed 
\enddemo

\vskip .2in 
\noindent 
{\bf Section 8: Components} 
\vskip .1in 

Throughout this section, $(\ca L,\D,S)$ is a regular locality on $\ca F$. That is, $\ca L$ is a proper 
locality on $\ca F$, and $\D=\d(\ca F)$. By 7.9, every partial subnormal subgroup $\ca K$ of $\ca L$ 
is then a regular locality, and $O^p(\ca K)\norm\norm\ca L$ by 7.11. We may write $\d(\ca L)$ for 
$\d(\ca F)$, and thus write also $(\ca K,\d(\ca K),S\cap\ca K)$ for the regular locality $\ca K$. 
A locality $\ca L$ is {\it simple} if $\ca L$ has exactly two partial normal subgroups.  

\definition {Definition 8.1} A partial subnormal subgroup $\ca K\norm\norm\ca L$ is a {\it component} of 
$\ca L$ if $\ca K=O^p(\ca K)$ and $\ca K/Z(\ca K)$ is simple. 
\enddefinition 

Write $Comp(\ca L)$ for the set of components of $\ca L$, and let 
$$ 
E(\ca L)=\<Comp(\ca L)\> 
$$ 
be the partial subgroup of $\ca L$ generated by the union of the components of $\ca L$. Recall that 
for subsets $X,Y$ of $\ca L$, the notation $[X,Y]=\1$ indicates that the commutators $[x,y]$ are defined and 
are equal to $\1$ for all $x\in X$ and all $y\in Y$.

\proclaim {Lemma 8.2} Let $\ca K\in Comp(\ca L)$ and let $\ca N\norm\ca K$. Then $Z(\ca K)=O_p(\ca K)$, 
and either $\ca N=\ca K$ or $\ca N\leq Z(\ca K)$. 
\endproclaim 

\demo {Proof} As $\ca K$ is itself a regular locality, we may as well assume that $\ca K=\ca L$. We have 
$Z(\ca L)\leq O_p(\ca L)$ by 2.7. Now let $\ca N\norm\ca L$ be given with $\ca N\nleq Z(\ca L)$. Here 
$Z(\ca L)\ca N\norm\ca L$ by I.5.1. As $\ca L/Z(\ca L)$ is simple, the Correspondence Theorem (I.4.8) yields 
$Z(\ca L)\ca N=\ca L$. 

Set $T=S\cap\ca N$. Then 
$$ 
O^p(\ca N)S\geq Z(\ca L)O^p(\ca N)T=Z(\ca L)\ca N=\ca L.  
$$ 
As $O^p(\ca N)\norm\ca L$, by 7.11, we conclude that $O^p(\ca L)\leq O^p(\ca N)$. Then 
$\ca N=\ca L$, since $\ca L=O^p(\ca L)$. Thus, either $\ca N\leq Z(\ca L)$ or 
$\ca N=\ca L$. The special case where $\ca N=O_p(\ca L)$ then yields $O_p(\ca L)\leq Z(\ca L)$, and 
completes the proof. 
\qed 
\enddemo 

\proclaim {Lemma 8.3} Let $\ca K\norm\norm\ca L$. 
\roster 

\item "{(a)}" If $\ca K\leq\ca N\norm\ca L$ then $\ca K\norm\norm\ca N$. 

\item "{(b)}" $O_p(\ca K)\leq O_p(\ca L)$, and if $\ca K\norm\ca L$ then $O_p(\ca K)\norm\ca L$.

\item "{(c)}" $O_p(\ca L)\ca K\norm\norm\ca L$. 

\endroster  
\endproclaim 

\demo {Proof} Point (a) is immediate from the observation that the intersection of partial normal 
subgroups is again a partial normal subgroup. 

In proving (b): suppose first that $\ca K\norm\ca L$. Then $O_p(\ca K)\norm\ca L$ by 7.11, with 
$\{O_p(\ca K)\}$ in the role of $\Bbb K$. Now (b) follows by an obvious induction argument on the length of 
a subnormal chain from $\ca K$ to $\ca L$. 

Set $\bar{\ca L}=\ca L/O_p(\ca L)$ and let $\bar{\ca K}$ be the image of $\ca K$ in $\bar{\ca L}$ under 
the canonical projection $\r:\ca L\to\bar{\ca L}$. The Correspondence Theorem (I.4.7) yields 
$\bar{\ca K}\norm\norm\bar{\ca L}$, and the preimage $\ca H$ of $\bar{\ca K}$ is a partial subnormal 
subgroup of $\ca L$. This yields (c). 
\qed
\enddemo

\proclaim {Lemma 8.4} Let $\ca K\in Comp(\ca L)$. Then: 
\roster 

\item "{(a)}" $Z(\ca K)\leq O_p(\ca L)$. 

\item "{(b)}" $[O_p(\ca L),\ca K]=\1$.   

\item "{(c)}" $\ca K=O^p(O_p(\ca L)\ca K)$. 

\endroster 
\endproclaim 

\demo {Proof} We have $Z(\ca K)\leq O_p(\ca K)$ by 4.8, and then point (a) follows from 6.3(a). 

Let $X$ be a subgroup of $\ca K$ of order prime to $p$. Then $O_p(\ca L)X$ is a subgroup of 
$\ca L$. Set $Z_0=O_p(\ca L)$ and recursively define $Z_n:=[Z_{n-1},X]$ for $n\geq 1$. Then 
$Z_n\leq O_p(\ca L)$ for all $n$, and since $\ca K\norm\norm\ca L$ we obtain $Z_n\leq\ca K\cap O_p(\ca L)$ 
for $n$ sufficiently large. Then $Z_{n+1}=1$ since $O_p(\ca K)\leq Z(\ca K)$; and then $Z_1=1$ by 
coprime action. In particular, we now have $[O_p(\ca L),O^p(N_{\ca K}(P)]=1$ for all $P\in\d(\ca K)$. One 
easily verifies (via I.2.9 for example) that $C_{\ca L}(O_p(\ca L))$ is a partial subgroup of $\ca L$, and 
thus 
$$ 
\<O^p(N_{\ca K}(P))\mid P\in\d(\ca K)\>\leq C_{\ca L}(O_p(\ca L)). 
$$ 
Then $O^p(\ca L)\leq C_{\ca L}(O_p(\ca L))$ by an application of 4.3, with $\ca K$ in the role of $\ca L$.  
As $(x\i,g\i,x,g)\in\bold D$ for all $x\in O_p(\ca L)$ and all $g\in\ca L$, we obtain (b). 

As $O_p(\ca L)\ca K$ is subnormal in $\ca L$ by 6.3(b), $O_p(\ca L)\ca K$ is a regular locality, and we 
may therefore take $\ca L=O_p(\ca L)\ca K$ in proving (c). Then (b) yields $\ca K\norm\ca L$. Now 
$O^p(\ca L)\leq\ca K$ as $\ca KS=\ca L$, and then $O^p(\ca K)\leq O^p(\ca L)$ since 
$$ 
\ca K=\ca L\cap\ca K=O^p(\ca L)S\cap\ca K=O^p(\ca L)T. 
$$ 
The reverse inclusion $O^p(\ca L)\leq O^p(\ca K)$ is given by  
$$ 
O^p(\ca K)S=(O^p(\ca K)T)S=\ca KS=\ca L. 
$$ 
\qed 
\enddemo

\proclaim {Theorem 8.5} Let $(\ca L,\D,S)$ be a regular locality which is not a group of characteristic $p$. 
Then $Comp(\ca L)\neq\nset$. Let $(\ca K_1,\cdots,\ca K_m)$ be a non-redundant list of the components of 
$\ca L$, and let $E(\ca L)$ be the partial subgroup of $\ca L$ generated by $\bigcup Comp(\ca L)$. Then 
$$ 
F^*(\ca L)=O_p(\ca L)E(\ca L)=O_p(\ca L)\ca K_1\cdots\ca K_m, 
$$ 
and 
$$ 
[\ca K_i,\ca K_j]=1\ \ (1\leq i\neq j\leq m). 
$$ 
Further, $E(\ca L)=O^p(F^*(\ca L))$, and $[O_p(\ca L),E(\ca L)]=1$.  
\endproclaim  

\demo {Proof} If $O_p(\ca L)$ is a large partial normal subgroup of $\ca L$ then 
$O_p(\ca L)=F^*(\ca L)\in\D$, and then $\ca L$ is a group of characteristic $p$, contrary to hypothesis. 
Thus $O_p(\ca L)$ is not large. Set $\ca H=O_p(\ca L)^\perp$, and let $\bold X$ be the set of partial 
subnormal subgroups $\ca K$ of $\ca H$ such that $\ca K$ is not a $p$-group. Regard $\bold X$ as a poset 
via inclusion, and let $\ca K$ be minimal in $\bold X$. Then $\ca K=O^p(\ca K)$, and $\ca K/O_p(\ca K)$ 
is simple. As $O_p(\ca K)\leq O_p(\ca L)$ by 8.3(b), and since $\ca K\leq\ca H=C_{\ca L}(O_p(\ca L))$ 
by 7.7(c), we obtain $\ca K\in Comp(\ca L)$. Here $\bold X\neq\nset$ as $\ca H\in\bold X$, so we have 
shown that $Comp(\ca L)\neq\nset$.  

Set $\ca M=F^*(\ca L)$ and let $\ca K\in Comp(\ca L)$. Let $\ca N\norm\ca L$ be a partial normal subgroup of 
$\ca L$ which is minimal with respect to the condition $\ca K\norm\norm\ca N$. If $\ca K=\ca L$ then 
$\ca M=\ca K$ and $\{\ca K\}=Comp(\ca L)$ by 8.2, and there is then nothing to prove. Thus we may assume 
that $\ca K\neq\ca L$. Then $\ca N\neq\ca L$, and induction on $|\ca L|$ then yields $\ca K\leq F^*(\ca N)$. 
Then $\ca K\leq F^*(\ca L)$ by 7.7(b). Thus $E(\ca L)\leq\ca M$, and we may therefore assume at this 
point that $\ca L=\ca M$.  

With $\ca N$ as above we have $\ca L=\ca N\ca N^\perp\neq\ca N$, and $\ca L\neq\ca N^\perp$. As 
$\ca L=F^*(\ca L)=F^*(\ca N)F^*(\ca N^\perp)$, induction on $|S|$ yields $\ca K\norm\ca L=E(\ca L)$. 
Then $\ca L=\ca K\ca K^\perp$ is a product of pairwise commuting components, and it suffices now to 
show that $\ca K^\perp$ contains each component $\ca K'$ of $\ca L$ other than $\ca K$. 

Let $\ca K'\in Comp(\ca L)$, and suppose that $\ca K\neq\ca K'$. Then $\ca K\cap\ca K'\leq Z(\ca K)$ by 
8.2. As $\ca L/Z(\ca L)$ is regular by 7.2, we may assume $Z(\ca L)=\1$. Let $\ca K'$ be a component 
with $\ca K\neq\ca K'$. Then $\ca K'\norm\ca L$ and $[S\cap\ca K,S\cap\ca K']=\1$. As 
$\ca K'=O^p(\ca K')$ we then have $\ca K'\leq O^p(C_{\ca L}(S\cap\ca K)$ by II.7.4, 
and then $\ca K'\leq\ca K^\perp$ as desired, by definition 5.5. 
\qed 
\enddemo 

\proclaim {Corollary 8.6} Let $(\ca L,\D,S)$ be a regular locality and let $\ca N$ be a partial normal 
subgroup of $E(\ca L)$. Then $\ca N=E(\ca N)Z$, where $Z\leq Z(E(\ca L))$ and where $E(\ca N)$ is a 
product of components of $\ca L$. 
\endproclaim 

\demo {Proof} We may assume without loss of generality that $\ca L=E(\ca L)$. Then $O_p(\ca L)=Z(\ca L)$ 
by the final statement in 8.5. We have $\ca N=F^*(\ca N)$ by 7.7(b), and thus $\ca N=O_p(\ca N)E(\ca N)$ 
where $E(\ca N)$ is a product of components of $\ca N$. The components of $\ca N$ are by definition 
components of $\ca L$, and $O_p(\ca N)\leq O_p(E(\ca L))$ by 7.7(d). Set $Z=O_p(\ca N)$. Thus 
$Z\leq Z(\ca L)$ and $\ca N=E(\ca N)Z$. 
\qed 
\enddemo

\vskip .2in 
\noindent 
{\bf Section 9: Im-partial subgroups and $E$-balance} 
\vskip .1in 

Let $(\ca L,\D,S)$ be a regular locality on $\ca F$, and let $X\leq S$ be fully 
normalized in $\ca F$. We shall see that there is then a regular locality $\ca L_X$ on $N_{\ca F}(X)$ 
which can be constructed in a somewhat indirect way from the partial group $N_{\ca L}(X)$. One of the goals 
of this section is to show that if $X\leq F^*(\ca L)$ then $\ca L_X$ can be produced directly as a subset of 
$\ca L$, and that the inclusion map $\ca L_X\to\ca L$ is a homomorphism of partial groups. 

In the category of groups and homomorphisms one has the obvious equivalence between the notions of 
``subgroup of $G$" and ``image of a homomorphism into $G$". The situation is different in the category 
of partial groups, where partial subgroups are indeed images of homomorphisms, but where images of  
homomomorphisms need not be partial subgroups. As an example, take 
$\ca L$ to be the additive group of integers, and let $\ca H$ be the subset $\{-1,0,1\}$ of $\ca L$. 
Then $\ca H$ is the image of a homomorphism of partial groups (see I.1.2 and I.1.12), but $\ca H$ is 
not a partial subgroup of $\ca L$ (since partial subgroups of groups are subgroups).

\definition {Definition 9.1} Let $\ca H$ and $\ca M$ be partial groups. Then $\ca H$ is an {\it im-partial 
subgroup} of $\ca M$ if $\ca H$ is the image in $\ca M$ of a homomorphism of partial groups. Equivalently:  
$\bold D(\ca H)\sub\bold D(\ca L)$ and $Pi_{\ca H}$ is the restriction of $\Pi_{\ca M}$ to $\bold D(\ca H)$. 
Write $\ca K\pce\ca H$ to indicate that $\ca K$ is an im-partial subgroup of $\ca H$. 
\enddefinition 

It is obvious that the relation $\pce$ is transitive. 
Of the long list (analogous to I.1.8) of properties of im-partial subgroups that one might compile, the 
following lemma provides the few instances that will be needed here.

\proclaim {Lemma 9.2} Let $\ca M$ be a partial group, let $\ca H\pce\ca M$ be an im-partial subgroup, and 
let $\ca K\leq\ca M$ be a partial subgroup. Then: 
\roster 

\item "{(a)}" $\ca H\cap\ca K$ is an im-partial subgroup of $\ca K$ and of $\ca M$. 

\item "{(b)}" $\ca H\cap\ca K$ is a partial subgroup of $\ca H$, and if $\ca K\norm\ca M$ then  
$\ca H\cap\ca K$ is a partial normal subgroup of $\ca H$. 

\endroster 
\endproclaim 

\demo {Proof} Let $w\in\bold W(\ca H\cap\ca K)\cap\bold D(\ca H)$. Then 
$\Pi_{\ca M}(w)=\Pi_{\ca H}(w)\in\ca H$. As $\bold D(\ca H)\sub\bold D(\ca M)$ and $\ca K\leq\ca M$ we 
have also $\Pi(w)\in\ca K$. Thus $\ca H\cap\ca K\leq\ca H$ and $\ca H\cap\ca K\pce\ca K$. 
Transitivity of $\pce$ yields $\ca H\cap\ca K\pce\ca M$. Now suppose that $\ca K\norm\ca M$, and let 
$x\in\ca H\cap\ca K$ and $h\in\ca H$ such that $(h\i,x,h)\in\bold D(\ca H)$. Then 
$\Pi_{\ca H}(h\i,x,h)\in\ca H$. But also $(h\i,x,h)\in\bold D(\ca M)$, so that $x^h\in\ca K$. Thus 
$\ca H\cap\ca K\norm\ca H$.  
\qed 
\enddemo

\definition {Example 9.3} Let $(\ca H,\G)$ be a localizable pair in a locality $(\ca L,\D,S)$. Then 
$\ca H_\G$ is an im-partial subgroup of $\ca L$. 
\enddefinition

For the remainder of this section $(\ca L,\D,S)$ will be a regular locality on $\ca F$. Denote by $\ca L^s$ 
the expansion of $\ca L$ to a proper locality whose set of objects is $\ca F^s$, and by $\ca L^c$ the 
restriction of $\ca L^s$ to a proper locality on $\ca F$ whose set of objects is $\ca F^c$. Set 
$$ 
\ca F^*(\ca F)=\ca F_{S\cap F^*(\ca L)}(F^*(\ca L))\ \ \text{and} \ \ 
E(\ca F)=\ca F_{S\cap E(\ca L)}(E(\ca L)). 
$$ 
For any proper locality $\ca L'$ on $\ca F$ that can be obtained from $\ca L$ by a process of expansions 
and restrictions, define $F^*(\ca L')$ and $E(\ca L')$ by means of Theorem II.A2. 
\vskip .1in

For any subgroup $X\leq S$ write $\ca F_X$ for $N_{\ca F}(X)$. If $X$ is fully normalized in $\ca F$ 
then by II.6.3 there exists a proper locality $(\ca L^c_X,(\ca F_X)^c,N_S(X))$ on $\ca F_X$, and in fact the 
proof of II.6.3 involves the construction of $\ca L^c_X$ as in im-partial subgroup of the locality $\ca L^c$. 
Write $\ca L^s_X$ for the expansion of $\ca L^c_X$ to a proper locality on $\ca F_X$ whose set of objects is 
$(\ca F_X)^s$, and $\ca L_X$ for the restriction of $\ca L_X^s$ to a regular locality on $\ca F_X$.

\proclaim {Lemma 9.4} Let $X$ be fully normalized in $\ca F$, and let $P$ be a subgroup of 
$N_S(X)$ containing $X$. Suppose that $P\in(\ca F_X)^s$. Then $P\in\ca F^s$. 
\endproclaim 

\demo {Proof} By II.6.3 $\ca F_X$ is the fusion system of a proper locality. Subgroups of $N_S(P)$ which 
are fully normalized in $\ca F_X$ are then fully centralized in $\ca F_X$, by II.1.13. As 
$P\in(\ca F_X)^s$, there exists an $\ca F_X$-conjugate $Q$ of $P$ such that $Q$ is fully centralized in 
$\ca F_X$, and then $Q$ is fully centralized in $\ca F$ by II.1.16. Set $R=O_p(C_{\ca F}(Q))$. Then, 
since $Q\in(\ca F_X)^s$, II.6.4 yields $R\in C_{\ca F_X}(Q)^c$. As $C_{\ca F_X}(Q)=C_{\ca F}(Q)$, II.6.4 
yields also $Q\in\ca F^s$, and thus $P\in\ca F^s$. 
\qed 
\enddemo 

\proclaim {Corollary 9.5} Let $X$ be fully normalized in $\ca F$, and set $\G=(\ca F_X)^s$. Then 
$(N_{\ca L^s}(X),\G)$ is a localizable pair, and the locality $(N_{\ca L}(X)_\G,\G,N_S(X))$ given by 1.2 is 
isomorphic to $\ca L^s_X$. 
\endproclaim 

\demo {Proof} We note first of all that $\ca F_X$ is the fusion system of a proper locality by 6.3, 
hence $\G$ is $\ca F_X$-closed by 6.7. Secondly, $N_S(X)$ is a maximal $p$-subgroup of $N_{\ca L}(S)$ as 
a straightforward consequence of I.2.11(b) and of the hypothesis that $X$ is fully normalized in $\ca F$. 
Thirdly, for any $w\in\bold W(N_{\ca L}(X))$ with $N_{S_w}(X)\in\G$ we have $N_{S_w}(X)\in\ca F^s$ by 9.4. 
This verifies the three conditions for a localizable pair in definition 1.1, and so 1.2 applies and yields a 
locality $(N_{\ca L}(X)_\G,\G,N_S(X))$. By 1.4(b) this locality is in fact a proper locality on $\ca F_X$. 
Further, $(N_{\ca L}(X)_\G$ is an expansion of the locality $\ca L^c_X$ constructed in I.6.3, and is 
therefore isomorphic to $\ca L^s_X$ by Theorem II.A1. 
\qed 
\enddemo 

In what follows, we shall always identify $\ca L^s_X$ with the locality $N_{\ca L^s}(X)_\G$ in 9.5, and 
$\ca L_X$ with the restriction of $\ca L_X^s$ to $\d(\ca F_X)$. 

\proclaim {Corollary 9.6} We have $\ca L_X\pce\ca L^s_X\pce\ca L^s$. 
\endproclaim 

\demo {Proof} Immediate from 9.3. 
\qed 
\enddemo

\proclaim {Lemma 9.7} Let $X$ be fully normalized in $\ca F$. Then $X\in\ca F^s$ if and only if 
$E(\ca L_X)=\1$. 
\endproclaim 

\demo {Proof} Set $P=O_p(\ca L_X)$, and note that $P=O_p(\ca F_X)$ by II.2.3. Suppose first that 
$E(\ca L_X)=\1$. Then $F^*(\ca L_X)=P$ by an application of 8.5 to the regular locality $\ca L_X$. Thus $P$ 
is a large partial normal subgroup of 
$\ca L_X$, and so $C_{\ca L_X}(P)\leq P$. Then $P\in\ca F^c$ by II.6.2, and thus $X\in\ca F^s$. 
Conversely, suppose that $X\in\ca F^s$, so that $P\in\ca F^c$. Then $S\cap E(\ca L_X)=\1$, and then  
the regular locality $E(\ca L_X)$ is trivial. 
\qed 
\enddemo

\proclaim {Lemma 9.8} Let $\ca H$ be a product of components of $\ca L$, and let $\ca K$ be 
the product of those components of $\ca L$ which are not contained in $\ca H$. Assume that $\ca H\norm\ca L$, 
and let $D$ be a subgroup of $S\cap\ca H^\perp$ such that $D$ is fully normalized in $\ca F$ and such that 
$S\cap\ca K\leq D$. Then $S\cap\ca H\leq E(\ca L_D)$. 
\endproclaim 

\demo {Proof} Set $\G=\{P\in\D\mid D\norm P\}$, set $\ca E=\ca F_{S\cap\ca H}(\ca H)$, and set $V=S\cap\ca K$. 
Notice that $E(\ca L)=\ca H\ca K$ and that $[\ca H,\ca K]=\1$, by 8.5. Then 2.10 applies to the regular 
locality $E(\ca L)$, and 2.10(e) yields the following result. 
\roster 

\item "{(1)}" Let $U\leq S\cap\ca H$. Then $UV\in\D$ and $UD\in\G$ if and only if $U\in\ca E^s$. 

\endroster 
In particular, $\G$ is non-empty, and then 1.2 shows that $(N_{\ca L}(D),\G)$ is a localizable pair. 
Write $\ca L_{(D,\G)}$ for the locality $(N_{\ca L}(D)_\G,\G,N_S(D))$.

Let $w\in\bold W(\ca H)$ and set $P=N_{S_w}(D)$. Then $D\norm P$, and (1) shows that $P\in\G$ if and only 
if $S_w\cap\ca H\in\ca E^s$. Thus $P\in\G$ if and only if $w\in\bold D(\ca H)$, and this shows that $\ca H$ 
is a partial subgroup of $\ca L_{(D,\G)}$. As $\ca H\norm\ca L$ we obtain: 
\roster 

\item "{(2)}" $\ca H\norm\ca L_{(D,\G)}$. 

\endroster 

As $\ca H\norm\ca L$ we have $S\cap\ca H$ strongly closed in $\ca F$, and hence $S\cap\ca H$ is strongly 
closed in $\ca F_D$. The hypothesis of II.6.10 is then fulfilled, with the subcentric locality 
$(\ca L_D^s,(\ca F_D)^s,N_S(D))$ on $\ca F_D$ in the role of $\ca L$ and with $S\cap\ca H$ in the role of 
$T$. We may therefore conclude that for each $P\in(\ca F_D)^{cr}$ there exists an $\ca F_D$-conjugate $Q$ of 
$P$ with $Q\cap\ca H\in\ca E^c$. Then $Q\in\G$ by (1), and then $(\ca F_D)^{cr}\sub\G$ as $\G$ is 
$\ca F_D$-closed. Then 1.4(b) yields: 
\roster 

\item "{(3)}" $\ca L_{(D,\G)}$ is a proper locality on $\ca F_D$. 

\endroster 

We may now apply Theorem II.A to $\ca L_{(D,\G)}$. Let $\ca M^+$ be a large partial normal subgroup of 
$\ca L_D^s$ containing $O_p(\ca F_D)$, and set $\ca M=\ca M^+\cap\ca L_{(D,\G)}$. Then 
$S\cap\ca M=S\cap\ca M^+$, and $\ca H\cap\ca M\norm\ca L_{(D,\G)}$. Set $\ca H_0=C_{\ca H}(\ca H\cap\ca M)$. 
As $\ca H\cap\ca M\norm\ca H$ where $\ca H$ is a regular locality, 7.7(c) yields $\ca H_0\norm\ca H$. 
Then $\ca H_0=Z(\ca H)\ca H_1$ where $\ca H_1=O^p(\ca H_1)$ is a product of components of $\ca H$, by 8.6.  
As $\ca L_{(D,\G)}$ is a subset of $N_{\ca L}(S\cap\ca H)$, it follows from 7.7(e) that $\ca L_{(D,\G)}$ 
acts on $\ca H$ by conjugation. This action plainly preserves $\ca H_0$ and $\ca H_1$. Thus 
$\ca H_1\norm\ca L_{(D,\G)}$. We now claim (see 3.7): 
\roster 

\item "{(4)}" $\S_{S\cap\ca H_1}(\ca F_D)\sub\G$. 

\endroster 
For this it suffices to show that $DXY\in\G$ for $X\in\ca F_{S\cap\ca H_1}(\ca H_1)^{cr}$ and  
$Y\in C_{\ca F_D}(S\cap\ca H_1)^{cr}$. For any such $X$ and $Y$ we have $XY\in\ca E^{cr}$ and 
$VXY\in E(\ca F)^{cr}$ by 2.10(d), so $VXY\in\D$. Thus $DXY\in\G$, and (4) holds. 

We have $[\ca H_1,S\cap\ca M]\leq\ca H_1\cap\ca M\leq Z(\ca H_1)$, and then $[\ca H_1,S\cap\ca M]=\1$ 
since $\ca H_1=O^p(\ca H_1)$. Let $\ca H_1^+$ be the partial normal subgroup of $\ca L^s_D$ whose 
intersection with $\ca L_{(D,\G)}$ is $\ca H_1$. Then $\ca H_1^+$ is the partial subgroup $\<\ca H_1\>$ 
of $\ca L^s_D$ generated by $\ca H_1$, by (4) and 3.14(b), and it follows that $[\ca H_1^+,S\cap\ca M]=\1$. 
As $\ca H_1^+=O^p(\ca H_1^+)$ by II.7.3, it follows from definition 5.5 that $\ca H_1^+\leq(\ca M^s)^\perp$ 
(where the ``$\perp$" operation is taken in the locality $\ca L^s_D$). As $\ca M^s$ is large in $\ca L^s_D$ 
we have $S\cap(\ca M^s)^\perp=Z(\ca M^s)$ by 6.4, and thus $S\cap\ca H_1\leq\ca M$. 
As $\ca H_1\cap\ca M\leq Z(\ca H_1)$ we conclude that $\ca H_1=\1$ and that $\ca H\leq\ca M$. Thus 
$S\cap\ca H\leq F^*(\ca L_D^s)$. As $S\cap F^*(\ca L_D^s)=S\cap F^*(\ca L_D)$ by Theorem II.A, the proof 
is complete. 
\qed 
\enddemo 

We are now in position to prove a version for regular localities of the ``$L$-balance" Theorem from 
finite group theory. For a discussion of this result for finite groups see [GLS], and for a fusion-theoretic 
version see [Asch].

\proclaim {Theorem 9.9} Let $\ca L$ be a regular locality $X\leq S$ be fully normalized in $\ca F$. Then 
$S\cap E(\ca L_X)\leq E(\ca L)$. Moreover, if $X\leq F^*(\ca L)$ then $\ca L_X\pce\ca L$ and 
$E(\ca L_X)\pce E(\ca L)$. 
\endproclaim 

\demo {Proof} Assume false, and let $X$ be a counterexample of maximal order. set $R=S\cap E(\ca L))$, and 
let $(\ca K_1,\cdots,\ca K_m)$ be a listing (in arbitrary order) of the components of $\ca L_X$. As 
$[\ca K_i,\ca K_j]=\1$ for all $i$ and $j$ with $i\neq j$, I.5.1 shows that $S\cap E(\ca L_X)$ is the 
product in any order of the intersections $S\cap\ca K_i$. Let $\ca H$ be the product of all of the 
components $\ca K_i$ such that $S\cap\ca K_i\nleq E(\ca L)$. As $\ca L_X\pce\ca L^s$ we have 
$\ca H\cap E(\ca L^s)\norm\ca H$ by 9.2(b), and then $\ca H\cap E(\ca L^s)\leq Z(\ca H)$ by 8.6. 

Set $A=N_R(X)$. Then $[S\cap\ca H,A]\leq Z(\ca H)$, so $A$ normalizes every $P$ in 
$\ca F_{S\cap\ca H}(\ca H)^{cr}$. For each such $P$, 9.2(b) yields 
$$ 
[N_{\ca H}(P),A]\leq N_{\ca H}(P)\cap E(\ca L)\leq Z(\ca H).     
$$ 
Thus $N_{\ca H}(P)$ centralizes the chain $Z(H)A\geq Z(A)\geq 1$, and then 
$[O^p(N_{\ca H}(P)),A]=\1$ by coprime action (II.2.5). By 5.2 $O^p(\ca H)$ is generated by 
the set of all $O^p(N_{\ca H}(P))$ with $P$ as above, and then $[\ca H,A]=\1$ as $\ca H=O^p(\ca H)$. Thus 
$[\ca H,AX]=\1$. 

Let $D\in(AX)^{\ca F}$ be fully normalized in $\ca F$. As $\ca F$ is 
inductive there exists an $\ca F$-homomorphism $\phi:N_S(AX)\to N_S(D)$ such that $(AX)\phi=D$. Set 
$Y=X\phi$ and $B=A\phi$. Then $B\leq N_R(Y)$, and equality then holds since $X$ is fully 
normalized in $\ca F$, and since $R$ is strongly closed in $\ca F$. As $N_S(X)\leq N_S(AX)$, $\phi$ 
maps $N_S(X)$ into $N_S(Y)$, and so $Y$ is fully normalized in $\ca F$. We observe that the 
$\ca F$-isomorphism $N_S(X)\to N_S(Y)$ induced by $\phi$ induces also an isomorphism 
$\ca F_X\to\ca F_Y$ of fusion systems, by II.1.5. Since $S\cap E(\ca L_X)$ and $S\cap E(\ca L_Y)$ 
depend only on $\ca F_X$ and $\ca F_Y$, we may now assume that $X=Y$ and that $D$ is fully normalized 
in $\ca F$. Moreover, from the preceding paragraph we have: 
\roster 

\item "{(1)}" $\ca H\leq C_{E(\ca L_X}(D)$. 

\endroster 
Note that $R\in\ca F^s$ by II.6.9. If $R\leq X$ then also $X\in\ca F^s$, and then 9.7 yields a contradiction 
to the non-triviality of $\ca H$. Thus $R\nleq X$, and so $X$ is a proper subgroup of $D$. The maximality of 
$|X|$ then implies that $D$ is not a counterexample to the lemma, and so: 
\roster 

\item "{(2)}" $S\cap E(\ca L_D)\leq E(\ca L)$. 

\endroster 

Set $S_{D,X}=N_S(D)\cap N_S(X)$, and let $\ca F_{D,X}$ be the fusion system $N_{N_{\ca F}(D)}(X)$ on 
$S_{D,X}$. It is a straightforward consequence of the definition (following II.1.4) of normalizers  
in fusion systems that $\ca F_{D,X}=N_{N_{\ca F}(X)}(D)$. Set $\G_{D,X}=(\ca F_{D,X})^s$, and write 
$(\ca L^s_{D,X},\G_{D,X},S_{D,X})$ for the proper locality on $\ca F_{D,X}$ given by two applications 
of 9.5. By 9.4, a word $w\in\bold W(N_{\ca L}(D)\cap N_{\ca L}(X))$ is in $\bold D(\ca L^s_{D,X})$ if and 
only if $S_w\cap S_{D,X}\in\ca F^s$. Thus: 
\roster 

\item "{(3)}" $\ca L^s_{D,X}\pce\ca L^s_D\pce\ca L^s$. 

\endroster 
We may now appeal to lemma 9.8 with $\ca L_X$ in the role of $\ca L$, in order to obtain 
$S\cap\ca H\leq E(\ca L^s_{D,X})$. We may write also: 
\roster 

\item "{(4)}" $S\cap\ca H\leq E(\ca L_{D,X})$, 

\endroster 
where $\ca L_{D,X}$ is the regular locality obtained by restriction from $\ca L^s_{D,X}$. 

Suppose that $D\leq O_p(\ca L)$. Then $X\leq O_p(\ca L)$, and so $X\leq E(\ca L)^\perp$. This yields $A=R$,  
so that $R\leq D\leq O_p(\ca L)$. Then $E(\ca L)=\1$, and $F^*(\ca L)=O_p(\ca L)$. Thus $O_p(\ca L)\in\D$, 
and $\ca L$ is a group of characteristic $p$. Then $\ca L_X=N_{\ca L}(X)$ is a group of characteristic $p$ by 
II.2.7(b), and $E(\ca L_X)=\1$. As there is nothing to prove in this case, we may assume that 
$D\nleq O_p(\ca L)$. 

By induction on $|\ca L|$ we may then assume that the Theorem holds for $\ca L_D$ in 
the role of $\ca L$. Thus:  
\roster 

\item "{(5)}" $S\cap E(\ca L_{D,X})\leq S\cap E(\ca L_D)$.  

\endroster 
Suppose next that $X=D$. Then $R\leq X$, so $O_p(\ca L)X\in\D\sub\ca F^s$, and then $X\in\ca F^s$ by II.6.9. 
Then $E(\ca L_X)=\1$ by 9.7, and there is again nothing to prove. Thus we may assume that $X$ is a proper 
subgroup of $D$, and the maximal choice of $X$ then yields $S\cap E(\ca L_D)\leq E(\ca L)$. This result, in 
combination with (4) and (5), now yields $S\cap\ca H\leq E(\ca L)$. By definition, $\ca H$ is the product 
of those components $\ca K$ of $\ca L$ such that $S\cap\ca K\nleq E(\ca L)$, we conclude that 
$\ca H=\1$, and that $S\cap E(\ca L_X)\leq E(\ca L)$. This completes the first part of the proof. 

Suppose now that $X\leq F^*(\ca L)$. Let $w\in\bold D(\ca L_X)$, and set $P=N_{S_w}(X)\cap F^*(\ca L)$. Then 
$X\leq P$, and $S_w\cap E(\ca L_X)\leq P$ by what has just been proved. Here 
$S_w\cap F^*(\ca L_X)\in\ca F_X^s$ by the definition of $\bold D(\ca L_X)$, and then 
$S_w\cap E(\ca L_X)\in\ca F_X^s$ by II.6.9. Thus $X\leq P\in\ca F_X^s$, and then $P\in\ca F^s$ by 9.4. 
Thus $P\in\d(\ca F)$. Recall that $\ca L_X\pce\ca L^s$ by definition. In the case that $w=(g)$ for some 
$g\in\ca L_X$ we now conclude that $g\in\ca L$. Thus $\ca L_X\sub\ca L$, and 
$\bold D(\ca L_X)\sub\bold D(\ca L)$. This yields $\ca L_X\pce\ca L$. 

We have already seen that $S\cap E(\ca L_X)\leq E(\ca L)$. As $E(\ca L_X)\pce\ca L$ we have 
$E(\ca L_X)\cap E(\ca L)\norm E(\ca L_X)$ by 9.2(b). Each component $\ca K$ of $\ca L_X$ has the 
property that $\ca K/Z(\ca K)$ is simple, so $\ca K$ has no proper partial normal subgroups containing 
$S\cap\ca K$. This shows that $E(\ca L_X)\sub E(\ca L)$, and then $E(\ca L_X)\pce E(\ca L)$ by 9.2(a). 
\qed 
\enddemo

\vskip .2in 
\noindent 
{\bf Section 10: Theorems C, D, and E} 
\vskip .1in 

Theorem C is given by Theorems 5.7 and 7.7 (which in fact contain a great deal more information), and 
Theorem D is Theorem 8.5. 

Let $(\ca L,\D,S)$ be a regular locality on $\ca F$. Point (a) of Theorem E is that if $\ca N\norm\ca L$ 
is a partial normal subgroup of $\ca L$ then $\ca N$ is a regular locality; which is already part of 
Theorem C. Point (b) of Theorem E, concerning the regular locality $\ca L_X$ on $N_{\ca F}(X)$, under the 
assumption that $X$ is fully normalized in $\ca F$ and that $X\leq F^*(\ca L)$, is given by Theorem 9.9. 
Point (c) of Theorem E is Proposition 7.2.

\enddocument